\newif\ifpreprint
\preprinttrue  

\ifpreprint
\documentclass[11pt]{article}

\usepackage{fullpage}
\usepackage{natbib}
\else 
\documentclass[twoside,11pt]{article}
\usepackage[abbrvbib,preprint]{jmlr2e}

\newtheorem{assumption}{Assumption}
\newtheorem{problem}{Problem}
\newtheorem{exercise}{Exercise}

\fi

\usepackage{enumitem}
\setlist[itemize]{itemsep=2pt, topsep=4pt, parsep=0pt}
\ifpreprint 
\usepackage{amsmath,amssymb,amsthm}
\usepackage{graphicx}

\usepackage[
colorlinks=true,
citecolor=cyan,
linkcolor=cyan,
urlcolor=cyan]{hyperref}
\fi
\usepackage{array}
\usepackage{mathtools}
\usepackage{appendix}

\usepackage{accents}
\usepackage{multirow}
\usepackage{algorithm}\usepackage[noend]{algpseudocode} \usepackage{subcaption}
\usepackage[font=small,labelfont=bf]{caption}
\usepackage[acronyms]{glossaries}   \glsdisablehyper  

\usepackage{csvsimple}	\usepackage{booktabs}   \usepackage{adjustbox}  \usepackage{placeins}  

\newcommand{\uncset}{{\mathcal{U}}}

\newcommand{\eg}{{\it e.g.}}
\newcommand{\ie}{{\it i.e.}}

\newcommand{\ones}{\mathbf 1}
\newcommand{\reals}{{\mathbf R} }

    \newcommand{\prob}{{\mathbf P}}

\newcommand{\Tr}{\mathop{\bf tr}}
\newcommand{\diag}{\mathop{\bf diag}}

\newcommand{\Expect}{{\mbox{\bf E}}}
\newcommand{\Co}{{\mathop {\bf Co}}} 
 
\newcommand{\dom}{\mathop{\bf dom}} 
\newcommand{\dif}{\text{d}}

\newcommand{\define}{\coloneqq}

\newcommand{\indicator}{\mathcal{I}}
\newcommand{\CVaR}{\mathop{\bf CVaR}}
\newcommand{\VaR}{\mathop{\bf VaR}}

\newcommand{\supp}{S}

\newacronym{LO}{LO}{linear optimization}
\newacronym{QO}{QO}{quadratic optimization}
\newacronym{MIQO}{MIQO}{mixed-integer quadratic optimization}
\newacronym{MIO}{MIO}{mixed-integer optimization}
\newacronym{MILO}{MILO}{mixed-integer linear optimization}
\newacronym{MINLO}{MINLO}{mixed-integer nonlinear optimization}
\newacronym{sBB}{sBB}{spatial branch and bound}
\newacronym{NLO}{NLO}{nonlinear optimization}
\newacronym{PWA}{PWA}{piecewise affine}
\newacronym{SVM}{SVM}{support vector machines}
\newacronym{ReLU}{ReLU}{rectified linear unit}
\newacronym{CPU}{CPU}{central processing unit}
\newacronym{GPU}{GPU}{graphics processing unit}
\newacronym{MPC}{MPC}{model predictive control}
\newacronym{ADMM}{ADMM}{alternating direction method of multipliers}
\newacronym{ADP}{ADP}{approximate dynamic programming}
\newacronym{FPGA}{FPGA}{field-programmable gate array}
\newacronym{DRO}{DRO}{distributionally robust optimization}
\newacronym{RO}{RO}{robust optimization}
\newacronym{SO}{SO}{stochastic optimization}
\newacronym{MRO}{MRO}{mean robust optimization}
\newacronym{LROPT}{LROPT}{learning for robust optimization}
\newacronym{LRO}{LRO}{learning for robust optimization}
\newacronym{DRP}{DRP}{disciplined robust programming}
\newacronym{DCP}{DCP}{disciplined convex programming}
\newacronym{DPP}{DPP}{disciplined parameterized programming}
\newacronym{DSP}{DSP}{disciplined saddle programming}
\newacronym{SG}{SG}{scenario generation}
\newacronym{SGDA}{SGDA}{stochastic gradient descent ascent}
\newacronym{ETO}{ETO}{estimate-then-optimize} 

\newcommand{\githubrepo}{\url{https://github.com/stellatogrp/cvxro_experiments}}

\ifpreprint 
\newtheorem{theorem}{Theorem}[section]
\newtheorem{lemma}{Lemma}[section]

\newtheorem{proposition}{Proposition}[section]
\newtheorem{assumption}{Assumption}[section]
\newtheorem{definition}{Definition}[section]
\newtheorem{example}{Example}[section]

\newtheorem{remark}{Remark}[section]

\fi 

\usepackage{xcolor}

\ifpreprint
\newcommand{\RC}[1]{{#1}}

\else
\newcommand{\RC}[1]{{\color{blue}#1}}

\fi

\graphicspath{{./figs/}}

\newcommand{\appxsets}{\ifpreprint Appendix~\ref{appx:gen_sets}\else Online Appendix~1\fi}
\newcommand{\appxdirect}{\ifpreprint Appendix~\ref{appx:bilinear}\else Online Appendix~7\fi}
\newcommand{\appxablation}{\ifpreprint Appendix~\ref{appx:ablation}\else Online Appendix~8\fi}
\newcommand{\appxdataport}{\ifpreprint Appendix~\ref{appx:data_port}\else Online Appendix~8.1.1\fi}
\newcommand{\appxdatainv}{\ifpreprint Appendix~\ref{appx:data_inv}\else Online Appendix~8.1.2\fi}
\newcommand{\appxhyper}{\ifpreprint Appendix~\ref{appx:hyper_ablation}\else Online Appendix~8.2\fi}
\newcommand{\appxdfl}{\ifpreprint Appendix~\ref{appx:dfl}\else Online Appendix~8.3\fi}

\begin{document}

\ifpreprint
\title{Learning Decision-Focused Uncertainty Sets \\
in Robust Optimization}
\author{Irina Wang, Bart Van Parys, and Bartolomeo Stellato}
\maketitle
\begin{abstract}

We propose a data-driven technique to automatically learn tractable contextual uncertainty sets in robust optimization, for a family of problems with both objective and constraint uncertainty. 
This decision-focused method reshapes the uncertainty sets by optimizing for both worst-case and average-case performance while guaranteeing constraint satisfaction \RC{with high probability} via a conditional value-at-risk constraint.
We solve the resultant constrained learning problem using a stochastic augmented Lagrangian method that relies on differentiating the solutions of the robust optimization problems with respect to the parameters of the uncertainty set.
\RC{Due to the nonsmooth and nonconvex nature of the augmented Lagrangian function, we apply the nonsmooth conservative implicit function theorem to establish that, under regularity and boundedness assumptions, the accumulation points of the resulting algorithm are feasible solutions of the constrained learning problem when the penalty parameters remain bounded, and stationary points of the constraint violation otherwise.}
\RC{Using a covering argument with Hoeffding's inequality, we show both marginal and conditional finite-sample probabilistic guarantees of constraint satisfaction for any approximately feasible solution.}
\RC{We also characterize when a decision-focused set is worth learning.}
\RC{In the proposed settings, numerical experiments show that our method outperforms traditional approaches in data-driven robust optimization in terms of out-of-sample performance.}

 \end{abstract}

\fi

\section{Introduction}
\label{sec:intro_p}
\setlength{\parskip}{0pt}
Over the past years, \gls{RO} has become a widely adopted efficient tool for decision-making under uncertainty.
A robust approach first defines an uncertainty set where the uncertain vector lives, and subsequently hedges its decision against any worst-case realization in this set~\citep{robustconvexopt,ben-tal_robust_2009,bertsimassurvey,robustadaptopt}.
To attain proper solution quality and practical robustness, the choice of uncertainty set is paramount.

Traditionally, uncertainty sets are shaped and sized to guarantee they contain an appropriately large fraction of the probability mass of the unknown distribution generating the samples (\ie, high coverage).
The main appeal of this approach is that the uncertainty set is independent of the subsequent optimization problem considered, and guarantees the worst-case performance regardless of the problem.
In this vein, many approaches to designing uncertainty sets assume structural information of the underlying distributions, and rely on these a priori assumptions to build guarantees of constraint satisfaction~\citep{ben-tal_robust_2000,bertsimassurvey,bandi_tractable_2012}.
\RC{In particular, many works bound the probability of constraint violation either a priori, from the size and shape of the set alone, or a posteriori, from the solution obtained~\citep{guarantees_2012, guarantees1}, but in neither case are the set parameters chosen for the quality of the resulting decisions.} 
\RC{Even in recent data-driven \gls{RO}, whether by approximating a high-probability region through quantile estimation~\citep{hong2022ms}, by learning a contextual uncertainty set that contains the realized uncertainty with high probability~\citep{knn}, or by requiring a \gls{DRO} ambiguity set to contain the data-generating distribution with high confidence~\citep{mohajerin_esfahani_data-driven_2018}, significant attention has been given to the probability distribution generating the data and less so to the quality of the resulting decisions.}

\RC{This decoupling of the uncertainty set construction from the subsequent optimization problem may lead to overly conservative decisions~\citep{reduce_conserve,gao2016distributionally}.
An established remedy is to take the optimization problem into account when designing the set.
As noted by~\cite{ben-tal_robust_2009}, high coverage is only a sufficient, and often restrictive, condition for a high probability of constraint satisfaction, and a long line of work has constructed uncertainty sets that target the probability of constraint violation directly and may have low coverage~\citep{bertsimas_price_2004,bertsimas_probabilistic_2019,scenario}.
In \gls{DRO}, likewise, a Wasserstein radius calibrated for coverage must grow with the dimension of the uncertainty~\citep{wass_rate}, whereas exploiting Lipschitz properties of the downstream objective removes this dependence~\citep{DBLP:gao2020}.}
Still, while these works take into account the problem constraints, and reduce conservatism compared to the decoupled approaches, they either do not optimize over uncertainty-set parameters to find the best one in terms of the performance of the resulting decisions,
or only do so in a limited capacity, that is, with very few tunable parameters.
\RC{These methods are generally non-contextual, and many consider constraints individually, which makes contextual problems with multiple uncertain constraints difficult to tune.}

\RC{In a separate vein, researchers have explored the idea of {\it decision-focused learning}, where a model predicting the unknown parameters of an optimization problem is trained not to minimize prediction error, but to optimize the quality of the decisions its predictions induce~\citep{predict_opt,end-to-end}.
However, these point-prediction formulations pass a forecast to a deterministic downstream problem, which does not enforce a probability guarantee on the uncertain constraints.
Our work brings this decision-focused perspective to the construction of uncertainty sets themselves.}

Lastly, the performance of the \gls{RO} problem may be judged on criteria that differ across applications, which may be optimized by different parameterizations of the uncertainty sets.
Traditionally, the goal of \gls{RO} uncertainty sets is to provide solutions that perform well with respect to the true worst-case realizations.
However, in practical applications, the decision-maker would prefer solutions that are suitable for typical cases as well.
Recently, a line of work on Pareto efficiency in \gls{RO}~\citep{pareto,pareto1,pareto3} \RC{identifies worst-case optimal solutions that, in addition, cannot be improved in any scenario of the uncertainty set without being worsened in another.}
\RC{Extending this idea of looking beyond the worst case, we construct uncertainty sets in a decision-focused manner,}
such that the resultant robust solutions perform well not only in terms of constraint satisfaction, but also in {\it both} worst-case and average-case performance.

\subsection{Our contributions}\label{sub:our_contributions}
In this work, we propose an end-to-end approach to calibrate uncertainty sets for a family of contextual optimization problems, \RC{which can have both uncertain objectives and constraints.}
This is a common scenario in many applications, including inventory management, where we solve similar optimization problems with varying prices and factors driving demand (\ie, context parameters), while satisfying the uncertain demand with high probability.
\RC{We construct the uncertainty set, which maps from contexts to context-based regions of uncertainty, by learning its contextual shape parameters using a custom performance metric.}
By directly taking into account the quality of the \gls{RO} solutions through the performance metric,
the learned set is {\it decision-focused}. 
Our specific contributions are:
\begin{itemize}
\item 
\RC{We learn the uncertainty set parameters by minimizing, over the family of parametric \gls{RO} problems, the expected value of a performance metric that combines worst-case and average-case performance, subject to a conditional value-at-risk ($\CVaR$) constraint that enforces the desired probability of constraint satisfaction.}
\item \RC{We present conditions under which we do and do not expect decision-focused learning to outperform traditional two-step procedures.}\item We present the learning problem as a \RC{nonsmooth} constrained bi-level problem, where the inner level is represented by the family of parametric \gls{RO} problems. 
\item \RC{We develop a stochastic augmented Lagrangian algorithm to solve the nonsmooth bi-level problem, which requires differentiating through the inner level solutions as well as the $\CVaR$ constraint.
Due to the nonconvex and nonsmooth nature of the problem, we show {\it path-differentiability} of the Lagrangian, and apply nonsmooth implicit differentiation.}
\RC{We show that, under regularity, boundedness, and step-size conditions, the accumulation points of our algorithm are feasible solutions of the constrained learning problem when the penalty parameters remain bounded, and are otherwise stationary points of the constraint violation.}\item \RC{We give both marginal and conditional finite-sample probabilistic guarantees of constraint satisfaction for the resulting solution.}
\item We benchmark our method on examples in portfolio optimization and multi-stage inventory management, outperforming data-driven robust and distributionally robust approaches in terms of out-of-sample performance.
The code is available at \githubrepo.
\end{itemize}
\subsection{Related work}\label{sub:related_work}
\paragraph{Data-driven robust optimization.}
\RC{In data-driven robust optimization, uncertainty sets are traditionally built from sample-based approximations of the unknown data-generating distribution.}
\RC{Approaches such as quantile estimation~\citep{hong2022ms}, bounding the type-$\infty$ Wasserstein distance~\citep{wassinf}, constructing intersections of multiple single-measure-based ambiguity sets~\citep{tanoumanddata}, performing uncertainty quantification using side-information~\citep{persak2023contextual}, and applying deep learning clustering techniques~\citep{GOERIGK2023106087} have focused on modeling the distribution.} 
\RC{In fact, the majority of these methods rely on a two-step procedure, which separates the construction of the uncertainty set from the resulting robust optimization problem, and often enforces a high coverage requirement on the selected uncertainty sets.}
\RC{In a finite-dimensional Bayesian setting with constraints concave in the distribution parameter, however,~\cite{gupta} shows that \gls{DRO} ambiguity sets with asymptotic guarantees can be far smaller than the confidence regions traditionally used.}
To reduce conservatism,~\cite{bertsimas_data-driven_2018} use hypothesis testing to pair a priori assumptions on the distribution with different statistical tests, and obtain various uncertain-constraint-dependent uncertainty sets with different shapes, computational properties, and modeling power.
These sets are not coverage-based, and offer improved performance for the targeted constraints.
Similarly,~\cite{vanparys2021data} use a \gls{DRO} approach with constraints on out-of-sample disappointment to find least conservative cost-estimators and optimizers for stochastic optimization problems.
We go
one step further, and construct parametric uncertainty sets to optimize {\it target performance metrics} for particular {\it families of contextual \gls{RO} problems}. 
To this end, we turn to ideas in differentiable optimization.
Recent developments in differentiable convex optimization have enabled a popular {\it end-to-end} approach of training a predictive model to minimize loss on a downstream optimization task, \RC{also coined {\it decision-focused learning}.
See the following paragraph for details.}
Our approach adapts this idea, by training the shape of the uncertainty set on its performance on the downstream parametric \gls{RO} problems. 
In a work with a similar goal, \cite{schuurmans2023distributionally} construct \gls{DRO} ambiguity sets which constrain only the worst-case distribution along the direction in which the expected cost of an approximate solution increases most rapidly.
Our work, however, considers a broader class of optimization problems and uncertainty sets, and \RC{handles several uncertain constraints jointly through a $\CVaR$ constraint on their maximum}.

\paragraph{Decision-focused learning.}
\RC{Decision-focused learning trains a predictive model end-to-end so that predictions are optimized for the quality of the resulting downstream decision, rather than for a prediction-error metric~\citep{amos1,end1,end2,end3,end-to-end,predict_opt}.
We point to the surveys by~\cite{end-survey} and~\cite{dfl} for a broad overview.
In the typical point-prediction subclass, the downstream problem is deterministic and takes a single scenario even when the underlying problem is stochastic, while other works learn distributional models~\citep{amos1} or decision proxies~\citep{schutte2025proxies}.
\cite{homemdemello2026} and~\cite{schutte2025proxies} characterize when such single-scenario approximations suffice for optimal decisions, in the average case and within specific problem classes.
\cite{schutte2024robust} propose robust loss functions for decision-focused learning, but do not test worst-case objective values or control constraint violation.}
\RC{A separate line of work asks when decision-focused training is worth it at all, by comparing it with the two-step approach.
\cite{hu2022fast} show that, in contextual linear optimization with a well-specified model, the two-step approach attains faster regret rates than integrated learning.
\cite{elmachtoub2025eto} show, through an asymptotic stochastic dominance argument, that the two-step approach dominates integrated estimation-optimization when the model class contains the truth, and that the order reverses under misspecification.
In these comparisons, decision-focused training pays off under misspecification. 
}
\RC{Our setting differs from this literature on multiple counts.
Our learned object is not a point prediction or a predictive distribution, but an uncertainty set, and the downstream problem is a robust optimization problem whose tail performance we control directly through the $\CVaR$ constraint, rather than through a single scenario.
This holds even when only the center of the set is trained and its shape is held at a first-step fit, as the center enters the downstream problem with a calibrated radius rather than as the scenario itself.
Furthermore, existing comparisons with the two-step approach consider problems with uncertainty only in the objective, and take into account only the expected value. 
When there are multiple uncertain constraints and the tail performance matters, we show that even the exact conditional moments can yield a suboptimal set, so the two-step set can be misspecified in its decision layer no matter how well its moments are estimated.}

\paragraph{End-to-end contextual optimization.}
\cite{bertsimas2020predictive} introduce a framework for prescriptive analytics that uses contextual information to directly prescribe optimal decisions, bypassing the intermediate step of estimating unknown parameters.
\RC{The decision-focused methods discussed above are typically deployed in this contextual setting, mapping context features to the cost or constraint parameters of a deterministic downstream problem.
As pointed out in the recent survey by~\cite{contex}, however, work in the end-to-end framework seldom addresses uncertainty in the downstream optimization problem, especially in the constraints.
Our approach aims to fill this gap, and we therefore focus here on the end-to-end works whose downstream problem is itself robust.}
\RC{Close to our setting,~\cite{cond-cov} propose an end-to-end training algorithm to produce contextual uncertainty sets that target a certain conditional coverage, for robust problems with an uncertain objective.
They construct penalty functions for the achieved objective value and the conditional coverage of the uncertainty set, and train the parameters of ellipsoidal uncertainty sets via stochastic gradient descent (SGD).
Decision protection then follows if the required conditional coverage is achieved, which they assess only empirically, and no convergence guarantees are given for the algorithm.
\cite{yeh} also treat an uncertain objective, with sets calibrated by conformal prediction, which carries formal marginal coverage guarantees. 
We, however, allow for not only an uncertain objective, but also for any number of uncertain constraints, and formulate our problem as a {\it constrained} optimization problem.
We do not enforce coverage constraints, which may introduce unnecessary conservatism, and instead constrain the $\CVaR$ of the constraint violations.
We also optimize over the average-case performance in addition to the worst-case performance, which results in less conservative decisions when the uncertainty realization is not the worst case. 
In addition, we characterize the accumulation points of our algorithm, and prove that the resulting decisions satisfy the uncertain constraints with high probability.}

\paragraph{Constrained learning.}
\RC{The idea of formulating learning as a constrained empirical-risk minimization problem, and updating dual multipliers via a Lagrangian-dual penalty to enforce constraint satisfaction, has been explored outside the uncertainty-set-learning setting.
\cite{fioretto2020lagrangian} train deep networks subject to problem-specific constraints by iteratively updating Lagrangian multipliers on the constraint violation, and~\cite{kotary2022fast} apply a similar Lagrangian-dual deep learning approach to obtain fast, constraint-satisfying approximations for job shop scheduling.
\cite{kotary2024deepal} extend this idea to the {\it learning to optimize} setting, training a model to predict dual solution estimates directly, and use augmented Lagrangian techniques to correct the poor convergence of a naive dual-ascent training scheme.
\cite{boero2025alcole} provide a supporting theoretical foundation for this line of work, showing that augmented Lagrangian-based constrained learning achieves strong duality and primal convergence to feasible, optimal solutions without the randomization or averaging required by earlier Lagrangian-dual training schemes.
These works constrain the predictor itself or learn to approximate the solution of an optimization problem, and do not differentiate through a per-instance optimization layer.
Our focus, however, is on developing contextual uncertainty sets for constrained \gls{RO} problems.
This yields a nonsmooth bi-level formulation in which the constrained \gls{RO} problem is solved at the inner level and differentiated through, alongside a $\CVaR$ constraint at the outer level.
We develop a stochastic augmented Lagrangian algorithm specific to our setting, with guarantees on its accumulation points and finite-sample probabilistic guarantees.
}

\paragraph{Pareto efficiency in \gls{RO}.}
\RC{The concept was introduced} by~\cite{pareto} for linear \gls{RO} problems, where they identify solutions that yield worst-case optimal performance and which are not dominated by any other such solutions in non-worst-case scenarios.
The concept has since been extended to general linear \gls{RO} problems in Euclidean spaces~\citep{pareto3}, and to linear adaptive robust optimization problems~\citep{pareto1}. 
The focus of this line of work is to find Pareto-dominating solutions among a set of \gls{RO} solutions for a fixed uncertainty set.
Tangentially, while not aimed directly at finding Pareto optimal solutions, the robust satisficing framework of~\cite{long2023robust} is also demonstrated to improve the Pareto efficiency compared to classic \gls{RO}, by searching among a {\it different set of solutions}. 
In our work, while also not directly searching for Pareto optimal solutions, we optimize a decision-focused contextual uncertainty set so that the unique \gls{RO} solution implied by the uncertainty set may perform better both in the worst case and on average.

\paragraph{Differentiable convex optimization.}
There has been extensive work on embedding differentiable optimization problems as layers within deep learning architectures~\citep{amos2021optnet}.
\cite{cvxpylayers2019} propose an approach to differentiating through disciplined convex programs, and implement their methodology in CVXPY.
Specifically, they implement differentiable layers for disciplined convex programs in PyTorch and TensorFlow 2.0, in the package cvxpylayers.
They convert the convex programs into conic programs, then implicitly differentiate the residual map of the homogeneous self-dual embedding associated with the conic program.
These developments have fueled the recent advances in end-to-end approaches, and are also the tool of choice in this work for differentiating through the inner level \gls{RO} solutions.

\paragraph{Nonsmooth implicit differentiation.}
Traditional techniques in differentiable convex optimization rely on the implicit function theorem, which requires \RC{a continuously differentiable residual with an invertible partial Jacobian}. 
Recently, a nonsmooth variant has been established, using conservative Jacobians, provided that a nonsmooth form of the classical invertibility condition is fulfilled~\citep{path-diff}.
As a key ingredient, definable and locally Lipschitz continuous functions are shown to admit conservative Jacobians, which follow the chain rule, and are thus compatible with backpropagation techniques in stochastic gradient descent~\citep{Davis,path-diff}. 
These techniques can be used to prove convergence to stationary points for nonsmooth problems.

\subsection{Layout of the paper}In Section~\ref{sec:intro}, we introduce the notion of the decision-focused set, \RC{explain its advantages,} and give a motivating example.
In Section~\ref{sec:learn}, we describe the data-driven algorithm and prove its convergence, and in Section~\ref{sec:guarantees}, provide probabilistic guarantees.
In Section~\ref{sec:hyper}, we discuss hyperparameter tuning, and in Section~\ref{sec:examples}, present various numerical examples.
For completeness, we give example convex reformulations of robust optimization problems for common uncertainty sets, with our shaped parameters, in the appendices.

\setlength{\parskip}{0pt plus 1pt}
\section{The optimization problem}
\label{sec:intro}
Consider a contextual stochastic optimization setting with context parameter $x \in \reals^p$ and uncertain parameter $u \in \reals^m$. 
The context parameter $x$ represents information available when we make our decisions (\eg, market factors for portfolio optimization, the supplier's price for inventory problems), while the uncertain parameter $u$ is not known (\eg, future returns in portfolio optimization), but may be affected by the context parameter.
We set $\prob_x$ as the distribution of $x$ with domain $\mathcal{D}_x \subseteq \reals^p$, $\prob_u$ as the distribution of $u$ with domain $\mathcal{D}_u \subseteq \reals^m$, and $\prob_{(u, x)}$ as their joint distribution, all of which are unknown.
Importantly, since $x$ may affect $u$, the two can be correlated.
We therefore denote by $\prob_{(u \mid x = \bar{x})}$ the conditional distribution of $u$ for a context $\bar{x} \in \mathcal{D}_x$.

\subsection{The conditional formulation}
\label{sec:cond}
For any context $\bar{x}$, the decision-maker, who has some uncertain objective function $f(z,u,\bar{x})$, would like to ensure that their decision $z(\bar{x}) \in \reals^n$ performs well in both the average case and the worst case, \RC{where we use the following {\it conditional value-at-risk} ($\CVaR$) formulation as a proxy of the worst-case value.
Specifically, for $\eta \in (0,1)$, the $\CVaR$ at level $\eta$ is the average of the worst $\eta$ fraction of outcomes, and is given by
\begin{equation*}
\CVaR(f(z(\bar{x}),u,\bar{x}),\eta \mid \bar{x}) = \inf_{\alpha}\left\{\frac{\Expect_{(u \mid x=\bar{x})}\left[(f(z(\bar{x}),u,\bar{x})- \alpha)_+\right]}{\eta} + \alpha\right\},
\end{equation*}
where $(a)_+ = \max\{a,0\}$~\citep{cvaropt,cvar}.
Controlling the $\CVaR$ implies controlling the $(1-\eta)$-quantile of the performance, also known as the {\it value-at-risk} ($\VaR$), $\VaR(f(\cdot),\eta \mid \bar{x}) = \inf\{\Gamma \mid \prob_{(u \mid x = \bar{x})}(f(\cdot) \leq \Gamma) \geq 1-\eta\}$, through the following relation
\begin{equation*}
\CVaR(f(\cdot),\eta \mid \bar{x}) \leq \Gamma ~\Rightarrow~ \VaR(f(\cdot),\eta \mid \bar{x})\leq \Gamma ~\Leftrightarrow ~ \prob_{(u \mid x = \bar{x})}(f(\cdot)\leq \Gamma) \geq 1 - \eta.
\end{equation*}
We use $\CVaR$ instead of $\VaR$ due to its tractability, as well as its ability to control not only the {\it threshold} worst-case value, but also the {\it magnitude} of all the worst-case values beyond the threshold.}
The decision-maker would then consider minimizing the aggregate objective,
\begin{equation}
	\label{eq:performance_metric}
	\Expect_{(u \mid x = \bar{x})}\left(\gamma f(z(\bar{x}),u,\bar{x})\right) + \RC{\CVaR(f(z(\bar{x}),u,\bar{x}),\eta \mid \bar{x})},
\end{equation} 
where the first term describes the average case, the second term describes the worst case, and $\gamma$ controls their trade-off. 
Ideally, the problem for any context $\bar{x} \in \mathcal{D}_x$ can then be formulated as
\begin{equation}
	\label{eq:bilevel_intro_opt_cond_new}
	\begin{array}{ll}
		\underset{z(\bar{x})}{\mbox{minimize}} & \Expect_{(u \mid x = \bar{x})}\left(\Phi_\gamma(z(\bar{x}),u,\bar{x})\right) \\
		\mbox{subject to} & \RC{\CVaR(g(z(\bar{x}),u,\bar{x}),\eta \mid \bar{x})} \leq 0,\\
	\end{array}
\end{equation}
where the functions $\Phi_\gamma$ and $g$ are defined as 
\begin{equation}
	\label{eq:epigraph_reform}
\begin{aligned}
\Phi_\gamma(z(\bar{x}),u,\bar{x}) &= \gamma f(z(\bar{x}),u,\bar{x}) + \RC{z_{n+1}(\bar{x})},\\
g(z(\bar{x}),u,\bar{x}) &= \max_{l=1,\dots,L}g_l(z(\bar{x}),u,\bar{x}), \quad \text{where} \quad g_1(z(\bar{x}),u,\bar{x}) = f(z(\bar{x}),u,\bar{x}) - \RC{z_{n+1}(\bar{x})}.
\end{aligned}
\end{equation}
\RC{By augmenting $z(\bar{x})$ from dimension $n$ to $n+1$ with an epigraph variable represented by $z_{n+1}(\bar{x})$}, the objective function reformulates~\eqref{eq:performance_metric} such that the \RC{$\CVaR$} is moved into the constraint.
\RC{We enforce $f \leq z_{n+1}$ jointly with feasibility, which is stronger than imposing them separately.}
The function $\RC{g: \reals^{n+1} \times \reals^m \times \reals^{p} \to \reals}$ is the uncertain constraint that the decision-maker would like to constrain below 0 with high probability (that is, with probability $1-\eta$, as implied by the value-at-risk).
As given in~\eqref{eq:epigraph_reform}, we assume $g$ to be the maximum of $L$ functions, each of which is concave (or a maximum of concave functions) in $u$ and convex in both $z$ and $x$.
The concavity assumption is standard for methods in~\gls{RO}, and our given form encompasses a wide variety of functions that can be approximated by pointwise maxima.
This structure also allows us to express any finite number of uncertain constraints as a single joint constraint, such that the violation probability is considered holistically.
By the epigraph reformulation, the function \RC{$f - z_{n+1}$} is contained in $g$ as $g_1$, one of the $L$ constituent functions.

The conditional formulation given in~\eqref{eq:bilevel_intro_opt_cond_new} is powerful, as the constraint controls the risk for {\it each context}.
However, while the conditional guarantee is desirable, in real-life settings the conditional distribution may be difficult to estimate~\citep{conformal}.
The associated performance guarantees require stringent assumptions, which additionally depend on the approximation technique applied.
While we do provide performance guarantees for this setting (see Section~\ref{sec:cond_guarantee}), for the sake of simplicity, in the following we focus on the {\it joint formulation},
\begin{equation}
	\label{eq:bilevel_intro_opt_cond_joint}
	\begin{array}{ll}
		\underset{z(\cdot)}{\mbox{minimize}} & \Expect_{(u,x)}\left(\Phi_\gamma(z(x),u,x)\right) \\
		\mbox{subject to} & \RC{\CVaR(g(z(x),u,x),\eta)} \leq 0,\\
	\end{array}
\end{equation}
\RC{which replaces the conditional risk constraint of~\eqref{eq:bilevel_intro_opt_cond_new} by a joint one across contexts, and takes the expectation over contexts in the objective.}
Instead of optimizing for a particular context, we want to find a mapping $z(\cdot)$ that performs well {\it across all contexts}.
Such a mapping would then be applicable for any context within the problem family.

\subsection{The joint formulation and the bi-level problem}
\label{sec:joint}
In~\eqref{eq:bilevel_intro_opt_cond_joint}, we are interested in finding the best mapping $z(\cdot)$ from context parameters $x$ to decisions $z$.
To this end, we define the parametric set-valued mapping given by $x \mapsto \uncset(\theta,x)$, where 
$\uncset(\theta,x) \subseteq \reals^m$ is a convex contextual uncertainty set that depends on the {\it set parameter} $\theta \in \reals^q$.
For example, $\theta$ can represent the mapping between $x$ and the center, size, and shape of~$\uncset(\theta,x)$.
We then define the mapping $z(\cdot)$ as the solution to the \RC{\gls{RO}} problem
\begin{equation}
	\label{eq:robust_prob}
	z(\theta, x) \in \begin{array}[t]{ll}
		\mbox{argmin}  &\RC{z_{n+1}}\\
		\mbox{subject to} & g(z,u,x)  \le 0  \quad \forall u \in \uncset(\theta,x).
	\end{array}
\end{equation}
In this manner, we would like to find the optimal parameter $\theta$ by solving the bi-level optimization problem
\begin{equation}
	\label{eq:bilevel_intro_opt_theta}
	\begin{array}{ll}
		\underset{\theta}{\mbox{minimize}} & \Expect_{(u,x)}\left(\Phi_\gamma(z(\theta,x),u,x)\right) \\
		\mbox{subject to} & \RC{\CVaR(g(z(\theta,x),u,x),\eta)} \leq 0,\\
	\end{array}
\end{equation}
where the inner level problems, which become constraints of the outer level problem~\eqref{eq:bilevel_intro_opt_theta}, are represented by~\eqref{eq:robust_prob}, for all $x \in \mathcal{D}_x$. 
Problem~\eqref{eq:robust_prob} guarantees that, for $\theta$ and any given $x$, the solution $z(\theta, x)$ satisfies the constraint for any value of the uncertain parameter $u$ in $\uncset(\theta,x)$.
The outer level problem~\eqref{eq:bilevel_intro_opt_theta} then selects $\theta$ in a {\it decision-focused} manner, such that the solutions $z(\theta, x)$ perform well in terms of the performance metric and satisfy the \RC{$\CVaR$} constraint across the problem class.

\subsection{Data-driven approximation}
\label{ssec:constrained_learning}
As the distribution $\prob_{(u,x)}$ is unknown, we are unable to obtain the true expectation in~\eqref{eq:bilevel_intro_opt_theta}.
We do, however, assume access to historical data.
Suppose we are given
a pair of {\it training} datasets $W^{N} = (U^N,X^N) = \{(u^i,x^i)\}_{i=1}^N \subseteq \mathcal{D}_{(u,x)}$, consisting of a set of $N$ independent samples of the uncertain parameter $u$ and context parameter $x$, governed by $\prob^N$, their joint distribution.
As the true probabilistic guarantee in~\eqref{eq:bilevel_intro_opt_theta} cannot be obtained, we instead target
the following {\it finite-sample probabilistic guarantee}
\begin{equation}
  \label{eq:prob_guarantee}
  \prob^{N}\left(\prob_{(u,x)}(g(z(\theta,x),u,x) \le 0) \ge 1 - \eta\right)\geq 1 - \beta,
\end{equation}
where \RC{$1-\beta \in (0,1)$ is the confidence level across draws of the dataset $W^N$.}
In other words, \RC{with probability at least $1 - \beta$ over the random draw of the dataset $W^N$ from the true distribution,} the learned $\theta$ satisfies the chance constraint $\prob_{(u,x)}(g(z(\theta,x),u,x) \le 0) \ge 1 - \eta$.
To obtain such a guarantee, we first must ensure the convergence of a data-driven version of~\eqref{eq:bilevel_intro_opt_theta}.
To this end, \RC{we write
\begin{equation*}
r_i(\alpha,\theta) = \frac{(g(z(\theta,x^i),u^i,x^i)- \alpha)_+}{\eta} + \alpha, \quad i = 1,\dots,N,
\end{equation*}
for the $\CVaR$ term of sample $i$, and} we define a constraint on the sampled $\CVaR$,
\begin{equation}
\label{eq:sampled_cvar}
\widehat{\CVaR}(\theta) = \inf_{\alpha} \frac{1}{N}\sum_{i=1}^N \RC{r_i(\alpha,\theta)} \leq \kappa,
\end{equation}
and require this condition to hold for feasible $\theta$.
The parameter $\kappa < 0$ is a threshold margin that we compute the value of later in Section~\ref{sec:guarantees}, Theorem~\ref{thm:fin_sample}, to ensure the finite-sample guarantee~\eqref{eq:prob_guarantee}.
We can then formulate the stochastic objective and constraint functions as \RC{averages over the dataset $W^N$},
\begin{equation*}
F(\theta) = \frac{1}{N}\sum_{i=1}^N\Phi_{\gamma}(z(\theta,x^i),u^i, x^i) \quad \text{and}\quad H(\alpha,\theta) =\RC{\left(\frac{1}{N}\sum_{i=1}^N r_i(\alpha,\theta) - \kappa\right)_+}.
\end{equation*}
We use these approximations to define the sampled version of~\eqref{eq:bilevel_intro_opt_theta} as
\begin{equation}
	\label{eq:bilevel_sampled}
	\begin{array}{ll}
		\underset{\alpha,\theta}{\mbox{minimize}} & F(\theta) \\
		\mbox{subject to} & H(\alpha,\theta) = 0.\\
	\end{array}
\end{equation}
\RC{The constraint $H(\alpha,\theta) = 0$ is the sampled $\CVaR$ constraint~\eqref{eq:sampled_cvar} with $\alpha$ as an additional variable~\citep{cvar_var}, and the two are equivalent.
For fixed $\theta$, the function $\alpha \mapsto (1/N)\sum_{i=1}^N r_i(\alpha,\theta)$ is convex and piecewise affine, so for $\eta \in (0,1)$ its infimum is attained.
Therefore $\theta$ satisfies~\eqref{eq:sampled_cvar} if and only if $H(\alpha,\theta) = 0$ for some $\alpha$, and since $F$ does not depend on $\alpha$, problem~\eqref{eq:bilevel_sampled} and the problem with constraint~\eqref{eq:sampled_cvar} have the same optimal set parameters $\theta$.}
With this formulation, we thus tune $\alpha$ and the set parameter $\theta$ to improve the quality of the solutions $z(\theta,x^i)$ of the inner \gls{RO} problem~\eqref{eq:robust_prob}, measured by the performance on the outer level problem.

\begin{figure}[ht]\centering
	\includegraphics[width=0.7\textwidth]{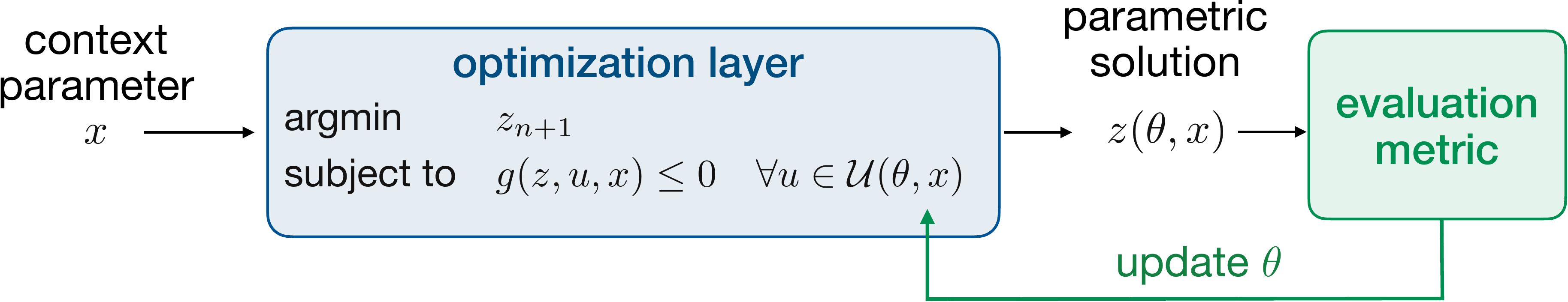}
	\caption{Our high-level procedure. 
The evaluation metric is the augmented Lagrangian function defined in Section~\ref{sec:training_procedure}, and is derived from the sampled problem~\eqref{eq:bilevel_sampled}.}
	\label{fig:arch}\end{figure}

Solving the stochastic problem~\eqref{eq:bilevel_sampled} is challenging because of the nonconvex bi-level structure introduced by the inner \gls{RO} problem~\eqref{eq:robust_prob}.
While it may be possible to directly reformulate it as a single-level problem, \RC{this would involve, for every training sample, bilinear constraints already for parameterizations affine in $\theta$, and piecewise polynomial ones for neural network maps.}
\RC{Such a problem is difficult even in the simplest case of $\ell_1$-norm uncertainty sets, as detailed in \appxdirect.}
We address these challenges in the following section, by introducing a stochastic first-order method to address problem~\eqref{eq:bilevel_sampled} with tractable iterations.
For the high-level architecture, see Figure~\ref{fig:arch}.

We give in \RC{\appxsets} examples of parameterized uncertainty sets and the corresponding convex reformulations of~\eqref{eq:robust_prob}, which we use to solve for $z(\theta,x^i)$. 
Learning $\theta$ using~\eqref{eq:bilevel_sampled} equates to learning a feasible and ideally optimal set-valued mapping from context parameters to a region of uncertainty realizations.

\subsection{\RC{Where a decision-focused set can help}}
\label{sec:gains}
\RC{Before we present our algorithm, we discuss the possible gains a decision-focused set can achieve.}
\RC{For contextual linear optimization, \cite{hu2022fast} show that, when the model class contains the truth, \gls{ETO} attains faster regret rates than decision-focused methods.
In that setting the optimal decision depends on the conditional distribution only through its mean, so a well-specified predictor of the mean is a sufficient statistic.}
\RC{The $\CVaR$ constraint of our formulation changes what the decision depends on, as the spread of the uncertainty enters, along with the division of the violation budget among the constraints.
Nevertheless, for the conditional formulation~\eqref{eq:bilevel_intro_opt_cond_new} with $\gamma = 0$, we show that the conclusion above carries over under Gaussian data for problems with a single linear uncertain constraint, where the conditional moments (mean and covariance) remain a sufficient statistic. 
When these conditions are not met, for instance with multiple linear constraints, decision-focused methods can exploit a structural gain. }\RC{Throughout this section we take $u \mid x \sim \mathcal{N}(\mu(x),\Sigma(x))$ and the ellipsoidal set
\begin{equation}
\label{eq:ellip_set_gains}
\uncset(\theta,x) = \{A(x)v + b(x) \mid \|v\|_2 \le 1\},
\end{equation}
with $\theta = (A(\cdot), b(\cdot))$ the mappings from the context to the shape and the center.
We define $\Phi$ as the standard normal distribution function, and $\varphi$ as the standard normal density.
We consider the \gls{ETO} set with $b(x) = \mu(x)$ and $A(x) = \rho\,\Sigma(x)^{1/2}$, that is, we assume it recovers the true moments, and the radius $\rho$ is chosen for a probabilistic guarantee.}

\paragraph{\RC{One uncertain constraint.}}
\RC{For a function affine in the uncertainty, $g_l = P_l(z,x)^T u + a_l(z,x)$, the robust constraint $g_l(z,u,x) \le 0$ for all $u$ in~\eqref{eq:ellip_set_gains} is
\begin{equation}
\label{eq:ellip_piece}
P_l(z,x)^T b(x) + \|A(x)^T P_l(z,x)\|_2 + a_l(z,x) \le 0.
\end{equation}
Under the \gls{ETO} set this reads $P_l^T\mu(x) + \rho\|\Sigma(x)^{1/2}P_l\|_2 + a_l \le 0$, which is equivalent to $\prob(g_l(z,u,x) \le 0 \mid x) \ge \Phi(\rho)$.
We observe the following. }

\begin{proposition}
\label{prop:single}
\RC{Consider the conditional problem~\eqref{eq:bilevel_intro_opt_cond_new} with only the epigraph piece, $L = 1$, so that the objective $f(z,u,x) = P(z,x)^T u + a(z,x)$ is uncertain and there are no other uncertain constraints.
For every context $\bar x$, the set of solutions $z$ that satisfy the $\CVaR$ constraint of~\eqref{eq:bilevel_intro_opt_cond_new} coincides with the set of solutions that satisfy the robust constraint $f(z,u,\bar x) \le z_{n+1}$ for all $u \in \uncset(\theta,\bar x)$, where $\uncset(\theta,\bar x)$ is the \gls{ETO} set with radius $\rho = \varphi(\Phi^{-1}(1-\eta))/\eta$.}
\end{proposition}
\begin{proof}
\RC{Given $\bar x$, the random variable $f(z,u,\bar x) - z_{n+1}$ is normal with mean $P^T\mu(\bar x) + a - z_{n+1}$ and standard deviation $\|\Sigma(\bar x)^{1/2}P\|_2$.
The $\CVaR$ of a standard normal variable at level $\eta$ is $\varphi(\Phi^{-1}(1-\eta))/\eta = \rho$, so the $\CVaR$ of $f(z,u,\bar x) - z_{n+1}$ is $P^T\mu(\bar x) + a - z_{n+1} + \rho \|\Sigma(\bar x)^{1/2}P\|_2$.
Requiring this quantity to be at most $0$ is~\eqref{eq:ellip_piece} with $b = \mu(\bar x)$ and $A = \rho\,\Sigma(\bar x)^{1/2}$.}
\end{proof}
\RC{The two feasible regions are identical for every context.
With $\gamma = 0$, the objective of~\eqref{eq:bilevel_intro_opt_cond_new} reduces to $z_{n+1}$ and is the same in both problems, so the true moments are already the best set parameters within the class.
In this case, no reshaping of $A$ and no shift of $b$ can enlarge the feasible region without violating the guarantee.}

\RC{For decision-focused ellipsoidal uncertainty sets to have an exploitable gain, we therefore need to depart from the conditions of Proposition~\ref{prop:single}. 
Specifically, one or more of the following must be satisfied: the conditional distribution is non-Gaussian, so that the estimated moments are no longer a sufficient statistic for the optimal decision; there are multiple uncertain constraints, so that the violation budget must be divided among them; the objective weighs both average and worst-case performance ($\gamma \neq 0$), so that the set parameters trade off the two criteria; or the guarantee is joint instead of conditional, so that the violation budget is also shared across the contexts.}

\RC{Below, we show that decision-focused learning can come out ahead in the presence of multiple uncertain constraints, through the ability to customize the violation budget.}

\paragraph{\RC{Multiple uncertain constraints.}}
\RC{With $L \ge 2$ pieces, the constraint of~\eqref{eq:bilevel_intro_opt_cond_new} implies the chance constraint $\prob(\max_l g_l(z,u,x) \le 0 \mid x) \ge 1 - \eta$.
The region of this chance constraint contains decisions that spend more of the violation budget $\eta$ on one piece and less on the others, and the optimal decision could be one such instance.
Specifically, at the optimum the {\it marginal cost} of risk is equal across the pieces, not the risk itself.
The \gls{ETO} set cannot favor one piece over another.
Its robust feasible region is
\begin{equation*}
\{z \mid \prob(g_l(z,u,x) \le 0 \mid x) \ge \Phi(\rho), \; l = 1,\dots,L\},
\end{equation*}
so every piece is held at the same marginal confidence, and the division of the risk among the pieces is fixed in advance rather than chosen for the objective.
An uncertainty set that reaches the jointly optimal decisions must protect the pieces unevenly, which a radius over $\Sigma(x)^{1/2}$ cannot do, but a reshaped $A(x)$ and shifted $b(x)$ can.}
\RC{We demonstrate this with a simple non-contextual example.}
\begin{example}
\label{ex:gains}
\RC{Let $u \in \reals^2$ have independent components $u_l \sim \mathcal{N}(\mu_l,\sigma_l^2)$, let the constraints be $g_l(z,u) = u_l - z_l$ for $l = 1,2$, and let the cost be $c^Tz$ with $c > 0$ and $\eta \in (0,1/2)$.
The goal is to satisfy the chance constraint $\prob(u_1 \le z_1,\, u_2 \le z_2) \ge 1 - \eta$, which by independence requires $\Phi(s_1)\Phi(s_2) \ge 1 - \eta$ in the standardized margins $s_l = (z_l - \mu_l)/\sigma_l$.
Since $1 - \eta > 1/2$, both margins are positive.
Under the \gls{ETO} set, we are restricted to the equal margins $\rho = s_1 = s_2$, so that $1-\eta = \Phi(\rho)^2$.
The optimal solution for the chance constraint, however, minimizes the cost $c^Tz = c^T\mu + c_1\sigma_1 s_1 + c_2\sigma_2 s_2$ subject to $\log\Phi(s_1) + \log\Phi(s_2) \ge \log(1-\eta)$, and its optimality conditions (with $\varphi$ and $\Phi$ the standard normal density and distribution function) are
\begin{equation}
\label{eq:gains_kkt}
c_l\sigma_l = \lambda\, \frac{\varphi(s_l)}{\Phi(s_l)}, \quad l = 1,2,
\end{equation}
for a multiplier $\lambda > 0$.
The optimal margins $s_1^\star,s_2^\star$ are not equal unless $c_1\sigma_1 = c_2\sigma_2$, so the \gls{ETO} set is strictly suboptimal for this chance constraint at the same confidence $1-\eta$.
The optimal margins $s_1^\star,s_2^\star$ are determined implicitly by~\eqref{eq:gains_kkt} together with $\Phi(s_1)\Phi(s_2) = 1-\eta$, and the chance-optimal decision is reached by the set with $A = \diag(s_1^\star\sigma_1, s_2^\star\sigma_2)$ and $b = \mu$, which is decision-focused rather than estimation-focused.
Geometrically, the reshaping widens the set along the axis with the smaller $c_l\sigma_l$, where margin is cheap.
For instance, with $c_1\sigma_1 = 1$, $c_2\sigma_2 = 4$, and $\eta = 0.1$, the \gls{ETO} margins are $s_1 = s_2 = 1.63$, while the optimal margins are $s_1^\star = 2.13$ and $s_2^\star = 1.37$, and the cost above $c^T\mu$ drops by about $6.6\%$.}
\end{example}

\RC{We pose the example on the chance constraint, which our $\CVaR$ constraint implies, because then the independence of the two constraints lets us characterize the optimal margins through~\eqref{eq:gains_kkt}. 
In Section~\ref{sec:intro_sets}, we empirically show that under the $\CVaR$ constraint, a decision-focused set can also outperform an~\gls{ETO} set.
We use continuous contexts, a joint formulation~\eqref{eq:bilevel_intro_opt_theta}, and additional average-case weight with $\gamma = 1$.
There, we cannot compute the optimal set parameters as above, so we learn them instead. 
}

\subsection{Motivating example}
\label{sec:intro_sets}
\RC{We consider a two-product newsvendor problem.
Each day the vendor orders $z \in \reals^2$ products at cost $k \in \reals^2$, allowing
fractional values for simplicity, and sells them at price $p \in \reals^2$ until either the
demand $u \in \reals^2$ or the inventory $z$ is exhausted.
The cost is the order cost minus the revenue,
\begin{equation*}
f(z,u,k,p) = k^Tz - p^T\min\{z,u\}
= k^Tz + \max\{-p^Tz,\; -p_1z_1 - p_2u_2,\; -p_1u_1 - p_2z_2,\; -p^Tu\},
\end{equation*}
a maximum of four affine functions in $u$.
Writing it in epigraph form gives the
\gls{RO} problem
\begin{equation}\label{eq:news_intro}
	\begin{array}{ll}
		\mbox{minimize} & t\\
		\mbox{subject to} & f(z,u,k,p) \le t \quad \forall u \in \uncset(\theta,x),\\
		& z \geq 0,
	\end{array}
\end{equation}
with context $x = (k,p)$.
We construct $\uncset(\theta,x)$ without knowing the
distribution of $(u,x)$.}

\begin{figure}[H]
	\centering
	\includegraphics[width=0.9\textwidth]{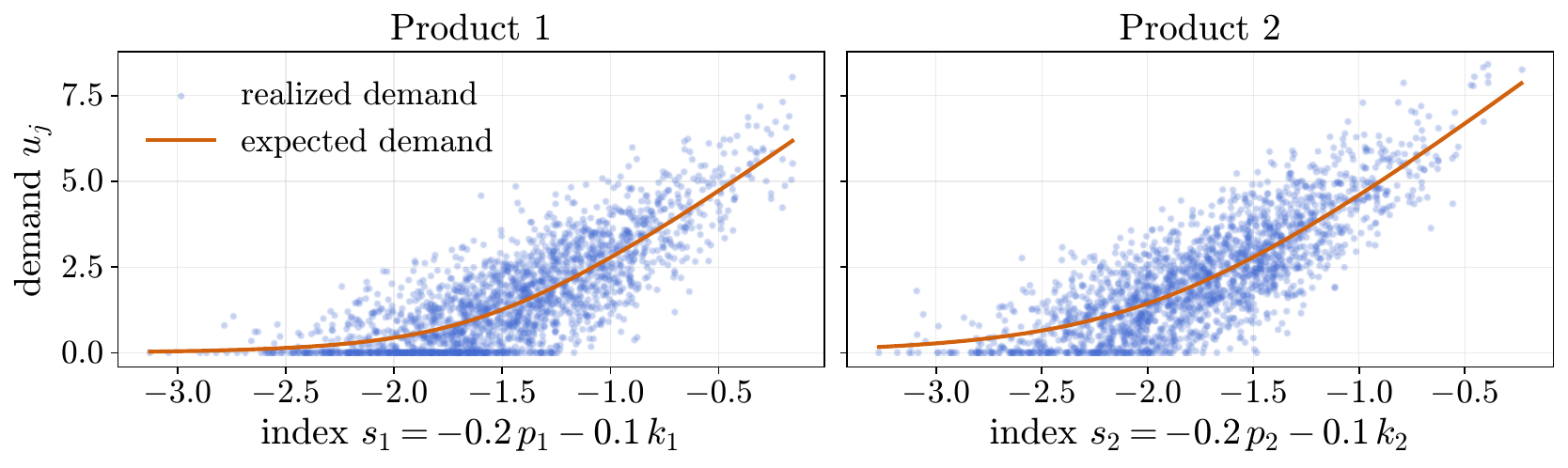}
	\caption{\RC{Demand against the index $s(x) = -0.2p - 0.1k$ for the $2000$ samples of
	one repetition, per product.
The curve is the conditional mean of
	demand~\eqref{eq:news_intro_dgp} as a function of the index. 
}
	}
	\label{fig:news_intro_context_demand}
\end{figure}

\paragraph{\RC{Problem setup.}}
\RC{For each sample, the context consists of $k = \max\{0.5,\; (4,5) + \delta_1\}$ and
$p = k + |\delta_2|$ componentwise, with $\delta_1,\delta_2$ independent normal with mean
$0$ and standard deviation $1.5$, so margins are positive almost surely. 
Demand is a latent Gaussian with context-dependent mean and identity covariance, decreasing with price,
\begin{equation}
\label{eq:news_intro_dgp}
u = \min\{20,\; \max\{0,\; \mathcal{N}(m(x), I)\}\}, \qquad
m(x) = (3,4) + s(x) + 3.6\,\sigma_s\left(\mathrm{softplus}(\tilde s(x)) - \bar{\varsigma}\right),
\end{equation}
where $s(x) = -0.2p - 0.1k$ is an affine price index, and all operations are componentwise.
The last term is a monotone distortion with mean zero, where $\tilde s(x) = (s(x) - \mu_s)/\sigma_s$, and $\bar{\varsigma}$ is the mean of $\mathrm{softplus}(\tilde s)$, which adds slight nonlinearity to the data.
All problem data is clipped to $[0,20]$ for a compact support.
See Figure~\ref{fig:news_intro_context_demand} for the resultant demand against the price index. 
We draw $2000$ context-uncertainty pairs, split into $600$ training, $400$ validation, and
$1000$ test samples, and repeat the experiment $10$ times.}

\paragraph{\RC{Methods.}}
\RC{We compare our decision-focused set, \gls{LRO}, against a set minimizing the mean-squared error of the mean and covariance mappings (MSE), and against non-robust decision-focused learning (DFL).
The first two each build a context-dependent ellipsoid with center $b(x)$ and shape $A(x)$, scaled by a radius $\rho$, with full definitions given in Section~\ref{sec:examples}.
The MSE set is an \gls{ETO} set, fitted to the data alone: its mean mapping targets the demand, and its covariance mapping targets a $k$-nearest-neighbors ($k$NN) local covariance estimate.
\gls{LRO} initializes both maps at the MSE fit and trains them through the downstream decisions.
We additionally use a {\it center-only} variant of \gls{LRO} that trains only the center, with the shape held at the MSE fit, to isolate the effect of mean-training.
DFL trains a contextual point prediction of $u$ with the same end-to-end pipeline as \gls{LRO} but with $\rho \approx 0$, and has two variants: one with only the average cost penalized, and another with an additional $\CVaR$ penalty, to represent risk-aware decision-focused learning.
For all mappings, we either use a multilayer perceptron (MLP) model, or a linear model. 
Since demand is nonnegative, every set has support $\{u \geq 0\}$.
For each method, we fit the set on the training data, calibrate its size on
the validation data for at most a $10\%$ probability of certificate violation (out-of-sample objective exceeds in-sample objective), and report
test values.
The exception is DFL, which has no set to size, so its violation probability cannot be calibrated.}

\begin{figure}[ht]
	\centering
	\includegraphics[width=\textwidth]{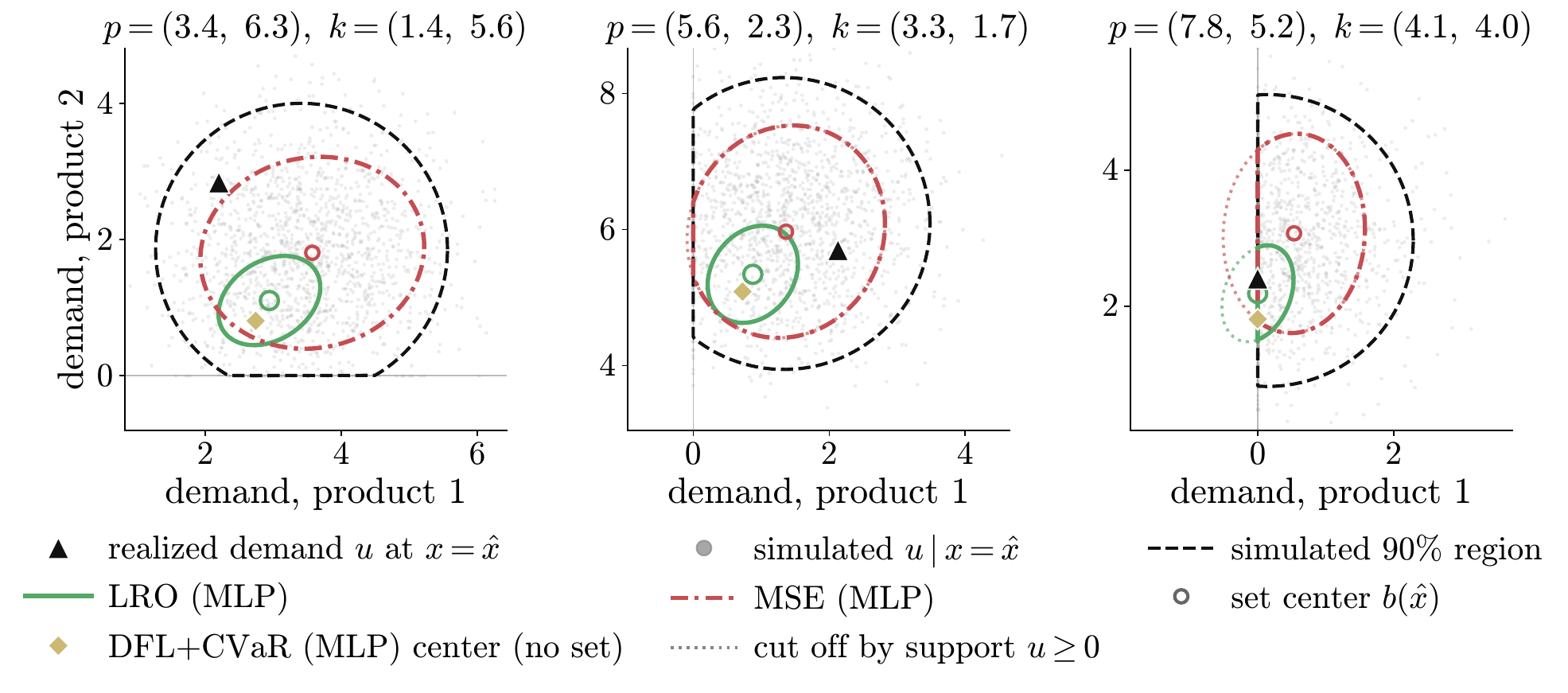}
	\caption{\RC{Uncertainty sets at three test contexts, ordered by price, all with the MLP model.
Solid and dash-dotted curves: the \gls{LRO} and MSE sets, each at its own calibrated size; dotted:
	the part removed by the support $u \geq 0$; rings: the set centers; dashed: the $90\%$
	region of the demand distribution at that context; and gray points: samples from it.
	The black triangle is the test-set demand actually realized.
The DFL method (yellow diamond) has only a point predictor, no set. 
}}
	\label{fig:news_intro_sets}
\end{figure}

\begin{figure}[ht]
	\centering
	\includegraphics[width=\textwidth]{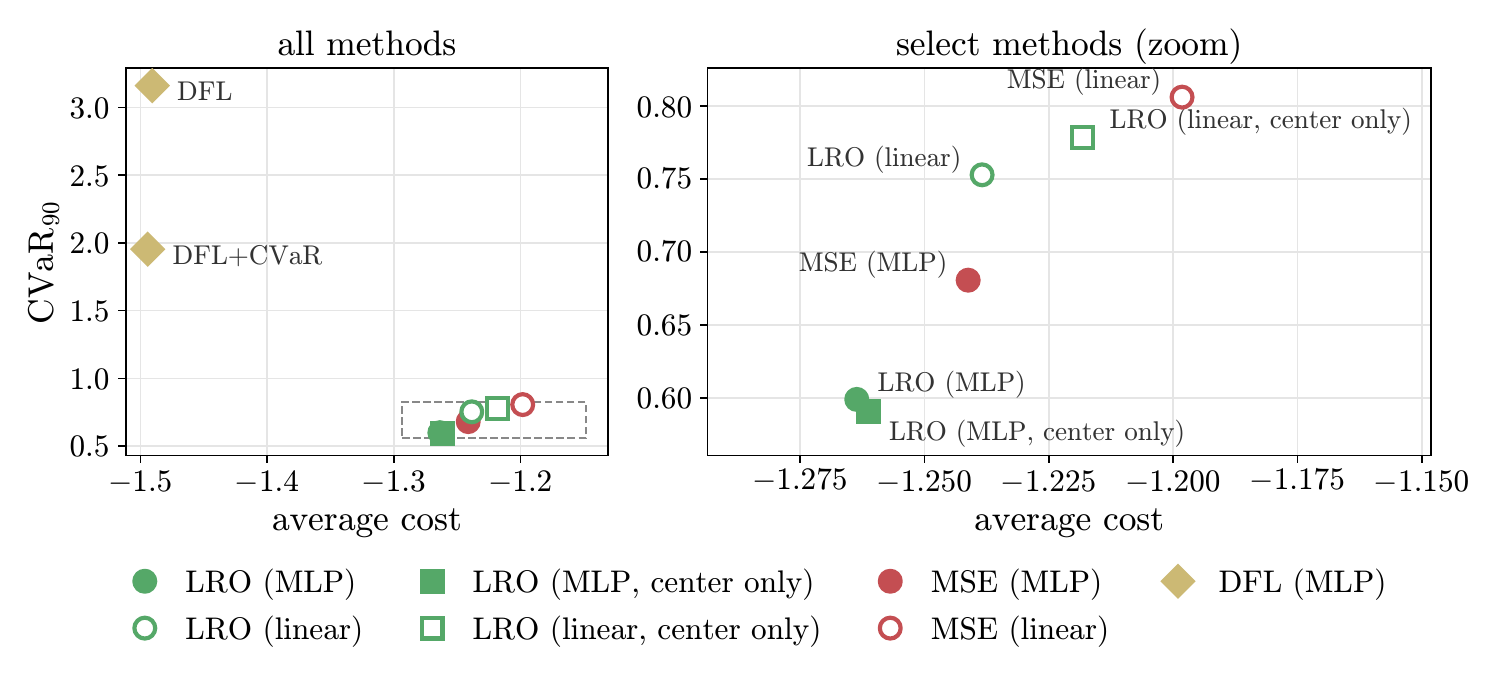}
	\caption{\RC{Out-of-sample average cost against $\CVaR$ of cost, averaged over $10$ runs.}}
	\label{fig:news_intro_avg_cvar}
\end{figure}

\paragraph{\RC{Comparison.}}
\RC{Figure~\ref{fig:news_intro_sets} shows the calibrated \gls{LRO} and MSE uncertainty sets, as well as the DFL point predictor, on three separate test contexts.
The \gls{LRO} set is the smaller of the two, containing only a portion of the simulated demand, selectively guarding the lowest plausible regions.
In fact, even the true demand (triangle) falls outside the set in the first two plots.
This is permissible by our method, since we care only about the {\it decision} and not the {\it coverage}.
The MSE set is an \gls{ETO} set and achieves high coverage, as expected.
The DFL point predictor moves toward lower demand, in the same direction as the \gls{LRO} center but further, since the point itself must guard against a low demand.}
\RC{Figure~\ref{fig:news_intro_avg_cvar} shows the resultant average-case cost against $\CVaR$ values (with risk $\eta = 0.1$) at a matched violation probability of about $10\%$ for every calibrated method.
DFL is the exception, violating at about $27\%$ with a $\CVaR$ penalty, and $38\%$ without.
Without an uncertainty set, it obtains the best average cost but suffers a large $\CVaR$, even when its training objective penalizes the $\CVaR$.
This point-prediction approach therefore trades tail risk for average cost. 
For this problem, such a trade-off may be acceptable, but when the uncertain constraint is physical, a violation rate of this size makes the decisions unusable.
Of the set-based robust methods, \gls{LRO} with the MLP models performs the best on both average cost and $\CVaR$, and the linear \gls{LRO} variants also outperform their MSE counterpart, illustrating the benefits of robust decision-focused learning.
In addition, while center-only and joint center+shape \gls{LRO} are approximately equivalent for the MLP model, shape training provides additional gains for the linear model, and is therefore important to consider.}

\RC{We have demonstrated that departing from the settings of Proposition~\ref{prop:single} allows \gls{LRO} to outperform \gls{ETO}.
The next sections detail our technique for finding these \gls{LRO} sets.}

\section{Learning the uncertainty set parameterization}
\label{sec:learn}
We now detail the procedure to solve the sampled problem~\eqref{eq:bilevel_sampled}, and prove that \RC{its accumulation points are either feasible solutions or stationary points of the constraint violation}.

\subsection{The stochastic augmented Lagrangian algorithm}
\label{sec:training_procedure}
Since our problem is nonsmooth and nonconvex due to the bi-level structure (see \RC{\appxdirect}), we search for feasible stationary points instead of globally optimal points.
\RC{An augmented Lagrangian method does so in two loops.
For a fixed multiplier $\lambda \in \reals$ and penalty parameter $\mu \in \reals$, the inner loop minimizes an augmented Lagrangian function, which adds to the objective a multiplier term and a quadratic penalty on the constraint, classically $F(\theta) + \lambda H(\alpha,\theta) + (\mu/2) H(\alpha,\theta)^2$ for problem~\eqref{eq:bilevel_sampled}.
The outer loop then updates $\lambda$ and $\mu$ to increasingly penalize constraint violation.
Evaluating the function requires solving the inner \gls{RO} problem~\eqref{eq:robust_prob} for every context in $W^N$, which could be expensive with large $N$, so we work with mini-batches instead.
At every step in the inner loop, we draw a batch of samples uniformly at random, indexed by a set $B \subseteq \{1,\dots,N\}$ of fixed size $|B|$, and we write $\Expect_B$ for the expectation over this random batch.
Our augmented Lagrangian is
\begin{equation}
\label{eq:lagrangian}
L(\alpha,\theta,\lambda,\mu) = F(\theta) + \lambda \Expect_B\big[\hat{H}_B(\alpha,\theta)\big] + \frac{\mu}{2}\Expect_B\big[\hat{H}_B(\alpha,\theta)^2\big],
\end{equation}
where $\hat{H}_B$ is the batch estimate of the constraint function, obtained by replacing the average over all $N$ samples in $H$ by the average over $B$,
\begin{equation*}
\hat{H}_B(\alpha,\theta) = \Big(\frac{1}{|B|}\sum_{i \in B} r_i(\alpha,\theta) - \kappa\Big)_+.
\end{equation*}
To minimize $L$, our inner loop evaluates at each step a batch estimate of $L$,
\begin{equation*}
\hat{L}_B(\alpha,\theta,\lambda,\mu) = \frac{1}{|B|}\sum_{i \in B}\Phi_{\gamma}(z(\theta,x^i),u^i, x^i) + \lambda \hat{H}_B(\alpha,\theta) + \frac{\mu}{2}\hat{H}_B(\alpha,\theta)^2,
\end{equation*}
whose Jacobian is an unbiased estimate of a Jacobian of $L$, since $L = \Expect_B[\hat{L}_B]$.
Although the penalty terms of $L$ differ from the full-sample $H$ and $H^2$, as $(\cdot)_+$ does not commute with $\Expect_B$, Jensen's inequality for $(\cdot)_+$ gives $\bar{H}(\alpha,\theta) = \Expect_B[\hat{H}_B(\alpha,\theta)] \geq H(\alpha,\theta) \geq 0$.
Therefore, the condition $\bar{H} = 0$ is stronger than $H = 0$, and a point with $\bar{H}(\alpha,\theta) = 0$ is also feasible for~\eqref{eq:bilevel_sampled}.
Our outer loop can then surrogate $H$ by $\bar{H}$.}

\RC{Algorithm~\ref{alg:learning_inner} gives the inner loop, which for $t \to \infty$ seeks a stationary point of $L$ at fixed $\lambda$ and $\mu$, defined below.
At each iteration, it solves the \gls{RO} problem~\eqref{eq:robust_prob} for all contexts in the drawn batch, computes the conservative Jacobian of $\hat{L}_{B^t}$, and takes a stochastic gradient step in $(\alpha,\theta)$.
Algorithm~\ref{alg:learning} gives the outer loop, which updates $\lambda$ and $\mu$ to increasingly penalize the violation value $\bar{H}(\alpha^k,\theta^k)$, which is estimated by its value on the latest inner-loop batch.}
\RC{In Section~\ref{sec:convergence_analysis} on the convergence analysis, for $t \rightarrow \infty$ and $k \rightarrow \infty$, we prove that the accumulation points of $(\alpha,\theta)$ are either feasible points of the sampled problem~\eqref{eq:bilevel_sampled}, or stationary points of the constraint violation function $V(\alpha,\theta) = (1/2)\Expect_B[\hat{H}_B(\alpha,\theta)^2]$.}
\begin{definition}[\RC{Stationary point of $L$}]
\label{def:stationary_lag}
Given $\lambda$ and $\mu$, a point $(\alpha,\theta)\in \reals^{q+1}$ is called a stationary point of the augmented Lagrangian function~\eqref{eq:lagrangian} if $0 \in J_L(\alpha,\theta)$, where $J_L(\alpha,\theta)$ denotes the conservative Jacobian (Definition~\ref{def:conserve}) of $L$ in the variables $(\alpha,\theta)$.
\end{definition}

\RC{Our algorithms} require conservative Jacobians, which are a generalized form of Jacobians for functions that are almost everywhere differentiable~\citep{conservative}, and are suitable for our case due to the nonsmooth nature of our problem.
In the following section, we detail the assumptions required for the existence of these conservative Jacobians.

\begin{algorithm}[ht]
	\caption{Stochastic gradient descent for the inner loop}
	\label{alg:learning_inner}
	\begin{algorithmic}[1]
	 \State {\bf given} $\theta^0$, $\alpha^0$, $\lambda, \mu,\{\delta^t\}_{t \in \mathbf{N}}, W^N, \RC{|B|}$
\For{$t=1,\dots,t_{\rm max}$}
	\State sample \RC{a batch $B^t \subseteq \{1,\dots,N\}$ of size $|B|$}
	 \State $z^{t}(\theta^{t-1},x^i) \gets $ solve inner \gls{RO} problem~\eqref{eq:robust_prob} for all \RC{$i \in B^t$}
	 \State $J_{\hat{L}_{\RC{B^t}}} \gets J_{\hat{L}_{\RC{B^t}}}(\alpha^{t-1},\theta^{t-1})$  \Comment{compute Jacobian of \RC{the batch estimate} $\hat{L}_{\RC{B^t}}$}
\State \RC{select one element $(v_{\alpha}^{t}, v_{\theta}^{t}) \in J_{\hat{L}_{B^t}}(\alpha^{t-1},\theta^{t-1})$}
  \State $\alpha^{t} \gets \alpha^{t-1} - \delta^{t} v_{\alpha}^{t}$
\State $\theta^{t} \gets \theta^{t-1} - \delta^t v_{\theta}^{t}$
\EndFor
	 \State {\bf return} $(\alpha^{t_{\rm max}}, \theta^{t_{\rm max}})$
	\end{algorithmic}
\end{algorithm} 

\begin{algorithm}[ht]
	\caption{Safeguarded stochastic augmented Lagrangian algorithm to solve~\eqref{eq:bilevel_sampled}}
	\label{alg:learning}
	\begin{algorithmic}[1]
	 \State {\bf given} $\theta^0 = (A^{\text{init}},b^{\text{init}})$, $\alpha^0$, $\lambda^0 \RC{\in [\lambda^{\min},\lambda^{\max}]}$, \RC{$\mu^0 > 0$}, $\kappa$, $\sigma > 1$, \RC{$\tau \in (0,1)$}, $\eta \in (0,1)$, $H_{\text{best}} = \infty$
	 \For{$k=1,\dots,k_{\rm max}$}
  \State $(\alpha^{k}, \theta^{k}) \gets$ minimize $L(\alpha, \theta, \lambda^{k-1}, \mu^{k-1})$ over $(\alpha, \theta)$ via Algorithm~\ref{alg:learning_inner}
  \If {$|\RC{\bar{H}}(\alpha^{k},\theta^{k})| \leq \tau |H_{\text{best}}|$}
	 \State $\lambda^{k} \gets  \Pi_{[{\lambda^{\min}}, \lambda^{\max}]}(\lambda^{k-1} + \mu^{k-1}\RC{\bar{H}}(\alpha^{k},\theta^{k}))$ 
  \Comment{update Lagrange multiplier}
    \State $H_{\text{best}} \gets |\RC{\bar{H}}(\alpha^{k},\theta^{k})|$
    \State $\mu^{k} \gets  \mu^{k-1}$
    \Else
	 \State $\mu^{k} \gets  \sigma\mu^{k-1}$ \Comment{update penalty}
    \State $\lambda^{k} \gets \lambda^{k-1}$
    \EndIf
\EndFor
	 \State {\bf return} $\RC{(\alpha^{k_{\rm max}}, \theta^{k_{\rm max}})}$
	\end{algorithmic}
\end{algorithm} 

\subsubsection{Path differentiability and conservative Jacobians}
In order to apply Algorithms~\ref{alg:learning_inner} and~\ref{alg:learning}, we need to establish the existence of $J_L$ and $J_{\hat{L}_{\RC{B}}}$, the conservative Jacobians of the augmented Lagrangian functions $L$ and $\hat{L}_{\RC{B}}$. 
The conservative Jacobian is defined as follows.

\begin{definition}[Conservative Jacobian~\citep{path-diff}]
\label{def:conserve}
\RC{Given a locally Lipschitz continuous function $F:\reals^{n} \rightarrow \reals^m$, we call $J_F:\reals^{n} \rightrightarrows \reals^{m\times n}$ its {\it conservative Jacobian} if $J_F$ has a closed graph, is locally bounded, and is nonempty with \ifpreprint $$\frac{\dif}{\dif t}F(\gamma(t)) = A\,\dot{\gamma}(t) \quad \text{for all } A \in J_F(\gamma(t)) \text{ and almost every } t,$$\fi
for any absolutely continuous curve $\gamma$ in $\reals^n$.} 
A locally Lipschitz function $F$ that admits a conservative Jacobian $J_F$ is called {\it path differentiable}.
\end{definition}
Notably, locally Lipschitz functions that are {\it definable} are path differentiable~\citep{path-diff}. 
To obtain the conservative Jacobians~$J_L$ and~$J_{\hat{L}_{\RC{B}}}$, we then require the following. 
\begin{assumption} [Definability]
    \label{ass:define}
All functions $f$ and $g_l$ in the robust optimization problem~\eqref{eq:robust_prob} are {\it definable} in $z$~\citep{definable,definable2,path-diff}.
\RC{In addition, the maps $(\theta,x) \mapsto (A(x),b(x))$ that parameterize the uncertainty set $\uncset(\theta,x)$, and any problem data that depend on the context $x$, are definable in $(\theta,x)$.}
\end{assumption} 
\begin{assumption} [Local Lipschitzness]
    \label{ass:lip}
    All functions $f$ and $g_l$ in the robust optimization problem~\eqref{eq:robust_prob} are locally Lipschitz in $z$.
\RC{In addition, the maps $(\theta,x) \mapsto (A(x),b(x))$ that parameterize the uncertainty set $\uncset(\theta,x)$, and any problem data that depend on the context $x$, are locally Lipschitz in $(\theta,x)$.}
A function $\ell$ is locally Lipschitz if, for each point $z$ in its domain, there exists a neighborhood $\mathcal{V}$ of $z$ such that $\ell$ is Lipschitz on $\mathcal{V}$.
\end{assumption}
The class of definable functions is very wide and includes semi-algebraic functions\RC{, that is,} functions whose graphs can be described by a finite number of polynomial equalities and inequalities.
This includes polynomials, the Euclidean norm, maxima and minima of finitely many polynomials, as well as piecewise polynomial functions, which incorporate a large number of functions we would expect to use in stochastic optimization problems. 
For a more detailed overview, see the works by~\cite{conservative} and~\cite{path-diff}.
For a foundational work, see~\cite{definable}. 

Assuming the definability of the objective and constraint functions, it follows that $F$ and $H$, and subsequently $L$, as compositions of definable pieces, are definable in the variables $(z,\alpha)$~\citep{path-diff}.
The same follows for the local Lipschitzness condition, which is preserved under composition.
\RC{Both model classes of the motivating example (Section~\ref{sec:intro_sets}) satisfy the conditions on the maps, since the linear maps are polynomial in $(\theta,x)$ and a ReLU network is piecewise polynomial, hence semi-algebraic and locally Lipschitz.}
By the chain rule for path differentiable functions, $L$ is then path differentiable in $\theta$ if $z(\theta,x)$ is itself path differentiable.
We observe that $z(\theta,x)$ is the solution map of the inner robust problem~\eqref{eq:robust_prob}, which admits a convex conic reformulation\RC{, given in} Appendix~\ref{sec:dpp}.
\RC{Path differentiability of the solution map then requires the following regularity condition on the conic reformulation.}

\begin{assumption} [Regularity of the inner conic program]
    \label{ass:unique}
    \RC{For every $(\theta,x) \in \mathcal{D}_\theta \times \mathcal{D}_x$, the convex conic reformulation of the \gls{RO} problem~\eqref{eq:robust_prob}, whose cone is a product of zero cones, nonnegative orthants, and second-order cones (Appendix~\ref{sec:dpp}), has a primal-dual solution that satisfies primal constraint nondegeneracy~\citep[Section 4.6]{bonnans2000} and the strong second-order sufficient condition~\citep{bonnans2005}.}
\end{assumption}
\RC{These are the standard conditions under which the solution of a conic program is locally unique and locally Lipschitz in its data~\citep[Chapters 4 and 5]{bonnans2000}.}
\RC{The first holds for generic problem data, and in specific cases the second condition can be ensured by a quadratic penalty $\varepsilon\|z\|^2$ over the decisions $z$ (excluding the epigraph variable), which is common practice in differentiable optimization and decision-focused learning~\citep{amos2021optnet,end2,dfl}.
We keep the assumption explicit and delay details to Appendix~\ref{sec:dpp}.}
We now establish the path differentiability of our solution. 
\begin{theorem}[Path differentiability] 
\label{thm:subdiff}
	\RC{Let Assumptions~\ref{ass:define} to~\ref{ass:unique} hold.
Then, for every $(\theta,x) \in \mathcal{D}_\theta \times \mathcal{D}_x$, the conic reformulation~\eqref{eq:conic} of the inner optimization problem~\eqref{eq:robust_prob} has a unique primal-dual solution, and its primal part $z(\theta,x)$ is locally Lipschitz in $(\theta,x)$ and path differentiable in $\theta$, with nonempty conservative Jacobian $J_{z}(\theta)$.}
\end{theorem}
The proof is given in Appendix~\ref{appx:subdiff}.
By application of the chain rule, we can then establish the existence of the conservative Jacobians of the augmented Lagrangian function $L$.
The exact forms are given in Appendix~\ref{appx:jacobians}.
\begin{proposition} [Conservative Jacobians] 
\label{prop:jacobians}
Suppose Assumptions~\ref{ass:define} to~\ref{ass:unique} hold, and for each given input pair $(\theta, x)$, the solution $z(\theta,x)$ of the optimization problem~\eqref{eq:robust_prob} is path differentiable with conservative Jacobian $ J_{z}(\theta)$.
Then, the Lagrangian function $L(\alpha,\theta,\lambda,\mu)$, from equation~\eqref{eq:lagrangian}, and its batched version\RC{s $\hat{L}_B(\alpha,\theta,\lambda,\mu)$}, are path differentiable\RC{, with conservative Jacobians $J_{\hat{L}_B}(\alpha,\theta)$ and $J_L(\alpha,\theta) = \Expect_B[J_{\hat{L}_B}(\alpha,\theta)]$}.
\end{proposition}

\subsubsection{Convergence analysis}
\label{sec:convergence_analysis}
\RC{To analyze the augmented Lagrangian algorithm, given in Algorithm~\ref{alg:learning}, we begin with the inner loop, Algorithm~\ref{alg:learning_inner}.}
In addition to the path-differentiability of the Lagrangian function established above, we require that the iterates are bounded, and that the step sizes are nonnegative, square summable, but not summable.
\begin{assumption}[Boundedness]
\label{ass:bounded}
There exist a constant $\alpha^{\max} > 0$ and a compact set \RC{$\mathcal{D}_{\alpha} \subseteq [-\alpha^{\max},\alpha^{\max}]$}, such that for $\alpha^0 \in  \mathcal{D}_{\alpha}$, \RC{$\alpha^t \in \mathcal{D}_{\alpha}$} almost surely.
Similarly, there exist constants $\theta^{\max}, z^{\max}, x^{\max} > 0$, compact sets \RC{$\mathcal{D}_{\theta} \subseteq \{\theta \in \reals^{q} \mid \|\theta\|_2 \leq \theta^{\max}\}$}, $\mathcal{D}_x \subset \reals^p$, $\mathcal{Z} \subset \reals^{n+1}$, such that for $\theta^0 \in \mathcal{D}_{\theta}$ and $x \in \mathcal{D}_x$, we have \RC{$\theta^t \in \mathcal{D}_{\theta}$}, $\|x\| \leq x^{\max}$, and $z(\theta^t,x) \in \mathcal{Z} = \{z \in \reals^{n+1} \mid \|z\| \le z^{\max}\}$ almost surely.
\end{assumption}
\begin{assumption}[Step sizes]
\label{ass:step_size}
The step sizes $\{\delta^t\}_{t \in \mathbf{N}}$ are chosen such that
\begin{align*}
   &\sum_{t=1}^\infty \delta^t = \infty, \quad \sum_{t=1}^\infty (\delta^t)^2 < \infty,\quad \delta^t = o(1/\log(t)).
\end{align*}
\end{assumption}
\RC{Under these assumptions, Theorem 3 of~\cite{path-diff} shows that, at fixed $\lambda$ and $\mu$ and for almost every initialization and scaling of the step sizes, the values $L(\alpha^t,\theta^t,\lambda,\mu)$ along Algorithm~\ref{alg:learning_inner} converge almost surely as $t \to \infty$, and every accumulation point of $\{(\alpha^t,\theta^t)\}$ is a stationary point of $L$.}
\RC{The penalty parameter $\mu^k$ of Algorithm~\ref{alg:learning}, however, may either remain bounded or diverge as $k \to \infty$, and this determines the nature of the accumulation points $(\bar\alpha,\bar\theta)$ of $\{(\alpha^k,\theta^k)\}$.}
\begin{theorem}
\label{thm:convergence}
Suppose Assumptions~\ref{ass:define}~to~\ref{ass:step_size} hold\RC{, and that each inner loop is run to a stationary point of $L(\cdot,\cdot,\lambda^{k-1},\mu^{k-1})$ in $\mathcal{D}_\alpha \times \mathcal{D}_\theta$.
Then, almost surely, every accumulation point $(\bar\alpha,\bar\theta)$ of $\{(\alpha^k,\theta^k)\}$ is a feasible point of the sampled problem~\eqref{eq:bilevel_sampled} if the sequence of penalty parameters $\{\mu^k\}$ is bounded, and a stationary point of $V(\alpha,\theta) = (1/2)\Expect_B[\hat{H}_B(\alpha,\theta)^2]$ otherwise.}
\end{theorem}
\RC{The proof is delayed to Appendix~\ref{appx:converge_stationary}.
As in the smooth safeguarded theory~\citep[Theorem 6.3]{auglag_safe}, the safeguards $[\lambda^{\min},\lambda^{\max}]$ of Algorithm~\ref{alg:learning} keep the multiplier $\lambda^k$ bounded, so that when the penalty $\mu^k$ diverges it dominates the augmented Lagrangian, and the accumulation points are at least stationary points of the constraint violation function $V$. }
\begin{remark}
	\label{rem:conv}
\RC{Theorem~\ref{thm:convergence} does not establish that the whole sequence $\{\theta^k\}$ converges, only properties of its accumulation points.
Full-sequence convergence would require local uniqueness and a stability analysis of the minimizers of the nonconvex and nonsmooth augmented Lagrangian in $(\lambda,\mu)$, which needs stronger assumptions that are difficult to verify here.}
	
\RC{Furthermore, while Theorem~\ref{thm:convergence} assumes each inner loop reaches a stationary point of $L$, and its conclusions involve accumulation points of the outer iterates, the given algorithms truncate the inner loop at a finite $t_{\text{max}}$ and the outer loop at a finite $k_{\text{max}}$, and practical implementation may also have summable step sizes. 
The convergence result for the inner loop therefore gives no bound on the stationarity error of the point returned after $t_{\text{max}}$ steps, and likewise none is given for the outer loop at $k_{\text{max}}$.
However, this mirrors standard practice in the asymptotic convergence theory that we build upon.
For instance, the stochastic subgradient and path-differentiable convergence results of~\cite{Davis} and~\cite{path-diff}, and the safeguarded augmented Lagrangian theory of~\cite{auglag_safe}, are similarly stated for infinite iteration limits, and applied in practice with a finite, truncated implementation.
To work with our practical iteration budget, our implementation replaces the stochastic gradient step of Algorithm~\ref{alg:learning_inner} by Adam~\citep{kingma2015adam} and the multiplier update of Algorithm~\ref{alg:learning} by a PID controller~\citep{stooke2020responsive}, which offer better performance but have weaker convergence guarantees.
For efficiency, we also estimate $\bar{H}$ on the latest batch rather than compute $H$ exactly, and do not project our iterates onto $\mathcal{D}_\alpha \times \mathcal{D}_\theta$.
In our examples in Section~\ref{sec:examples}, we then search over the radius $\rho$ of the returned set for a value that meets the empirical constraint, which may still provide the guarantees of the following section. 
	The radius calibration also serves a second purpose.
The batched value $\bar{H} = 0$ if and only if the sampled $\CVaR$ constraint~\eqref{eq:sampled_cvar} is satisfied at level $\eta|B|/N$ instead of $\eta$, which may be conservative, and recalibrating the radius on a larger validation set removes this conservatism.}
\end{remark}

\section{Probabilistic guarantees of constraint satisfaction}
\label{sec:guarantees}
\RC{A feasible solution of the sampled problem~\eqref{eq:bilevel_sampled} does not by itself provide the probabilistic guarantee of constraint satisfaction~\eqref{eq:prob_guarantee}, which requires generalization to unseen data.} 
\RC{In this section, we establish this guarantee, and in fact establish it for any (approximately) empirically feasible solution satisfying $H(\alpha^\star,\theta^\star) \leq \epsilon$, for some $\epsilon \geq 0$.}
We prove two results: one with respect to the joint distribution, as given in~\eqref{eq:prob_guarantee}, and another with respect to the conditional distribution, for the setting introduced in Section~\ref{sec:cond}.

\subsection{Guarantees with respect to the joint distribution}
To arrive at our desired guarantees, where the inner probability is with respect to the unknown joint distribution~$\prob_{(u,x)}$, we need to analyze the generalization behavior of the robust solutions on out-of-sample data. 
In addition to Assumption~\ref{ass:bounded}, we need a boundedness assumption on the function from the $\CVaR$ definition.
\begin{assumption}[Bounded range]
  \label{assup:bound_range} The function $(g(z(\theta,x),u,x)- \alpha)_+/\eta + \alpha$ maps to values in $(-c^{\min}, c^{\max})$, where $0 \leq c^{\min}, c^{\max} < \infty$, for all $\theta,\alpha,u,x$ in its support.
\end{assumption}

Under Assumption~\ref{assup:bound_range}, we then let $C = c^{\min} + c^{\max}$ and define the normalized function
\begin{equation}
	\label{eq:g_setting}
\psi_{w}(u, x) = \frac{1}{C} \left(c^{\min} + \frac{(g(z(\theta,x), u, x) - \alpha)_{+}}{\eta} + \alpha\right),
\end{equation}
where $w = (\theta, \alpha) \in \mathcal{W} \subset \reals^{q+1}$.
Note that when Assumption~\ref{ass:bounded} holds, the set $\mathcal{W} = \mathcal{D}_{\theta} \times \mathcal{D}_{\alpha}$ is compact \RC{and contains every iterate of Algorithm~\ref{alg:learning}}.
The function $\psi_{w}$ then maps to values in $(0,1)$, and we consider the entire function class
$\Psi =  \{\psi_{w}\}_{w \in\mathcal{W}},$
which we use to acquire the finite-sample probabilistic guarantee in Theorem~\ref{thm:fin_sample}.
\RC{We control the deviation of empirical means from their expectations uniformly over $\Psi$ with the following bound.
\begin{lemma}[Uniform deviation bound]
\label{lem:unif_dev}
Let Assumption~\ref{ass:bounded} hold, and let $\psi_w(u,x) \in [0,1]$ be $L_\psi$-Lipschitz in $w \in \mathcal{W}$, uniformly in $(u,x)$.
Let $(u^i,x^i)$, $i = 1,\dots,N$, be independent, not necessarily identically distributed, and let $\bar{\psi}_w = (1/N)\sum_{i=1}^N \psi_w(u^i,x^i)$ and $m_w = (1/N)\sum_{i=1}^N \Expect[\psi_w(u^i,x^i)]$.
Then, for every $\tau > 0$,
\begin{equation*}
\prob\left(\sup_{w \in \mathcal{W}}\left|\bar{\psi}_w - m_w\right| > \tau\right) \leq \beta(N,\tau) = 2\left(1 + \frac{16 L_\psi w^{\max}}{\tau}\right)^{q+1}\exp\left(-\frac{N\tau^2}{2}\right),
\end{equation*}
where $w^{\max} = \theta^{\max} + \alpha^{\max}$, with $\theta^{\max}$ and $\alpha^{\max}$ the constants of Assumption~\ref{ass:bounded}.
\end{lemma}
\begin{proof}
If $L_\psi = 0$, then $\bar{\psi}_w - m_w$ does not depend on $w$ and Hoeffding's inequality gives the bound directly, so let $L_\psi > 0$.
Let $\delta = \tau/(4L_\psi)$ and let $\mathcal{N} \subseteq \mathcal{W}$ be a $\delta$-net of $\mathcal{W}$.
By Assumption~\ref{ass:bounded}, every $w = (\theta,\alpha) \in \mathcal{W}$ satisfies $\|w\|_2 \leq \|\theta\|_2 + |\alpha| \leq w^{\max}$, so $\mathcal{W}$ lies in the ball of radius $w^{\max}$ in $\reals^{q+1}$.
This ball admits a $(\delta/2)$-net of at most $(1 + 4w^{\max}/\delta)^{q+1}$ points~\citep[Section 4.2]{vershynin2018}, and replacing each point of that net that lies within $\delta/2$ of $\mathcal{W}$ by a point of $\mathcal{W}$ at distance at most $\delta/2$ gives $|\mathcal{N}| \leq (1 + 4w^{\max}/\delta)^{q+1}$.
For any $w \in \mathcal{W}$, pick $w' \in \mathcal{N}$ with $\|w - w'\|_2 \leq \delta$.
By the Lipschitz condition, $|\bar{\psi}_w - \bar{\psi}_{w'}| \leq L_\psi\delta$ and $|m_w - m_{w'}| \leq L_\psi\delta$, so $|\bar{\psi}_w - m_w| \leq |\bar{\psi}_{w'} - m_{w'}| + \tau/2$.
For each fixed $w' \in \mathcal{N}$, the summands $\psi_{w'}(u^i,x^i)$ are independent with values in $[0,1]$, so Hoeffding's inequality gives $\prob(|\bar{\psi}_{w'} - m_{w'}| > \tau/2) \leq 2\exp(-N\tau^2/2)$.
A union bound over the at most $(1 + 4w^{\max}/\delta)^{q+1}$ points of $\mathcal{N}$ multiplies this probability by $|\mathcal{N}|$, which gives $\beta(N,\tau)$.
\end{proof}
}
\begin{theorem}[\RC{Finite-sample guarantee}]
	\label{thm:fin_sample}
Suppose Assumptions~\ref{ass:define} to~\ref{ass:bounded} \RC{and~\ref{assup:bound_range}} hold\RC{, and that $\psi_w$ is Lipschitz in $w \in \mathcal{W}$ with constant $L_\psi$, uniformly in $(u,x)$.
Let $\tau > 0$, $\epsilon \geq 0$, and let the margin satisfy $\kappa \leq -(\epsilon + C\tau)$}.
Then, \RC{with probability at least $1-\beta$ over the draw of $W^N$, every $(\alpha,\theta) \in \mathcal{W}$ with $H(\alpha,\theta) \leq \epsilon$ satisfies the chance constraint, that is,}
\begin{equation}
	\label{eq:prob_guarantee_fin}
	\prob^{N}\left (\RC{\prob_{(u,x)}(g(z(\theta,x),u,x)\leq 0) \geq 1 - \eta \ \text{ for all } (\alpha,\theta) \in \mathcal{W} \text{ with } H(\alpha,\theta) \leq \epsilon}\right) \geq 1 - \beta,
	\end{equation}
	where \RC{$\beta = \beta(N,\tau)$ is the bound of Lemma~\ref{lem:unif_dev}}.
\end{theorem}

\begin{proof}
\RC{Let $w^\star = (\theta^\star, \alpha^\star) \in \mathcal{W}$ be any point, possibly data-dependent, with $H(\alpha^\star,\theta^\star) \leq \epsilon$.
By the definition of $H$, this gives $(1/N)\sum_{i=1}^N r_i(\alpha^\star,\theta^\star) \leq \kappa + \epsilon$, and by the normalization~\eqref{eq:g_setting},
$$\frac{1}{N}\sum_{i=1}^N \psi_{{w}^{\star}}(u^i, x^i) \leq \frac{c^{\min} + \kappa+ \epsilon}{C}.$$
Since $w^\star$ depends on the data, we need the deviation bound uniformly over $\mathcal{W}$.
The samples in $W^N$ are independent, so Lemma~\ref{lem:unif_dev} gives, with probability at least $1 - \beta(N,\tau)$ over the draw of $W^N$,
$$\Expect_{(u,x)}[\psi_{w^\star}(u,x)] \leq \frac{c^{\min}+ \kappa+\epsilon}{C} + \tau.$$
Undoing the normalization, with $r_{w}(u,x) = (g(z(\theta,x),u,x)- \alpha)_+/\eta + \alpha$, we obtain
$$\Expect_{(u,x)}[r_{w^\star}(u,x)] = C\,\Expect_{(u,x)}[\psi_{w^\star}(u,x)] - c^{\min} \leq \kappa + \epsilon + C\tau \leq 0.$$
As the $\CVaR$ is the infimum of this expectation over $\alpha$, the $\CVaR$ of $g(z(\theta^\star,x),u,x)$ at level $\eta$ is nonpositive on this event, which implies the chance constraint.
The event of Lemma~\ref{lem:unif_dev} does not depend on $w^\star$ and has probability at least $1-\beta(N,\tau)$, so it covers every such point simultaneously, which gives~\eqref{eq:prob_guarantee_fin}.}
\end{proof}
\RC{This guarantee covers any data-dependent point returned by Algorithm~\ref{alg:learning} that lies in $\mathcal{W}$ and satisfies $H \leq \epsilon$.
A radius recalibrated after training is covered when the recalibrated parameter remains in $\mathcal{W}$ and still satisfies $H \leq \epsilon$.
The bound $\beta(N,\tau)$ is nonvacuous once $N\tau^2/2 > \log 2 + (q+1)\log(1 + 16 L_\psi w^{\max}/\tau)$, so it requires $N$ to be large relative to the dimension $q$ of the set parameter.
Note that for our numerical experiments, which do not use a large $N$ compared to $q$, this guarantee is only qualitative.}
\RC{The Lipschitz condition on $\psi_w$ holds} when the constraint function $g$ is Lipschitz in $z$, and $z$ is uniformly Lipschitz in $\theta$ over $\mathcal{D}_\theta$.
For instance, when $g$ is affine or a maximum of affine functions in the variable $z$, which is a common scenario in many optimization problems (and is the case for our numerical experiments in Section~\ref{sec:examples}), it is Lipschitz in $z$ for all $(u,x) \in \mathcal{D}_{(u,x)}$,
\ie,
\[
|g(z,u,x) - g(z',u,x)| \leq \bar{L}\|z - z'\|_2 \quad \text{for all } z,z'\in \mathcal{Z},~(u,x) \in \mathcal{D}_{(u,x)},
\]
with some Lipschitz constant $\bar{L}$.
\RC{By Theorem~\ref{thm:subdiff}, $z(\theta,x)$ is locally Lipschitz in $(\theta,x)$.}
Since $\mathcal{D}_\theta$ and $\mathcal{D}_x$ are both compact by Assumption~\ref{ass:bounded}, \RC{and a locally Lipschitz function on a compact set is Lipschitz on it, there is} a uniform constant $L_z$ such that
\[
\|z(\theta,x) - z(\theta',x)\|_2 \leq L_z\|\theta - \theta'\|_2 \quad \text{for all } \theta,\theta'\in \mathcal{D}_\theta,~x\in \mathcal{D}_x.
\]
It follows that $g$ is Lipschitz in $\theta$ with constant $\bar{L}L_z$, and $\psi_w \in \Psi$ is Lipschitz in $w=(\theta,\alpha)$ with Lipschitz constant $\RC{L_\psi} \leq (1/C)((\bar{L}L_z+1)/\eta + 1)$.

\subsection{Guarantees with respect to the conditional distribution}
\label{sec:cond_guarantee}
To establish performance guarantees with respect to the conditional distribution $\prob_{(u \mid x = \bar{x})}$ instead of the joint distribution $\prob_{(u,x)}$, we now consider the optimization problem~\eqref{eq:bilevel_intro_opt_cond_new},
for which an approximated sampled formulation (an alternative to~\eqref{eq:bilevel_sampled}) becomes
\begin{equation}
	\label{eq:bilevel_sampled_cond}
	\begin{array}{ll}
		\underset{\alpha,\theta}{\mbox{minimize}} & \tilde{F}(\theta,\bar{x}) \\
		\mbox{subject to} & \tilde{H}(\alpha,\theta,\bar{x}) = 0,\\
	\end{array}
\end{equation}
where
$$\tilde{F}(\theta,\bar{x}) = \sum_{i=1}^N w(\bar{x},x^i)\Phi_{\gamma}(z(\theta,x^i),u^i, x^i),$$
and 
 \RC{$$\tilde{H}(\alpha,\theta,\bar{x}) = \left(\sum_{i=1}^N w(\bar{x},x^i)\left((g(z(\theta,x^i),u^i,x^i)- \alpha)_+/\eta + \alpha\right) - \kappa\right)_+,$$}
consistent with the aggregate form of $H(\alpha,\theta)$ above.
The weights $w(\bar{x},x^i)$ are constructed such that $\sum_{i=1}^Nw(\bar{x},x^i) = 1$ and $w(\bar{x},x^i) \geq 0$ for $i = 1,\dots,N$, and are a measure of the proximity of the data point $x^i$ to $\bar{x}$.
With this construction, the conditional distribution $\prob_{(u \mid x = \bar{x})}$ is approximated as
$$\tilde{\prob}_{(u \mid x = \bar{x})} = \sum_{i=1}^Nw(\bar{x},x^i)\delta_{u^i}.$$

As discussed by~\cite{bertsimas2020predictive}, the weights could be constructed using $k$-nearest-neighbors ($k$NN), kernel methods, nonnegative local linear methods, trees, or random forests.
Here, we focus on the $k$NN setting, and note that the remaining methods could similarly be applied, albeit with different assumptions.

In the $k$NN setting, we choose the weights
$w(\bar{x},x^i) = (1/k)\ones\{x^i\text{ is a $k$NN of } \bar{x}\}$,
where ties are broken randomly, and the nearest neighbors are determined using the Euclidean distance in $\reals^p$. 
In their work,~\cite{bertsimas2020predictive} derive asymptotic guarantees for their constructed contextual prescriptors.
We, however, require a finite-sample guarantee analogous to~\eqref{eq:prob_guarantee_fin} for a particular context $\bar{x}$.
We therefore require the following assumptions:
\begin{assumption}[Regular marginal density]
	\label{ass:positive}
	The distribution $\prob_x$ of the context variable $x$ admits a density $\mu_x$ with respect to the Lebesgue measure\RC{, and there exist a radius $R > 0$ and a constant $\mu_{\min} > 0$ such that for all $r \in (0, R]$,
	$\int_{B(\bar{x},r)} \mu_x(x)\, dx \geq \mu_{\min}\, V_p\, r^p$},
	where $B(\bar{x},r) = \{x \in \reals^p \mid \|x - \bar{x}\|_2 < r\}$ is the open ball of radius $r$ centered at $\bar{x}$, and $V_p$ is the volume of the unit ball in $\reals^p$.
\end{assumption}
\begin{assumption}[H\"{o}lder continuous conditional expectation]
	\label{ass:holder}
There exist constants $L > 0$ and $s > 0$ such that for all $w = (\theta, \alpha) \in \mathcal{W}$ and all $\bar{x}, \bar{x}' \in \mathcal{D}_x$,
\begin{equation*}
    \left| \Expect_{(u \mid {x = \bar{x}})}[\psi_w(u, \bar{x})] - \Expect_{(u \mid {x = \bar{x}}')}[\psi_w(u, \bar{x}')] \right| \leq L \| \bar{x} - \bar{x}' \|_2^s.
\end{equation*}
\end{assumption}
The density assumption ensures that the conditional distribution at $\bar{x}$ is meaningful and that the local mass of the ball around $\bar{x}$ is well-behaved, while the continuity assumption ensures that borrowing information from neighboring contexts introduces only a controlled amount of bias.
We can then derive the finite-sample guarantee for the conditional distribution. 
\begin{theorem}[Conditional guarantee]
	\label{thm:fin_sample_cond}
Suppose Assumptions~\ref{ass:define} to~\ref{ass:bounded} \RC{and~\ref{assup:bound_range} to~\ref{ass:holder}} hold\RC{, that $\prob_x$ is supported on $\mathcal{D}_x$, and that $\psi_w$ is $L_\psi$-Lipschitz in $w \in \mathcal{W}$, uniformly in $(u,x)$.
Fix $\tau_{\mathrm{var}} > 0$ and $\epsilon \geq 0$, let $\tau_{\mathrm{bias}}$ be as below, and choose the margin $\kappa \leq -(\epsilon + C(\tau_{\mathrm{var}} + \tau_{\mathrm{bias}}))$}.
\RC{Then, for a particular context $\bar{x} \in \mathcal{D}_x$, with probability at least $1-\beta_{\mathrm{cond}}$, every $(\alpha,\theta) \in \mathcal{W}$ with $\tilde{H}(\alpha,\theta,\bar{x}) \leq \epsilon$ satisfies the conditional chance constraint at $\bar{x}$, that is,}
\begin{equation}
	\label{eq:prob_guarantee_cond}
	\begin{aligned}
	\prob^{N}\Big(&\RC{\prob_{(u \mid x = \bar{x})}(g(z(\theta,\bar{x}),u,\bar{x})\leq 0) \geq 1 - \eta} \\
	&\RC{\text{for all } (\alpha,\theta) \in \mathcal{W} \text{ with } \tilde{H}(\alpha,\theta,\bar{x}) \leq \epsilon}\Big) \geq 1 - \beta_{\mathrm{cond}},
	\end{aligned}
\end{equation}
where $\beta_{\mathrm{cond}} = \beta_{\mathrm{var}} + \beta_{\mathrm{bias}}$.
The variance term is \RC{$\beta_{\mathrm{var}} = \beta(k,\tau_{\mathrm{var}})$, the bound of Lemma~\ref{lem:unif_dev} with $k$ samples in place of $N$}.
The value $\tau_{\mathrm{bias}}$ is chosen as $\tau_{\mathrm{bias}}= L r^s$ for an \RC{$r \in \left[\left(k/(\mu_{\min} V_p N)\right)^{1/p},\; R\right]$, an interval that is nonempty when $k \leq \mu_{\min} V_p R^p N$}, and the corresponding bias probability is
$$\beta_{\mathrm{bias}} = \exp\left(-\RC{\mu_{\min}}\, V_p\, r^p\, N\right)\left(\frac{e\,\RC{\mu_{\min}}\, V_p\, r^p\, N}{k}\right)^k.$$
\end{theorem}

\RC{The proof is given in Appendix~\ref{appx:cond_proof}.}
We note that due to the bias term, and due to the use of $k$ instead of $N$ in the variance term, this result adds conservatism compared to Theorem~\ref{thm:fin_sample}.

\section{Hyperparameter tuning}
\label{sec:hyper}
As our learning problem is stochastic, nonsmooth, and nonconvex, our feasible solution is not guaranteed to be a locally or globally optimal solution. 
Furthermore, as the theoretical guarantees are often conservative, they require a very large number of samples~$N$. In practice, we tune for a target guarantee using a validation dataset, employing hyperparameter tuning to select the best among various feasible solutions.
\RC{As opposed to traditional methods with limited hyperparameters (for instance, only the radius), our flexibility allows for improved solutions, with the trade-off of increased tuning time.}

\RC{Our modeling choices have varying degrees of importance.
The first is the model class, that is, whether the maps are linear, a multilayer perceptron (MLP), or another neural network.
As seen in the motivating example, this choice strongly affects the solution.
The second is the relative strength (learning rate) of shape-tuning compared to center-tuning, since shape parameters outnumber center parameters.
We represent this as $c$, and when $c=0$, we recover the center-only variant of \gls{LRO}.
The motivating example showed both $c=0$ and $c=1$, with effects that depend on the model class.
Next, we have the hyperparameters given in Algorithms~\ref{alg:learning_inner} and~\ref{alg:learning}, which we list below in Table~\ref{tab:default} with suggested values and ranges.}
\RC{These suggestions are obtained through exploratory testing on multiple examples.}
\RC{For any experiment, we recommend prioritizing the selection of a model class and $c$, then tuning over the most impactful hyperparameters below, such as the radius $\rho$, the risk $\eta$, and the objective multiplier $\gamma$, either through grid search or bisection.} 
\begin{table}[ht]
\RC{\footnotesize
  \centering
	\begin{tabular}{lllllllllllll}\toprule
	$\rho$& $c$ & $\alpha^0$ & $\lambda^0$ & $\mu^0 $ & $\gamma$ &$\kappa$ & $\sigma$ & $\tau$ & $k_{\text{max}}$ & $t_{\text{max}}$ & $\eta$ & $|B|$\\
		\midrule
   $[10^{-5},10]$& $[0,1]$ &0  & $1$ & $1$  & $[0,5]$ & $ -10^{-5} $ & 1.005 & 0.95 & 10 & 15 & $(0,0.2]$& 30\\
	\bottomrule
		\end{tabular}
		\caption{Hyperparameter values and suggested ranges for the training algorithms.}
		\label{tab:default}
}\end{table}\RC{As noted in Remark~\ref{rem:conv},} to limit the training cost in practice, we cap the number of inner iterations~$t_{\rm max}$ and outer iterations~$k_{\rm max}$, as the per-iteration gains on the validation dataset are diminishing. 
For the step sizes $\delta^t$ given in Algorithm~\ref{alg:learning_inner}, we employ an incrementally decreasing schedule, where we multiply the step size by some value $\hat{\delta}$ every $\hat{k}$ steps.
\RC{We set them as $\hat{k} = 20$ and $\hat{\delta} = 0.5$, with an initial step size $\delta^1 = 0.001$, so that the updates decay fast enough to match our relatively limited iteration budget.
These numbers can be adjusted based on the application and the overall tuning budget.}

\section{Examples}
\label{sec:examples}
We run the following experiments on \ifpreprint the Princeton Institute for Computational Science and Engineering (PICSciE) facility \fi \RC{with up to 20 cores of 2.4 GHz AMD EPYC 9654 processors}.
We solve the inner optimization problems~\eqref{eq:robust_prob} with the Clarabel optimizer with default settings.
The learning procedure relies on the Python package cvxpylayers~\citep{cvxpylayers2019}, and
the code to replicate the experiments (along with baselines) can be accessed at
\begin{center}
	\githubrepo
\end{center}
\RC{For all experiments, we split the data into three components: training ($30\%$), validation ($20\%$), and testing ($50\%$), each independently sampled with continuous, essentially unique contexts.
We tune all hyperparameters based on the performance on the validation dataset.}
\RC{Specifically, filtering on a validation probability of certificate violation below 10\%, which describes the probability that the instance is infeasible and/or the in-sample objective underestimates the out-of-sample cost, we keep the three configurations with the lowest 90th percentile validation value (objective value on the feasible instances of the validation dataset), and among them choose the one with the lowest average validation value.}
\RC{We report the corresponding test-set values, from feasible test-set instances.}

\RC{We test the following methods.
We omit the generic DFL variants tested in the motivating example of Section~\ref{sec:intro_sets}, which have high percentages of infeasibilities, and report them only in \appxdfl.
Every method but two (the Wasserstein baselines, whose sets are balls around
training draws) constructs an ellipsoidal uncertainty set
\begin{equation}
	\label{eq:ell_set}
	\uncset(x) = \{u = b(x) + A(x)v \mid \|v\|_2 \leq \rho\},
\end{equation}
with differing center and shape maps $b(\cdot)$ and $A(\cdot)$.
We primarily use two model classes for the maps: linear, or multilayer perceptrons (MLPs) with one hidden layer of ReLU units. 
Below, we detail each method along with its tuning procedure.}
\RC{Every baseline with a $k$NN estimate (MSE, CMV, LKNN, and W-KNN) tunes its neighbor count $k$ jointly with $\rho$ on the validation data, over $5\%$, $10\%$, $15\%$, $20\%$, $25\%$, and $30\%$ of the training samples, while \gls{LRO} keeps its pretraining fixed at the $10\%$ value.}
\begin{itemize}
\item {\bf Learning for decision-focused \gls{RO} (LRO) (our method)}.
\RC{Both center and shape maps are trained for the downstream decisions, and we write $\uncset(\theta,x)$ for the set.
In the linear class, the maps are
\begin{equation}
\label{eq:lin_map}
    A^{\rm lro}(x) = W_Ax+ h_A,\quad b^{\rm lro}(x) = W_bx+h_b,
\end{equation}
where $W_A \in \reals^{m \times m \times p}$ is a tensor, with $(W_Ax)_{ij} = \sum_{k=1}^p (W_A)_{ijk}x_k$, so that $W_Ax \in \reals^{m \times m}$, and $h_A \in \reals^{m \times m}$, $W_b \in \reals^{m \times p}$, and $h_b \in \reals^m$.
In the MLP class, each map is its own network on the context, with one hidden layer of ReLU units, 16 for the portfolio example, and 32 for the inventory and motivating newsvendor examples. 
The parameters of the uncertainty set are $\theta$, the collected weights of the two maps of the chosen class, which we learn through the procedure given in Section~\ref{sec:training_procedure}. 
 In support of the regularity required by Assumption~\ref{ass:unique}, during training we add a quadratic penalty $\varepsilon\|z\|^2$, with $\varepsilon = 10^{-4}$, where $z$ represents the decision variables excluding the epigraph.
 For initialization, we pretrain the maps of either class on the training data $(U^N,X^N)$, detailed below in the MSE method. 
 For each model class, we test $c \in \{0,0.03,0.1,1\}$.
 For each value of $c$, we set $\eta = 0.1$, and jointly tune over $\gamma \in \{0.2,1,5\}$ and $\rho$, using the validation dataset to select among them.
Candidate $\rho$ are probed by bisection on the validation violation rate, with target 0.1, to a tolerance of 0.01.
If no probed radius meets the target while remaining solvable, the search returns the solvable radius with the lowest validation violation.
For deployment, we select the single best non-zero $c$ (by the same validation-set selection rule), and report it together with the {\it center-only} \gls{LRO} variant ($c=0$), which isolates whether the center map provides the gain. 
The remaining hyperparameters are set at the values in Table~\ref{tab:default}.
In \appxablation{}, we show ablations over $c$, $\eta$, and $\gamma$.}\RC{\item {\bf Mean-squared-error \gls{RO} (MSE)}.
This is the \gls{ETO} set that \gls{LRO} initializes from.
The parameters $\theta^0$ minimize the mean-squared error of the mean and covariance mappings on the training data $(U^N,X^N)$: the center map against the uncertainty realizations, and the shape map against the local covariance square roots $\hat{\Sigma}_k^{1/2}(X^N)$ of its residuals, defined below in the CMV method.
That is, we solve the two problems,
$\text{minimize}~\|b^{\rm lro}(X^N) - U^N\|^2$ and $\text{minimize}~\|A^{\rm lro}(X^N) - \hat{\Sigma}_k^{1/2}(X^N)\|^2$,
in closed form for the linear class and by gradient descent for the MLP class, with hyperparameters set through a priori (data-only) testing.
We tune the neighbor count $k$ (percentages given above) of the shape target jointly with $\rho$, which is probed by the same bisection as \gls{LRO}.}
\item{\bf Contextual Mean-Variance \gls{RO} (CMV)}.
\RC{We also apply an \gls{ETO} set inspired by~\cite{knn}, differing from MSE only in the shape map, which is local rather than a global fit.
The center $b^{\rm cmv}(x)$ is the same least-squares mean map of either model class, and
$A^{\rm cmv}(x) = \hat{\Sigma}_k^{1/2}(x),$
where, for each context $x$, we find the $k$ realizations in $X^N$ that are closest (in the Euclidean distance) to $x$, and take the square root of the covariance matrix of the corresponding $k$ center residuals $u^j - b^{\rm cmv}(x^j)$, so that the shape describes the spread about the fitted center.
We tune $k$ jointly with $\rho$, as for MSE.}\RC{\item {\bf Purely local $k$-nearest-neighbor (LKNN)}.
Here, both the center $\hat{\mu}_k(x)$ and shape $\hat{\Sigma}_k^{1/2}(x)$ are estimated directly from the $k$ nearest training neighbors of $x$.
We tune $k$ jointly with $\rho$, as for MSE.
When $k \leq m$, the local covariance is singular, and we add the ridge $0.01\,(\Tr \hat{\Sigma}_k(x)/m)\,I$ to it.}\RC{\item {\bf Oracle \gls{RO} (Oracle)}.
We also report the \gls{ETO} set built from the known conditional moments of the data-generating process, $A(x) = \Sigma(x)^{1/2}$ and $b(x) = \mu(x)$, with $\rho$ tuned on the validation data exactly as for \gls{LRO}.
These are the moments before clipping to a compact set, but the clipping is designed so that the effect on the moments is negligible.
We additionally report \gls{LRO} initialized at these moments and trained from there (\gls{LRO} (oracle)), which separates what decision-focused training adds from what estimation error costs.
For this variant, we select the best $c$ from $\{0,0.03,0.1,1\}$ on the validation data, so the center-only $c=0$ is also a candidate.
Neither set is generally available to a practitioner.}
\item {\bf Mean-Variance \gls{RO} (MV)}.
\RC{The standard non-contextual set~\citep{ben-tal_robust_2000}, where both maps are constant, given by the empirical mean and covariance of the training-set uncertainty $U^{N}$.
That is, $b(x) \equiv \hat{\mu}$ and $A(x) \equiv \hat{\Sigma}^{1/2}$.
We tune the set size $\rho$ on the validation data by the same bisection as \gls{LRO}.}\RC{\item {\bf Wasserstein \gls{DRO} (W)}.
The second non-contextual baseline is a Wasserstein-1 \gls{DRO} ambiguity set around the pooled training draws $U^N$~\citep{mohajerin_esfahani_data-driven_2018}. 
For tractability we use the \gls{MRO} approximation of~\cite{mro} (see \appxsets) with $K = 30$ clusters.
 We tune the set size $\rho$ on the validation data by the same bisection as \gls{LRO}.}\RC{\item {\bf Local Wasserstein \gls{DRO} (W-KNN)}.
The contextual counterpart of W and the Wasserstein analogue of LKNN, where for each context $x$, the ambiguity set is a Wasserstein-1 ball around the uncertainty realizations of the $k$ training contexts closest to $x$ in the Euclidean distance, again through the \gls{MRO} approximation, with $\min(10, k)$ clusters.
We tune $k$ jointly with the shared radius $\rho$, as for MSE. }\end{itemize}
For all methods, the hyperparameter tuning is done separately for each repetition of the experiments.
For each setting, we repeat the experiment 10 times, each with a new independent dataset, and report the performance metrics defined in Table~\ref{tab:perfmetrics}.
\RC{Note that the $\VaR_{90}$, $\CVaR_{90}$, and average-case values are computed only on the feasible test-set instances, where feasibility takes into account all uncertain constraints except for the epigraph-objective.} 

\begin{table}[ht]
\RC{\small
  \centering
\adjustbox{max width=\linewidth}{
	\begin{tabular}{ll}\toprule
		Name & Description\\
		\midrule
        $\hat{\eta}$ & Test-set empirical probability of certificate violation\\
        infeasibility & Test-set fraction of instances violating an uncertain constraint\\
        ~& other than the epigraph-objective by more than $10^{-5}$\\
		  $\VaR_{90}$ & Test-set 90th percentile objective value\\
        $\CVaR_{90}$ & Test-set 10\% conditional value-at-risk\\
        average-case & Test-set average objective value\\
        tune time & Time in seconds for training and hyperparameter tuning\\
				inference time & Average time in milliseconds per test-set context for solving the \\
				~& \gls{RO} problem, with the selected hyperparameters\\
	\bottomrule
		\end{tabular}}
		\caption{Performance metrics, averaged over 10 runs.}
		\label{tab:perfmetrics}
}\end{table}\subsection{Portfolio optimization}
 \label{sec:port}
	\RC{From Section~\ref{sec:gains}, we expect decision-focused learning to be more effective with multiple uncertain constraints.
We therefore apply it to a {\it constrained} portfolio example.  }
\paragraph{Problem description.}
\RC{We consider a robust risk-budgeting portfolio problem. 
We select a long-only portfolio $z \in \reals^n$ of $n$ stocks, grouped into $k$ sectors by the indicator matrix $S \in \{0,1\}^{n \times k}$. 
We maximize the worst-case return $u^Tz$, where $u \in \reals^n$ are the uncertain asset returns. 
In addition, we impose loss budgets, such that neither the loss of any sector nor the loss of any stock may exceed its budget, $\beta^{\rm sec} \in \reals^k_{++}$ and $\beta \in \reals^n_{++}$. 
Writing the objective in epigraph form with a single joint uncertain constraint, represented as a maximum of affine functions in $u$, gives the \gls{RO} problem
\begin{equation}
	\label{eq:port_1}
	\begin{array}{ll}
		\mbox{minimize} & L\\[2pt]
		\mbox{subject to} & \max\big\{g_0(z,L,u),\;
		\textstyle\max_{j \leq k} g^{\rm sec}_j(z,u),\;
		\textstyle\max_{i \leq n} g_i(z,u)\big\} \leq 0
		\quad \forall u \in \uncset(\theta,x)\\[4pt]
		& \ones^Tz = 1, \quad z \geq 0,
	\end{array}
\end{equation}
with variables $z \in \reals^n$ and $L \in \reals$, where the epigraph and risk-budget rows are
$$g_0(z,L,u) = -u^Tz - L, \qquad
g^{\rm sec}_j(z,u) = -u^T(S_j \odot z) - \beta^{\rm sec}_j, \qquad
g_i(z,u) = -u_iz_i - \beta_i,$$
with $\odot$ the elementwise product and $S_j$ the $j$th column of $S$. 
We intersect every
uncertainty set with the support $\{u \geq -\ones\}$, since a gross return is nonnegative.}

\paragraph{Data.}
\RC{We consider $n = 10$ stocks in $k = 3$ sectors and $n = 30$ stocks in $k = 5$ sectors, with $m = n$.
Every sample draws its own context $x = (x_{\rm sec}, x_{\rm vol}, x_{\rm mom})$ uniformly on $[-1,1]^{k+2}$, so $p = k+2$.
These consist of a state per sector measuring how favorable current conditions are, a market volatility state, and a momentum state.
Conditional on the context, the log returns $\log(\ones + u)$ are Gaussian, and the arithmetic returns $u$ are clipped to $[-1,3]$.
Stocks in favored sectors earn higher expected returns, with a spillover effect to neighboring sectors, and momentum lifts all stocks together.
Returns are correlated through a common market shock whose strength grows with the volatility, on top of stock-specific noise that is larger in favored sectors.
\appxdataport{} gives the full data-generating process.
Both the mean and covariance are contextual, with a partially non-affine mapping.
The loss budgets are set so that about $3$ budget rows are
active at $n=10$ and about $6$ at $n=30$, so the problem is effectively multi-constrained.
We consider datasets of $1000$ and $2000$ samples, so that the training sets have $N = 300$ and $N = 600$ samples.
}

\begin{table}[htp]
\RC{\centering
  \fontsize{10}{12}\selectfont\setlength{\tabcolsep}{4pt}
 \adjustbox{max width=\linewidth}{
	\begin{tabular}{llllllll}   \toprule
		\multicolumn{8}{l}{$n=10$, $N=300$}\\
		\midrule
		Method & $\hat{\eta}$ & infeas & $\VaR_{90}$ & $\CVaR_{90}$ & average-case & tune (s) & inf. (ms)\\
        \midrule
        \gls{LRO} (linear, $c=1$)  & $0.1230$ & $0.1176$ & $0.0569 \pm 0.0091$ & $\bf{0.0943 \pm 0.0077}$ & $\bf{-0.0897 \pm 0.0107}$ & $3\cdot (3\cdot 169)$ & $14.1$\\
        \gls{LRO} (linear, $c=0$)  & $0.1104$ & $0.1042$ & $\bf{0.0557 \pm 0.0089}$ & $0.0970 \pm 0.0072$ & $-0.0884 \pm 0.0054$ & $3\cdot 153$ & $15.6$\\
        MSE (linear)  & $0.1138$ & $0.1054$ & $0.0627 \pm 0.0048$ & $0.1043 \pm 0.0044$ & $-0.0732 \pm 0.0071$ & $6\cdot 40$ & $16.5$\\
        CMV (linear)  & $0.1222$ & $0.1124$ & $0.0635 \pm 0.0066$ & $0.1041 \pm 0.0111$ & $-0.0741 \pm 0.0076$ & $6\cdot 42$ & $16.6$\\
        \midrule
        \gls{LRO} (MLP, $c=1$)  & $0.1266$ & $0.1222$ & $0.0622 \pm 0.0069$ & $0.1019 \pm 0.0089$ & $-0.0858 \pm 0.0055$ & $3\cdot (3\cdot 285)$ & $14.5$\\
        \gls{LRO} (MLP, $c=0$)  & $0.1216$ & $0.1106$ & $0.0631 \pm 0.0073$ & $0.1043 \pm 0.0065$ & $-0.0818 \pm 0.0108$ & $3\cdot 278$ & $15.6$\\
        MSE (MLP)  & $0.1226$ & $0.1102$ & $0.0688 \pm 0.0094$ & $0.1102 \pm 0.0097$ & $-0.0690 \pm 0.0067$ & $6\cdot 38$ & $15.2$\\
        CMV (MLP)  & $0.1232$ & $0.1122$ & $0.0661 \pm 0.0051$ & $0.1065 \pm 0.0074$ & $-0.0705 \pm 0.0055$ & $6\cdot 41$ & $15.9$\\
        \midrule
        \gls{LRO} (oracle, $c=0$)  & $0.1090$ & $0.0980$ & $0.0588 \pm 0.0065$ & $0.1015 \pm 0.0070$ & $-0.0829 \pm 0.0059$ & $4\cdot (3\cdot 148)$ & $13.8$\\
        Oracle  & $0.1062$ & $0.0970$ & $0.0643 \pm 0.0081$ & $0.1037 \pm 0.0088$ & $-0.0710 \pm 0.0060$ & $1\cdot 38$ & $14.8$\\
        MV  & $0.1032$ & $0.0956$ & $0.0866 \pm 0.0045$ & $0.1140 \pm 0.0041$ & $-0.0316 \pm 0.0020$ & $1\cdot 41$ & $17.0$\\
        W  & $0.0756$ & $0.0398$ & $0.1126 \pm 0.0035$ & $0.1474 \pm 0.0031$ & $-0.0282 \pm 0.0036$ & $1\cdot 44$ & $7.0$\\
        W-KNN  & $0.1036$ & $0.0508$ & $0.0975 \pm 0.0026$ & $0.1366 \pm 0.0025$ & $-0.0444 \pm 0.0045$ & $6\cdot 207$ & $77.6$\\
        LKNN  & $0.1144$ & $0.1084$ & $0.0743 \pm 0.0058$ & $0.1066 \pm 0.0041$ & $-0.0543 \pm 0.0044$ & $6\cdot 38$ & $15.9$\\
		\bottomrule \\[-8pt]
		\toprule
		\multicolumn{8}{l}{$n=10$, $N=600$}\\
		\midrule
		Method & $\hat{\eta}$ & infeas & $\VaR_{90}$ & $\CVaR_{90}$ & average-case & tune (s) & inf. (ms)\\
        \midrule
        \gls{LRO} (linear, $c=1$)  & $0.1103$ & $0.1045$ & $\bf{0.0577 \pm 0.0078}$ & $\bf{0.0970 \pm 0.0091}$ & $\bf{-0.0907 \pm 0.0080}$ & $3\cdot (3\cdot 193)$ & $16.0$\\
        \gls{LRO} (linear, $c=0$)  & $0.1039$ & $0.0956$ & $0.0579 \pm 0.0072$ & $0.0981 \pm 0.0099$ & $-0.0860 \pm 0.0053$ & $3\cdot 185$ & $15.8$\\
        MSE (linear)  & $0.1078$ & $0.0994$ & $0.0645 \pm 0.0091$ & $0.1038 \pm 0.0111$ & $-0.0710 \pm 0.0068$ & $6\cdot 74$ & $16.0$\\
        CMV (linear)  & $0.1120$ & $0.1043$ & $0.0627 \pm 0.0084$ & $0.1010 \pm 0.0119$ & $-0.0735 \pm 0.0068$ & $6\cdot 74$ & $16.4$\\
        \midrule
        \gls{LRO} (MLP, $c=1$)  & $0.1054$ & $0.1022$ & $0.0618 \pm 0.0128$ & $0.1009 \pm 0.0138$ & $-0.0852 \pm 0.0069$ & $3\cdot (3\cdot 275)$ & $17.8$\\
        \gls{LRO} (MLP, $c=0$)  & $0.1084$ & $0.1022$ & $0.0628 \pm 0.0088$ & $0.1029 \pm 0.0121$ & $-0.0815 \pm 0.0115$ & $3\cdot 265$ & $16.7$\\
        MSE (MLP)  & $0.1133$ & $0.1050$ & $0.0644 \pm 0.0080$ & $0.1044 \pm 0.0112$ & $-0.0701 \pm 0.0058$ & $6\cdot 74$ & $15.8$\\
        CMV (MLP)  & $0.1118$ & $0.1040$ & $0.0662 \pm 0.0076$ & $0.1055 \pm 0.0083$ & $-0.0695 \pm 0.0066$ & $6\cdot 74$ & $16.1$\\
        \midrule
        \gls{LRO} (oracle, $c=0.03$)  & $0.1047$ & $0.0954$ & $0.0589 \pm 0.0062$ & $0.1003 \pm 0.0071$ & $-0.0825 \pm 0.0049$ & $4\cdot (3\cdot 188)$ & $16.2$\\
        Oracle  & $0.1008$ & $0.0925$ & $0.0648 \pm 0.0098$ & $0.1027 \pm 0.0091$ & $-0.0688 \pm 0.0054$ & $1\cdot 71$ & $14.8$\\
        MV  & $0.0991$ & $0.0944$ & $0.0835 \pm 0.0040$ & $0.1111 \pm 0.0053$ & $-0.0323 \pm 0.0019$ & $1\cdot 75$ & $15.5$\\
        W  & $0.0815$ & $0.0466$ & $0.1094 \pm 0.0050$ & $0.1470 \pm 0.0059$ & $-0.0281 \pm 0.0012$ & $1\cdot 49$ & $4.0$\\
        W-KNN  & $0.1090$ & $0.0495$ & $0.0964 \pm 0.0036$ & $0.1400 \pm 0.0058$ & $-0.0450 \pm 0.0046$ & $6\cdot 419$ & $80.8$\\
        LKNN  & $0.1014$ & $0.0979$ & $0.0720 \pm 0.0056$ & $0.1045 \pm 0.0069$ & $-0.0532 \pm 0.0024$ & $6\cdot 70$ & $15.8$\\
		\bottomrule
		\end{tabular}}
		\caption{Portfolio optimization with $n=10$ stocks, averaged over 10 runs.
See Table~\ref{tab:perfmetrics} for descriptions.
Values shown as mean $\pm$ half-IQR, where half-IQR $= (Q_{.75} - Q_{.25})/2$.
The best value for each category is in bold.
In tune, the base number is the average time for one hyperparameter setting, including training and the $\rho$ search, and the multiplier is the number of settings tuned over.
Inference is reported per test instance.}
		\label{table:port_n10}
}\end{table}

\begin{table}[ht]
\RC{\centering
  \fontsize{10}{12}\selectfont\setlength{\tabcolsep}{4pt}
 \adjustbox{max width=\linewidth}{
	\begin{tabular}{llllllll}   \toprule
		\multicolumn{8}{l}{$n=30$, $N=300$}\\
		\midrule
		Method & $\hat{\eta}$ & infeas & $\VaR_{90}$ & $\CVaR_{90}$ & average-case & tune (s) & inf. (ms)\\
        \midrule
        \gls{LRO} (linear, $c=0.1$)  & $0.1066$ & $0.1052$ & $0.0289 \pm 0.0112$ & $\bf{0.0671 \pm 0.0118}$ & $\bf{-0.1051 \pm 0.0159}$ & $3\cdot (3\cdot 863)$ & $82.8$\\
        \gls{LRO} (linear, $c=0$)  & $0.1062$ & $0.1026$ & $0.0314 \pm 0.0078$ & $0.0723 \pm 0.0076$ & $-0.0972 \pm 0.0123$ & $3\cdot 858$ & $73.6$\\
        MSE (linear)  & $0.1054$ & $0.1010$ & $0.0349 \pm 0.0073$ & $0.0734 \pm 0.0098$ & $-0.0907 \pm 0.0108$ & $6\cdot 200$ & $76.8$\\
        CMV (linear)  & $0.1088$ & $0.1054$ & $0.0357 \pm 0.0086$ & $0.0730 \pm 0.0118$ & $-0.0917 \pm 0.0111$ & $6\cdot 202$ & $79.3$\\
        \midrule
        \gls{LRO} (oracle, $c=0$)  & $0.1088$ & $0.1022$ & $\bf{0.0286 \pm 0.0064}$ & $0.0684 \pm 0.0082$ & $-0.1042 \pm 0.0026$ & $4\cdot (3\cdot 660)$ & $56.4$\\
        Oracle  & $0.0964$ & $0.0928$ & $0.0376 \pm 0.0047$ & $0.0753 \pm 0.0068$ & $-0.0877 \pm 0.0044$ & $1\cdot 164$ & $65.9$\\
        MV  & $0.0910$ & $0.0880$ & $0.0713 \pm 0.0056$ & $0.0989 \pm 0.0052$ & $-0.0302 \pm 0.0066$ & $1\cdot 200$ & $77.5$\\
        W  & $0.0922$ & $0.0356$ & $0.0964 \pm 0.0034$ & $0.1312 \pm 0.0115$ & $-0.0238 \pm 0.0048$ & $1\cdot 282$ & $50.5$\\
        W-KNN  & $0.1004$ & $0.0658$ & $0.0692 \pm 0.0094$ & $0.1087 \pm 0.0093$ & $-0.0547 \pm 0.0041$ & $6\cdot 742$ & $286.2$\\
        LKNN  & $0.1004$ & $0.0982$ & $0.0550 \pm 0.0022$ & $0.0873 \pm 0.0047$ & $-0.0568 \pm 0.0027$ & $6\cdot 218$ & $81.8$\\
		\bottomrule \\[-8pt]
		\toprule
		\multicolumn{8}{l}{$n=30$, $N=600$}\\
		\midrule
		Method & $\hat{\eta}$ & infeas & $\VaR_{90}$ & $\CVaR_{90}$ & average-case & tune (s) & inf. (ms)\\
        \midrule
        \gls{LRO} (linear, $c=0.1$)  & $0.1066$ & $0.1036$ & $\bf{0.0215 \pm 0.0038}$ & $\bf{0.0608 \pm 0.0029}$ & $\bf{-0.1154 \pm 0.0044}$ & $3\cdot (3\cdot 1043)$ & $81.3$\\
        \gls{LRO} (linear, $c=0$)  & $0.0986$ & $0.0955$ & $0.0283 \pm 0.0032$ & $0.0669 \pm 0.0045$ & $-0.1043 \pm 0.0035$ & $3\cdot 1012$ & $78.1$\\
        MSE (linear)  & $0.1033$ & $0.1002$ & $0.0325 \pm 0.0065$ & $0.0705 \pm 0.0037$ & $-0.0937 \pm 0.0040$ & $6\cdot 371$ & $80.6$\\
        CMV (linear)  & $0.1049$ & $0.1015$ & $0.0334 \pm 0.0052$ & $0.0716 \pm 0.0035$ & $-0.0939 \pm 0.0039$ & $6\cdot 376$ & $80.0$\\
        \midrule
        \gls{LRO} (oracle, $c=0$)  & $0.0981$ & $0.0911$ & $0.0288 \pm 0.0055$ & $0.0704 \pm 0.0042$ & $-0.1026 \pm 0.0047$ & $4\cdot (3\cdot 905)$ & $62.5$\\
        Oracle  & $0.1003$ & $0.0962$ & $0.0358 \pm 0.0034$ & $0.0749 \pm 0.0041$ & $-0.0885 \pm 0.0043$ & $1\cdot 307$ & $66.7$\\
        MV  & $0.1038$ & $0.1013$ & $0.0691 \pm 0.0038$ & $0.0960 \pm 0.0023$ & $-0.0310 \pm 0.0024$ & $1\cdot 374$ & $80.7$\\
        W  & $0.1064$ & $0.0507$ & $0.0893 \pm 0.0066$ & $0.1219 \pm 0.0067$ & $-0.0250 \pm 0.0029$ & $1\cdot 267$ & $25.4$\\
        W-KNN  & $0.1107$ & $0.0707$ & $0.0614 \pm 0.0045$ & $0.1018 \pm 0.0043$ & $-0.0597 \pm 0.0028$ & $6\cdot 1189$ & $228.2$\\
        LKNN  & $0.1005$ & $0.0982$ & $0.0522 \pm 0.0014$ & $0.0856 \pm 0.0031$ & $-0.0609 \pm 0.0016$ & $6\cdot 377$ & $81.5$\\
		\bottomrule
		\end{tabular}}
		\caption{Portfolio optimization with $n=30$ stocks, averaged over 10 runs.
Columns and conventions as in Table~\ref{table:port_n10}.}
		\label{table:port_n30}
}\end{table}

\paragraph{Results.}
\RC{Tables~\ref{table:port_n10} and~\ref{table:port_n30} compare all methods at a matched model class and the same validation target.
Across all four settings of $(n,N)$, \gls{LRO} wins on all three performance metrics against the \gls{ETO} methods (MSE, CMV, Oracle) of the same model class (linear or MLP or oracle).
The MLP results are generally suboptimal compared to the linear results, so we omit them for $n = 30$.
This demonstrates the mismatch between estimation accuracy and solution optimality, since the nonlinear contextual mapping is better represented by the MLP class, which nonetheless performs worse.
In fact, given the true mappings, the Oracle predictor is in general suboptimal compared to the MSE linear predictor.
The \gls{LRO} mapping initialized at the Oracle mapping outperforms the Oracle as well, indicating that our decision-focused gains are structural, and not a result of correcting estimation error. }\RC{As expected, on this contextual data, the non-contextual methods, MV and W, in general perform the worst.
The purely local LKNN outperforms them, and W-KNN does in the average case, but both are still suboptimal compared to the methods with global maps.
In particular, W-KNN is also expensive, with tuning times comparable to or even higher than those of \gls{LRO}, as it scales adversely with $m$ and $N$.
Per solve, a W-DRO problem is more expensive than an RO one with an ellipsoidal set.
Therefore, the inference times for all ellipsoidal sets (which are roughly equal) are lower than those of W-KNN.
The W tuning and inference times are both cheap, since in this problem the context does not enter the formulation, and a single solve per $\rho$ is sufficient.  }

\RC{Since \gls{LRO} has a fixed batch size, the number of inner solves per training step does not increase with $N$, but since the validation dataset grows in proportion to $N$, the overall tuning time of the methods typically increases with $N$.
Across the methods, performance generally improves with larger $N$, though in some cases with trade-offs between the tail and average case.
Comparing $n=30$ to $n=10$, the \gls{LRO} advantages hold, though the average-case gain reduces, for instance, from around 22\% (at $n=10, N=300$ for the linear \gls{LRO} vs.\ MSE) to around 16\% (at $n=30, N=300$ for the same pair), and from around 28\% to around 23\% on the $N=600$ pairs.
This is, however, partially due to the MSE average (denominator of the gain) increasing in magnitude.}

\RC{As the \gls{LRO} predictors from different model classes yield different results, we observe that our method is highly model-dependent.
If the tuning budget allows, we suggest trying several and selecting among them.
Of course, the Oracle is generally not available, and is shown only as a proof of concept.
The value of $c$ (magnitude of shape learning) is also an important design choice.
Here, the selected joint center and shape versions ($c>0$) improve on the center-only variants ($c=0$) in the average case and the $\CVaR$ in every setting, and most at $n=30$, though the win can be partially attributed to a larger realized $\hat{\eta}$.}
\RC{The \gls{LRO} oracle variant mostly selected $c=0$ as the best configuration, which shows that the importance of $c$ is model-based.}
\RC{In order to be relatively fair, we limited each \gls{LRO} setting to four values of $c$, $\eta = 0.1$, which is our $\hat{\eta}$ tuning target, and $\gamma \in \{0.2,1,5\}$, to test reasonable weights for the average-case objective.
In \appxhyper{}, we report what each hyperparameter setting contributes on its own, and additionally add several $(\eta,\gamma)$ pairs for ablation. 
}

\FloatBarrier
\subsection{Inventory management}
	\RC{From Section~\ref{sec:gains}, we also expect \gls{LRO} to outperform \gls{ETO} when the data is non-Gaussian.
In this example, we therefore test a multi-constrained non-Gaussian setting, and show that \gls{ETO} with the exact conditional mean and covariance could still be misspecified. }

\paragraph{Problem description.}
We consider a two-stage adjustable \gls{RO} problem on inventory management~\citep[Section 5.2.1]{pareto}.
There exists a retail network consisting of a single warehouse and $n$ different retail points, indexed by $i = 1, \dots, n$.
For simplicity, only a single product is sold.
There are a total of $C$ units of the product available, and each retail point starts with $0$ inventory and is capable of stocking at most $c_i$ units.
The transportation cost for distributing the items at the $i$th retail point is $t_i$ currency units per unit of inventory, and the operating cost is $h_i$ per unit.
The revenue per unit is $r_i$, which we take as a context parameter, as these values are known but may be subject to minor changes across problem instances, \RC{such as} supplier price changes. 
Customer demand $u \in \reals^n$ is uncertain, but we have access to contextual information \RC{$\xi \in [0,1]^{k+3}$, recording how active each of $k$ demand regimes currently is together with a promotion, a seasonality, and a competitor-price signal.
The total contextual information is thus $x = (r,\xi)$.}
The decision variables are the stocking decisions for all retail points, denoted $s \in \reals^n$.
In addition, the problem involves a vector $w(u)$ of realized sales at each point, which is dependent on the stocking decision $s$ and the unknown factors $u$ affecting the demand, \RC{given by}
$w_i(u) = \min\{s_i,u_i\} $
for each retail point.
The manager makes stocking decisions to maximize worst-case profits, which is equivalent to minimizing worst-case loss, 
$-r^Tw(u) + (t+h)^Ts$ for all $u\in \uncset(\theta,x)$.
For each problem instance, the \gls{RO} formulation is given as
\begin{equation}
	\label{eq:inv_1}
	\begin{array}{ll}
		\mbox{minimize} & \tau  \\
		\mbox{subject to} & -r^Tw(u) + (t+h)^Ts \le \tau, \quad \forall u\in \uncset(\theta,x)\\
		& w_i(u) \leq s_i, \quad i = 1, \dots, n, \quad \forall u\in \uncset(\theta,x)\\
		& w_i(u) \leq u_i, \quad i = 1, \dots, n, \quad \forall u\in \uncset(\theta,x)\\
			& \ones^Ts = C,\quad 0 \le s \leq c,\\
	\end{array}
\end{equation}
with variables $s\in \reals^n$, $w\in \reals^n$, and $\tau \in \reals$.
To get a tractable approximation, we replace the adjustable decisions $w_i(u)$ with their affine adjustable robust counterparts $v^0_i + v_i^Tu$, where $v^0_i \in \reals$ and $v_i \in \reals^m$ are auxiliary variables introduced as part of the linear decision rules~\citep{ben-tal_robust_2009}.
\ifpreprint The new formulation is
\begin{equation}
	\label{eq:inv_2}
	\begin{array}{ll}
		\mbox{minimize} & \tau + \varepsilon\|z\|^2\\
		\mbox{subject to} & -r^T\RC{v^0} - r^TVu + (t+h)^Ts \le \tau, \quad \forall u\in \uncset(\theta,x)\\
		& \RC{v^0_i} + v_i^Tu\leq s_i, \quad i = 1, \dots, n, \quad \forall u\in \uncset(\theta,x)\\
		& \RC{v^0_i} + v_i^Tu \leq u_i, \quad i = 1, \dots, n, \quad \forall u\in \uncset(\theta,x)\\
			& \ones^Ts = C, \quad 0 \le s \leq c,\\
	\end{array}
\end{equation}
where $V = [v_1 \cdots v_n]^T \in \reals^{n \times m}$. 
\fi
We treat the $2n+1$ uncertain constraints as a single joint uncertain constraint.

\paragraph{Data.}
\RC{As in the work by~\cite{pareto}, we let $n = 10$ retail points, which is also the size of the uncertain parameter, $m=n$.
Each sample draws its own context $x = (r, \xi)$, consisting of revenues $r$ perturbed around their nominal values, and additional information $\xi$ uniform on $[0,1]^{k+3}$, recording how active each of $k = 3$ demand regimes is, together with promotion, seasonality, and competitor-price signals.
The mean demand rises with the active regimes and the promotion, and falls with the revenue (which also indicates price).
The covariance is a regime-dependent diagonal plus a rank-one heavy-tailed seasonal common factor, tailored so that the heavy tail changes the conditional distribution but not its mean or covariance.
\appxdatainv{} gives the full data-generating process.
We use datasets of $1000$ samples, so that $N = 300$.}

\begin{table}[ht]
\RC{\centering
  \fontsize{10}{12}\selectfont\setlength{\tabcolsep}{4pt}
 \adjustbox{max width=\linewidth}{
	\begin{tabular}{llllllll}   \toprule
		\multicolumn{8}{l}{$n=10$, $N=300$}\\
		\midrule
		Method & $\hat{\eta}$ & infeas & $\VaR_{90}$ & $\CVaR_{90}$ & average-case & tune (s) & inf. (ms)\\
        \midrule
        \gls{LRO} (linear, $c=1$)  & $0.0972$ & $0.0972$ & $-302.10 \pm 6.32$ & $\bf{-251.26 \pm 10.34}$ & $-430.78 \pm 6.29$ & $3\cdot (3\cdot 881)$ & $34.0$\\
        \gls{LRO} (linear, $c=0$)  & $0.0906$ & $0.0906$ & $-253.69 \pm 11.61$ & $-199.36 \pm 12.46$ & $-423.54 \pm 4.81$ & $3\cdot 886$ & $33.4$\\
        MSE (linear)  & $0.0938$ & $0.0938$ & $-236.55 \pm 6.37$ & $-183.21 \pm 6.29$ & $-417.52 \pm 1.79$ & $6\cdot 85$ & $33.5$\\
        CMV (linear)  & $0.1018$ & $0.1018$ & $-234.49 \pm 6.65$ & $-180.13 \pm 6.96$ & $-419.16 \pm 1.77$ & $6\cdot 87$ & $33.5$\\
        \midrule
        \gls{LRO} (MLP, $c=1$)  & $0.0992$ & $0.0992$ & $-308.18 \pm 5.69$ & $-239.79 \pm 11.21$ & $-452.80 \pm 5.65$ & $3\cdot (3\cdot 917)$ & $35.0$\\
        \gls{LRO} (MLP, $c=0$)  & $0.1032$ & $0.1032$ & $-312.76 \pm 6.42$ & $-246.60 \pm 5.98$ & $-448.96 \pm 4.02$ & $3\cdot 924$ & $32.5$\\
        MSE (MLP)  & $0.0850$ & $0.0850$ & $-264.33 \pm 5.64$ & $-200.55 \pm 4.27$ & $-449.49 \pm 1.41$ & $6\cdot 88$ & $35.6$\\
        CMV (MLP)  & $0.0920$ & $0.0920$ & $-261.71 \pm 4.32$ & $-198.19 \pm 3.65$ & $-448.94 \pm 1.30$ & $6\cdot 86$ & $33.6$\\
        \midrule
        \gls{LRO} (oracle, $c=1$)  & $0.0964$ & $0.0964$ & $\bf{-315.93 \pm 8.95}$ & $-250.15 \pm 6.80$ & $\bf{-466.12 \pm 3.83}$ & $4\cdot (3\cdot 727)$ & $32.6$\\
        Oracle  & $0.0896$ & $0.0896$ & $-275.37 \pm 9.54$ & $-207.12 \pm 2.80$ & $-456.83 \pm 2.75$ & $1\cdot 91$ & $33.6$\\
        MV  & $0.1004$ & $0.1004$ & $-165.35 \pm 11.52$ & $-163.64 \pm 11.52$ & $-170.41 \pm 11.71$ & $1\cdot 87$ & $35.2$\\
        W  & $0.1030$ & $0.1026$ & $-197.35 \pm 8.55$ & $-194.56 \pm 8.70$ & $-205.01 \pm 8.43$ & $1\cdot 1055$ & $434.8$\\
        W-KNN  & $0.0994$ & $0.0982$ & $-257.23 \pm 11.89$ & $-226.98 \pm 13.04$ & $-340.36 \pm 16.35$ & $6\cdot 298$ & $118.1$\\
        LKNN  & $0.1024$ & $0.1024$ & $-195.46 \pm 14.70$ & $-155.79 \pm 15.82$ & $-374.32 \pm 13.65$ & $6\cdot 97$ & $35.8$\\
		\bottomrule
		\end{tabular}}
		\caption{Average inventory management performance over 10 runs.
Columns and conventions as in Table~\ref{table:port_n10}.}
		\label{table1}
}\end{table}

\paragraph{Results.}
\RC{Table~\ref{table1} compares all methods.
Again, across all model classes, \gls{LRO} comes out ahead.
For the linear class, \gls{LRO} outperforms all other linear methods, including its center-only variant, on all three performance metrics, and achieves the best $\CVaR$ overall.
Contrary to the portfolio experiment, we observe that shape training contributes most of the gain, as the center-only variant recovers only about a quarter of the $\CVaR$ improvement over MSE.
For the MLP class, \gls{LRO} outperforms MSE in $\VaR$ and $\CVaR$ (by $17\%$ and $20\%$, respectively), while matching it in the average case, at a slightly higher realized violation rate ($0.099$ against $0.085$).
Here, the center-only variant and the selected shape learning ($c=1$) perform similarly, with one slightly better in the tails and the other in the average case.}
\RC{For this problem, Oracle outperforms the remaining \gls{ETO} baselines, and the MLP \gls{ETO} baselines outperform their linear counterparts, which shows that good moment estimation can still be useful.
However, the \gls{LRO} method initialized from the oracle wins against it in all three metrics, so decision-focused learning exploits structural gains that the oracle cannot see.
\gls{LRO}-oracle also outperforms \gls{LRO}-MLP, which once again demonstrates the method's model-dependence.}

\RC{Similar to the portfolio example, the non-contextual methods MV and W perform badly here. 
The local W-KNN method achieves a good $\CVaR$, but has a suboptimal average-case cost.
LKNN does better than W-KNN in the average case, but is suboptimal in the tails.
Comparing MSE, CMV, and LKNN, we note that global maps are generally beneficial, and in particular the center map (going from LKNN to CMV) has the largest effect.
None of these methods, however, does as well as \gls{LRO}, which illustrates the benefits of decision-focused learning.
In \appxhyper{}, we show additional ablations for \gls{LRO}.
The consensus there, across both examples, is that no single $(\eta,\gamma)$ is best everywhere, but every setting improves on MSE in both tail metrics, so a tuning budget is better spent on the model class and $c$, which in general have larger effects.
}

\FloatBarrier
\section{\RC{Conclusions and discussion}}
We have developed a data-driven method to learn the contextual shape of uncertainty sets in a decision-focused manner, \RC{that is,} using the structure and the solution of the robust counterpart.
\RC{We cast learning as a stochastic optimization problem that balances average-case and worst-case performance, while a $\CVaR$ constraint limits how often the learned decisions violate the uncertain constraints.}
We \RC{characterize the accumulation points} of our algorithm using nonsmooth implicit differentiation, and provide finite-sample probabilistic guarantees on unseen data.
\RC{We also characterize when the decision-focused set is worth learning.
With Gaussian data, a single uncertain constraint, and a purely worst-case objective, the \gls{ETO} set built from the true conditional moments is already optimal within the ellipsoidal class, and no reshaping can enlarge the feasible region at the same guarantee.
Decision-focused learning requires the setting to depart from these conditions, that is, a non-Gaussian conditional distribution, several uncertain constraints that share the violation budget, an objective that weighs average-case against worst-case performance, or a guarantee that is joint instead of conditional.}
\RC{This motivates our algorithmic and experimental design.}
Through numerical examples under these settings, we demonstrate the performance gain achieved by our method compared to traditional methods, for both worst-case and average-case performance.
\RC{In particular, initializing \gls{LRO} at the oracle mapping and improving upon it suggests that the gains are structural rather than only a correction of estimation error.}

\RC{Several limitations remain, however, which are directions for future research.
First, our experiments are run with synthetic data, with hypothesized contextual mappings.
An ideal next step is to deploy this on a real-world dataset, which may have noisier signals. 
Furthermore, our algorithm is analyzed in the ideal asymptotic state, but applied under finite iterations.
Given the complexity of the formulation, convergence is difficult to prove, and even currently it is presented for accumulation points and not in the actual limit.
Another limitation is the lack of a global optimality guarantee.
The single-level formulation available to us, in \appxdirect, is nonconvex even in simple settings, with complementarity constraints and terms that are bilinear in $\theta$, and is difficult to solve globally.} 

\RC{Along these lines, it would be insightful to further analyze the geometric properties of the uncertainty sets found.
We have done so for a simple case in Section~\ref{sec:gains}, which is possible due to the independence of the constraints.
Identifying such structural properties more generally could yield sets in closed form, or at least guide the hyperparameter search of our algorithm. 
Nevertheless, this work demonstrates that better solutions exist and can be found, and, on the examples we study, the uncertainty sets we reach already outperform \gls{ETO} baselines.}

\newcommand{\myack}{Irina Wang and Bartolomeo Stellato are supported by
the NSF CAREER Award ECCS 2239771. 
Irina Wang is also supported by the Princeton Wallace Memorial Fellowship.
The simulations presented in this article were performed on computational resources managed and supported by Princeton Research Computing, a consortium of groups including the Princeton Institute for Computational Science and Engineering (PICSciE) and the Office of Information Technology's High Performance Computing Center and Visualization Laboratory at Princeton University.
We would like to thank Jason Klusowski for the useful feedback and discussions on the probabilistic guarantees, and Cole Becker for his contributions to the original ideation of this work. 
}

\ifpreprint
\section{Acknowledgments}
\myack
\fi

\begin{appendices}

\ifpreprint
\providecommand{\ecref}[1]{\ifpreprint Appendix~\ref{#1}\else \ref{#1}\fi}\section{Uncertainty sets and general reformulations}
\label{appx:gen_sets}
We derive convex reformulations of~\eqref{eq:robust_prob} for a general class of constraints $g$ and common uncertainty sets.
\RC{Throughout this appendix, $f(z)$ denotes a generic deterministic objective, since the uncertain objective of~\eqref{eq:robust_prob} enters through its epigraph constraint as in Section~\ref{sec:joint}.} 
We express $u$ as an affine transformation of an uncertain parameter $v$, \ie, $u = A(x)v + b(x)$, where $A(x) \in \reals^{m\times \hat{m}}$, $\hat{m} \leq m$, and $b(x) \in \reals^m$.
The proofs are in
\ifpreprint
Appendix~\ref{appx:reform}.
\else
\RC{the online appendices that follow}.
\fi
For notational simplicity, we suppress the dependency of $g$ on $x$ here and for all following sections on reformulations.
We consider two types of functions $g$, the maximum of concave functions (left), where each $g_l$ is concave in $u$, and its special case, the maximum of affine functions (right), where \RC{$P_l$ is affine and $a_l$ is convex in $z$, so that each piece is convex in $z$ for every $u$},
\[
g(z,u) = \max_{l=1,\dots,L} g_l(z,u), \qquad g(z,u) = \max_{l=1,\dots,L} (\RC{P_l(z)^T}u + a_l(z)).
\]

\paragraph{Box and ellipsoidal uncertainty.}
The ellipsoidal uncertainty set is defined as
\[
\uncset(\theta,x) = \left\{u = A(x)v+b(x) \mid \RC{\| v\|} \leq 1,\quad u \in S\right\},
\]
where $A(x)\in \reals^{m\times \hat{m}}$, $b(x)\in \reals^{m}$, and $S$ includes additional support information such as upper and lower bounds.
\RC{Here $\|\cdot\|$ is any norm, and $\|\cdot\|_*$ denotes its dual norm.
Taking $\|\cdot\| = \|\cdot\|_\infty$ gives the box uncertainty set.} The trainable parameter is $\theta = (A(\cdot),b(\cdot))$, the mapping from the context parameter to the reshaping matrices. Note that if $A = \rho I$, $b = 0$, and~$S = \reals^m$, this is equivalent to the traditional ellipsoidal uncertainty set ${ \{ u \mid\RC{\| u\|_2} \leq \rho \}}$ \RC{when $\|\cdot\| = \|\cdot\|_2$}. 
For the maximum of concave constraint, we have the following reformulation,
\begin{equation*}
	\begin{array}{ll}
	\text{minimize} & f(z)\\
	\mbox{subject to} & [-g_l]^\star(z,h_l)+ \sigma_S(\omega_l)  -b(x)^T(h_l + \omega_l) \\ & \qquad + \RC{\|A(x)^T(h_l + \omega_l) \|_*} \leq 0, \quad l = 1,\dots,L,
	\end{array}
	\end{equation*}
with auxiliary variables $\omega_l, h_l \in \reals^m$ for $l = 1,\dots,L$, and $g^\star(z,h)$ the conjugate function.
This result follows from the theorem of infimal convolutions and composites for conjugate functions~\citep{rockafellar_wets_1998,composites}.
In the special (maximum of affine) case, we instead have
\begin{equation*}
	\begin{array}{ll}
		\mbox{minimize} & f(z)\\
		\mbox{subject to} & a_l(z) + \sigma_S(\omega_l)+ b(x)^T(P_l(z) - \omega_l) \\ & \qquad +  \RC{\|A(x)^T(\omega_l -P_l(z))\|_*} \le 0, \quad l = 1,\dots,L.	\end{array}
\end{equation*}

\paragraph{Budget uncertainty.}
The budget uncertainty set requires both box and 1-norm constraints,
\[
\uncset_{\text{bud}}(\theta,x) = \{u = A(x)v+b(x) \mid \| v\|_\infty \leq \rho_1, \quad \| v \|_1 \leq \rho_2, \quad u \in S\},
\]
where $\theta = (A(\cdot),b(\cdot))$.
For the maximum of concave constraint, we have the following reformulation,
\begin{equation*}
	\begin{array}{ll}
		{\text{minimize}} & f(z) \\
		\mbox{subject to} & [-g_l]^\star(z,h_l)  + \sigma_S(\omega_l) -  b(x)^T(h_l + \omega_l) \\ & \qquad + \rho_1\|r_l - A(x)^T(h_l + \omega_l)\|_1  + \rho_2\|r_l\|_\infty \leq 0, \quad l =  1, \dots,L,	\end{array}
	\end{equation*}
with $r_l \in \reals^{\hat{m}}$, $\omega_l, h_l \in \reals^m$ for $l = 1,\dots,L$ as auxiliary variables.
The maximum of affine reformulation is
\begin{equation}
	\label{eq:budget_max}
	\begin{array}{ll}
		\mbox{minimize} & f(z)\\
		\mbox{subject to} &  a_l(z) + \sigma_S(\omega_l) + b(x)^T(P_l(z) - \omega_l) \\ & \qquad + \rho_1\RC{\|r_l - A(x)^T(\omega_l - P_l(z))\|_1} + \rho_2\|r_l\|_\infty \leq 0, \quad l = 1,\dots,L.
	\end{array}
\end{equation}

\paragraph{Polyhedral uncertainty.}
The polyhedral uncertainty set is defined as
\begin{equation*}
    \uncset(\theta,x) = \{ u = A(x)v+b(x) \mid Dv \leq d, \quad u \in S\},\end{equation*}
where $\theta = (A(\cdot),b(\cdot))$.
For the maximum of concave constraint, we have the following reformulation,
\begin{equation*}
\begin{array}{ll}
	{\text{minimize}} & f(z)\\
	\mbox{subject to} & [-g_l]^\star(z,\omega_l) \RC{{}+ d^Th_l} - b(x)^T(\omega_l + \gamma_l) + \sigma_S(\gamma_l) \leq 0, \quad l = 1,\dots,L\\
	 &A(x)^T(\omega_l+ \gamma_l) = -D^Th_l,\quad l = 1,\dots,L\\
	 & \RC{h_l \geq 0, \quad l = 1,\dots,L},
\end{array}
\end{equation*}
with auxiliary variables $\gamma_l, \omega_l \in \reals^m$, $h_l \in \reals^{r}$ for $l = 1,\dots,L$.
The maximum of affine reformulation is
\begin{equation*}
	\begin{array}{ll}
		\mbox{minimize} & f(z)\\
		\mbox{subject to} &  a_l(z) + d^Th_l \RC{{}+ \sigma_S(\gamma_l)} + b(x)^T( P_l(z)-\gamma_l)\le 0, \quad l = 1,\dots,L\\
		& A(x)^T(P_l(z)-\gamma_l) = D^Th_l , \quad l = 1,\dots,L\\
		& h_l \geq 0, \quad l = 1,\dots,L.
	\end{array}
\end{equation*}

\paragraph{Mean robust uncertainty.}
The \RC{\gls{MRO}}~\citep{mro} set
is defined as
\begin{equation*}
	\begin{aligned}
		\uncset_{\text{mro}}(\theta,x) & = \Big\{ u = A(x)(v_{11},\dots,v_{LK})+ b(x)\in \supp^{K\times L} \;\Big|\; \\
& \qquad\qquad \sum_{k=1}^K \sum_{l=1}^L \alpha_{lk} \RC{\| v_{lk} - \bar{d}_k \|^r \le \rho^r},\quad \alpha \in \Gamma \Big\},
	\end{aligned}
\end{equation*}
where $\Gamma = \{\alpha \mid \sum_{l=1}^L\alpha_{lk}= w_k,~\alpha_{lk}\geq 0 \text{ for all } k,l\}$.
Here, we assume the constraint function $g$ of the \gls{RO} problem to be a maximum of $L$ pieces, indexed by $l$.
\RC{For this set, the robust constraint is the worst-case weighted expectation $\sum_{k=1}^K\sum_{l=1}^L \alpha_{lk}\, g_l(z, A(x)v_{lk}+b(x)) \le 0$ over all $(v,\alpha)$ in the set, with each $v_{lk} \in \reals^m$, rather than a single $u$ substituted into~\eqref{eq:robust_prob}.} 
The support $S$ of each component of $u$ is assumed to be convex and closed. 
We assume $A$ is square, with $A(x) \in \reals^{m\times m}$\RC{, that is, $\hat{m} = m$}, and we have $\theta = (A(\cdot), b(\cdot))$. 
This set depends on $U^N$, a dataset of $N$ samples.
We partition $U^N$ into $K$ disjoint subsets, with $\bar{d}_k$ the centroid of the $k$th subset. The weight we place upon each subset is $w_k$, and is equivalent to the proportion of points in the subset. \RC{We choose $r \geq 1$ to be an integer exponent, and $\rho$ to be a positive radius.
We let $r_*$ be the conjugate exponent, with $1/r + 1/r_* = 1$ and $r_* = \infty$ when $r = 1$, and define $\phi(r_*) = (r_*-1)^{r_*-1}/r_*^{r_*}$.}
For the maximum of concave constraint, we have the following reformulation,
\begin{equation*}
	\begin{aligned}
	\begin{array}{ll}
	{\text{minimize}} & f(z)\\
	\mbox{subject to} &\lambda \RC{\rho^r} + \sum_{k=1}^K w_k s_k \leq 0\\
	 & [-g_l]^\star(\RC{z},h_{lk}) -b(x)^T(h_{lk} + \omega_{lk}) + \sigma_S(\omega_{lk}) - (A(x)^T(h_{lk} + \omega_{lk}))^T \bar{d}_k \\
& \qquad \qquad + \RC{\phi(r_*)\lambda \left\|A(x)^T(h_{lk} + \omega_{lk})/\lambda \right\|^{r_*}_*} \le s_k,\\ & \qquad\qquad l = 1,\dots,L, \quad k = 1,\dots,K\\
& \lambda \geq 0,
	\end{array}
	\end{aligned}
	\end{equation*}
with auxiliary variables \RC{$s \in \reals^K$ and} $h_{lk}, \omega_{lk} \in \reals^m$ for $l = 1,\dots,L$, $k = 1,\dots, K$, and $\lambda \in \reals$.
For the special case of the maximum of affine constraint, we have
\begin{equation}
	\label{eq:mro_max}
	\begin{aligned}
		\begin{array}{ll}
			\text{minimize} & f(z)\\
			\mbox{subject to} & \lambda\RC{\rho^r}  + \sum_{k=1}^{K} w_{k} s_{k}  \le 0\\
			& a_l(z) + \sigma_S(\omega_{lk}) - b(x)^T(\omega_{lk} - \RC{P_l(z)}) - (A(x)^T(\omega_{lk} - P_l(z)))^T\bar{d}_k \\
   &\qquad \qquad + \RC{\phi(r_*)\lambda \left\|A(x)^T(\omega_{lk} - P_l(z))/\lambda \right\|^{r_*}_*} \leq s_k,\\ & \qquad\qquad l = 1,\dots, L, \quad k = 1,\dots, K\\
			& \lambda \geq 0.
		\end{array}
	\end{aligned}
\end{equation}
\RC{Here, $\lambda\|y/\lambda\|_*^{r_*}$ denotes the closed perspective, which equals $0$ if $\lambda = 0$ and $y = 0$, and $+\infty$ if $\lambda = 0$ and $y \neq 0$.
Our experiments use $r = 1$, for which $r_* = \infty$ and $\phi(r_*)\lambda\|y/\lambda\|_*^{r_*}$ is the indicator of $\|y\|_* \le \lambda$~\citep[Section 2]{mro}.
In both reformulations, $\lambda\rho^r$ then becomes $\lambda\rho$, and the term in $\phi(r_*)$ is replaced by the constraint $\|A(x)^T(h_{lk} + \omega_{lk})\|_* \le \lambda$, respectively $\|A(x)^T(\omega_{lk} - P_l(z))\|_* \le \lambda$, for $l = 1,\dots,L$ and $k = 1,\dots,K$.}

 \fi

\section{\RC{Conic reformulation and regularity conditions}}
\label{sec:dpp}
The convex optimization problems~\eqref{eq:robust_prob} can be transformed to the cone program~\citep{cvxpylayers2019}
	\begin{equation}
 \label{eq:conic}
		\begin{array}{ll}
			\mbox{minimize} & c^Tz\\
			\mbox{subject to} & Gz+ s = h \\
			& s \in \mathcal{K},
		\end{array}
	\end{equation}
where $z \in \reals^{\bar{n}}$ is the primal variable, $s \in \reals^{\bar{m}}$ is a slack variable, and the set $\mathcal{K} \subseteq \reals^{\bar{m}}$ is a nonempty, closed, convex cone with dual cone $\mathcal{K}^\star \subseteq \reals^{\bar{m}}$.
\RC{For the maximum of affine constraints with $\ell_1$, $\ell_2$, $\ell_\infty$, or budget uncertainty sets, polyhedral supports, and \gls{MRO} sets with $r \in \{1,2\}$ (\appxsets), which cover our experiments, the cone $\mathcal{K}$ is a Cartesian product of zero cones, nonnegative orthants, and second-order cones.}
The dual problem is given as
	\begin{equation}
 \label{eq:conic_dual}
		\begin{array}{ll}
			\mbox{minimize} & h^Tr\\
			\mbox{subject to} & G^Tr+ c = 0 \\
			& r \in \mathcal{K}^\star,
		\end{array}
	\end{equation}
	where $r \in \reals^{\bar{m}}$ is the dual variable.
\RC{The data $(G,h,c)$ come from canonicalizing~\eqref{eq:robust_prob}, whose constituent functions $g_l$ we require to be concave in the uncertain parameter and to satisfy the disciplined parameterized programming (DPP) rules of~\citet{cvxpylayers2019}.
Canonicalization first collects the declared parameter values, the shape $A(x)$, the center $b(x)$, the radius $\rho$, and any context-dependent problem data, into a vector $\pi(\theta,x)$, and then forms $(G,h,c)$ as an affine function of $\pi$~\citep[Lemma 1]{cvxpylayers2019}.
The map $\pi$ is not affine in $(\theta,x)$, since the linear maps~\eqref{eq:lin_map} of Section~\ref{sec:examples} are already bilinear in the weights and the context, and the MLP maps are not affine in their weights.
It is, however, locally Lipschitz and definable by Assumptions~\ref{ass:define} and~\ref{ass:lip}, which cover both classes.}

\RC{The primal-dual solutions of~\eqref{eq:conic} and~\eqref{eq:conic_dual} are the roots in $(z,r)$ of the residual}
\begin{equation}
\label{eq:kkt_residual}
R(z,r;\theta,x) = \begin{bmatrix} c + G^Tr \\ h - Gz - \Pi_{\mathcal{K}}(h - Gz - r)\end{bmatrix},
\end{equation}
\RC{where $\Pi_{\mathcal{K}}$ is the Euclidean projection onto $\mathcal{K}$, since $s = h - Gz$ satisfies $s \in \mathcal{K}$, $r \in \mathcal{K}^\star$, and $s^Tr = 0$ if and only if $s = \Pi_{\mathcal{K}}(s - r)$.
A conservative Jacobian of $R$ with respect to $(z,r)$ is the set of matrices}
\begin{equation}
\label{eq:kkt_jacobian}
U = \begin{bmatrix} 0 & G^T \\ -(I - P)G & P\end{bmatrix}, \qquad P \in J_{\Pi_{\mathcal{K}}}(h - Gz - r),
\end{equation}
\RC{where $J_{\Pi_{\mathcal{K}}}$ is the Clarke Jacobian of the projection, which is conservative since the projection is Lipschitz and semialgebraic~\citep[Proposition 2 and Corollary 1]{conservative}, and explicit for the cones above~\citep{cones,cones1}.
Appendix~\ref{appx:subdiff} differentiates the solution map through $R$, which requires every $U$ to be nonsingular at the solution.
Lemma~\ref{lem:nonsingular} derives this from the two conditions of Assumption~\ref{ass:unique}, which we now state at a primal-dual solution $(z^\star,s^\star,r^\star)$ of~\eqref{eq:conic}.}

\RC{\emph{Primal constraint nondegeneracy.}
The range of $G$ and the lineality space of the tangent cone to $\mathcal{K}$ at $s^\star$ together span $\reals^{\bar{m}}$~\citep[Section 4.6]{bonnans2000}.
It holds for all data $(G,h,c)$ outside a set of measure zero~\citep{pataki2001generic}.
However, genericity does not cover the structured family produced by the learned parameterization, which may lie entirely in a measure zero set, so we keep this condition as an explicit assumption.}

\RC{\emph{Strong second-order sufficient condition.}
The Hessian of the Lagrangian of~\eqref{eq:conic} in $z$, plus the curvature term contributed by $\mathcal{K}$ at $(s^\star,r^\star)$, is positive definite on the affine hull of the critical cone~\citep{bonnans2005}.
We expect the condition to hold when every variable of the reformulation is penalized, since the penalty then supplies curvature in every direction.
In practice, we add the penalty $\varepsilon\|z\|^2$ over the decision variables, so that the decision is unique and the second-order condition holds along the decision directions.
The epigraph variable needs no penalty, since at every solution it equals the worst-case objective and is thus determined by the decision variables.
However, the conic reformulation also carries auxiliary variables along which the penalty gives no curvature, and they are pinned only when their constraints are active, for instance when every robust constraint is active at the solution.
Penalizing them as well biases the objective, so we leave them unpenalized, and keep the condition as an assumption.}

\RC{\begin{lemma}[Nonsingularity of $U$]
\label{lem:nonsingular}
Suppose the cone $\mathcal{K}$ is a Cartesian product of zero cones, nonnegative orthants, and second-order cones.
Let a primal-dual solution of~\eqref{eq:conic} satisfy primal constraint nondegeneracy and the strong second-order sufficient condition.
Then every matrix $U$ in~\eqref{eq:kkt_jacobian} is nonsingular at that solution.
\end{lemma}
\begin{proof}
Since $J_{\Pi_{\mathcal{K}}}$ is the Clarke Jacobian of the projection, the matrices $U$ in~\eqref{eq:kkt_jacobian} are the elements of the Clarke Jacobian of $R$ with respect to $(z,r)$.
On the zero-cone blocks of $\mathcal{K}$, the projection vanishes and the corresponding rows of $R$ are affine, so once these blocks are eliminated, problem~\eqref{eq:conic} is a linear symmetric cone program in the sense of~\citet{kong2011jacobian}.
Specifically, it is their dual problem, with their variables $(x,y,s)$ equal to our $(r,z,h-Gz)$, and $R$ is the natural-map residual~\citep{facchinei2003} of its KKT system with the slack eliminated through $s = h - Gz$.
Both eliminations are invertible affine changes of variables, the first because primal nondegeneracy gives the equality rows of $G$ full row rank, and the second because the slack $s = h - Gz$ is an affine function of $z$.
They therefore map Clarke Jacobian elements one-to-one, preserve nonsingularity, and pass both conditions to the reduced program.
Primal constraint nondegeneracy makes the multiplier $r^\star$ unique, so by \citet[Proposition 3.6]{kong2011jacobian} and its dual counterpart, which the authors state holds likewise, the strong second-order sufficient condition is equivalent to dual constraint nondegeneracy.
Primal and dual constraint nondegeneracy together are equivalent to nonsingularity of every element of the Clarke Jacobian of the natural-map residual, and to strong regularity~\citep[Theorem 3.9]{kong2011jacobian}, so every $U$ in~\eqref{eq:kkt_jacobian} is nonsingular at the solution.
\end{proof}}
\section{Proof of conservative Jacobians for the robust optimization problem}
\label{appx:subdiff}
\RC{To prove Theorem~\ref{thm:subdiff}, we differentiate the primal-dual solution of~\eqref{eq:conic} as a root of the residual $R$ in~\eqref{eq:kkt_residual}, which depends on $(\theta,x)$ only through the parameter vector $\pi(\theta,x)$ of Appendix~\ref{sec:dpp}.
We first differentiate the root in $\pi$ with the following nonsmooth implicit function theorem, which combines Clarke's implicit function theorem~\citep[Theorem 7.1.1]{clarke1990} with Corollary 1 of~\citet{path-diff}.}
\RC{\begin{proposition}[Nonsmooth implicit function theorem]
\label{prop:implicit}
    Let $R \colon \reals^{d} \times \reals^{d_w} \to \reals^{d}$ be path differentiable with conservative Jacobian $J_R$, and let $R(\bar{y};\bar{w}) = 0$.
Suppose that the set $J_R(\bar{y};\bar{w})$ is convex and that every matrix $[U~V] \in J_R(\bar{y};\bar{w})$, where $U \in \reals^{d \times d}$ is the block corresponding to $y$, has $U$ invertible.
Then there exist neighborhoods $\mathcal{Y}$ of $\bar{y}$ and $\mathcal{V}$ of $\bar{w}$ and a locally Lipschitz, path-differentiable function $y \colon \mathcal{V} \to \mathcal{Y}$ such that $y(w)$ is the unique root of $R(\cdot;w)$ in $\mathcal{Y}$ for every $w \in \mathcal{V}$, with conservative Jacobian
    \[
    J_y(w) = \{-U^{-1}V \mid [U~V] \in J_R(y(w);w)\}.
    \]
\end{proposition}}
\RC{\begin{proof}
By convexity, $J_R$ contains the Clarke Jacobian of $R$ at $(\bar{y},\bar{w})$~\citep[Theorem 1]{path-diff}, so Clarke's theorem gives the locally Lipschitz unique root $y(\cdot)$ and Corollary 1 of~\citet{path-diff} gives its conservative Jacobian.
\end{proof}}
\RC{We apply the proposition with $y = (z,r)$ and $w = \pi$.}

\RC{\emph{Path differentiability and convexity.}
The map $R$ composes smooth (bilinear) maps of $((z,r),\pi)$ with the projection $\Pi_{\mathcal{K}}$, whose Clarke Jacobian is conservative (Appendix~\ref{sec:dpp}).
By the chain rule~\citep{conservative}, $R$ is path differentiable with conservative Jacobian~\eqref{eq:kkt_jacobian} in $(z,r)$.
In every $[U~V]$ in $J_R$, with $V$ the derivative in $\pi$ through the affine data $(G,h,c)$, both blocks are affine in the selection $P \in J_{\Pi_{\mathcal{K}}}(h - Gz - r)$.
The Clarke Jacobian is convex, so $J_R$ is convex.}

\RC{\emph{Existence and uniqueness.}
By Assumption~\ref{ass:unique} and Lemma~\ref{lem:nonsingular}, for each $(\theta,x) \in \mathcal{D}_\theta \times \mathcal{D}_x$, the residual at $\pi(\theta,x)$ has a root $(z^\star,r^\star)$ at which every $U$ is invertible, so Proposition~\ref{prop:implicit} gives a neighborhood of $(z^\star,r^\star)$ that contains no other root.
The roots are the primal-dual solutions of the convex program~\eqref{eq:conic}, which form a convex set, and a convex set with an isolated point is a single point.
Hence $(z^\star,r^\star)$ is the unique primal-dual solution.
Since $\pi$ is continuous, $\pi(\theta',x')$ lies in the neighborhood $\mathcal{V}$ of the proposition for every $(\theta',x') \in \mathcal{D}_\theta \times \mathcal{D}_x$ near $(\theta,x)$.
The root $y(\pi(\theta',x'))$ is then the unique primal-dual solution at $(\theta',x')$, by the same argument at $(\theta',x')$.
Near $(\theta,x)$, the root map therefore coincides with the solution map, $(z(\theta',x'),r(\theta',x')) = y(\pi(\theta',x'))$.}

\RC{\emph{Conclusion.}
The map $\pi$ is locally Lipschitz and definable by Assumptions~\ref{ass:define} and~\ref{ass:lip}, hence path differentiable with a conservative Jacobian $J_\pi$~\citep{path-diff}.
By the chain rule~\citep{conservative}, $y \circ \pi$ is locally Lipschitz and path differentiable in $(\theta,x)$, with conservative Jacobian
\[
J_{y \circ \pi}(\theta,x) = \left\{ -U^{-1}V A \mid [U~V] \in J_R(y(\pi(\theta,x));\pi(\theta,x)),~A \in J_\pi(\theta,x) \right\}.
\]
Canonicalization~\eqref{eq:conic} keeps the decision variables of~\eqref{eq:robust_prob} as a subvector of the conic variable $z$.
Keeping the corresponding rows of these products, and the $\theta$ columns of $A$, gives $J_z(\theta)$, which proves Theorem~\ref{thm:subdiff}.}

\RC{In the implementation, the derivative with respect to $\pi$ is obtained by differentiating the solution of the homogeneous self-dual embedding, with the adjoint linear system solved iteratively by least squares (LSQR), as in~\citet[Appendix B]{cvxpylayers2019}.
In the case of singular systems, least squares is an efficient numerical heuristic.}

\section{Conservative Jacobians for the augmented Lagrangian}
\label{appx:jacobians}
\RC{We give the conservative Jacobians of Proposition~\ref{prop:jacobians}, and that of the constraint violation function $V$ used in Appendix~\ref{appx:converge_stationary}.
For each sample $i$, we write $z^i = z(\theta,x^i)$ for the solution of~\eqref{eq:robust_prob} and $J_{z^i}(\theta)$ for its conservative Jacobian from Theorem~\ref{thm:subdiff}.
For a batch $B$ and fixed $\lambda$ and $\mu$, we form the derivatives of $\hat{L}_B$ in $\alpha$ and $\theta$,
\allowdisplaybreaks \begin{align*}
v_\alpha &= \frac{(\lambda + \mu \hat{H}_{B}(\alpha,\theta))M_{B}}{|B|}\sum_{i \in B} \left(1 - M_i/\eta\right),\\
v_\theta &= \frac{1}{|B|}\sum_{i \in B} J_{\Phi_{\gamma}}(z^i)J_{z^i}(\theta) + \frac{(\lambda + \mu \hat{H}_{B}(\alpha,\theta))M_{B}}{\eta |B|}\sum_{i \in B} M_i J_{g}(z^i)J_{z^i}(\theta),
\end{align*}
where each $J$ is one selected element of the corresponding conservative Jacobian, $M_i \in \partial_+(g(z^i,u^i,x^i) - \alpha)$, and $M_B \in \partial_+\big((1/|B|)\sum_{i \in B} r_i(\alpha,\theta) - \kappa\big)$.
Here $\partial_+$ is the Clarke subdifferential of the positive part, with $\partial_+(a) = \{1\}$ for $a > 0$, $\partial_+(0) = [0,1]$, and $\partial_+(a) = \{0\}$ for $a < 0$.
We let $J_{\hat{L}_B}(\alpha,\theta)$ be the convex hull of the pairs $(v_\alpha,v_\theta)$ over all such selections, with the same $M_i$ and $M_B$ in both components.
The convex hull of a conservative Jacobian is again a conservative Jacobian~\citep[Section 2]{path-diff}, so $J_{\hat{L}_B}$ is a convex-valued conservative Jacobian of $\hat{L}_B$.}

\RC{We write $\partial^c f(z) \define \Co\{\lim_{z_i \rightarrow z} \nabla f(z_i) \mid z_i \in \dom(\nabla f)\}$ for the Clarke subdifferential~\citep{Clarke}.
The average $J_L = \Expect_B[J_{\hat{L}_B}]$ is a convex-valued conservative Jacobian of the augmented Lagrangian~\eqref{eq:lagrangian}, so it contains $\partial^c L$~\citep[Theorem 1]{path-diff}.
Similarly, let $J_{\hat{H}_B}$ be the convex hull, over the same selections, of the pairs of factors multiplying $\lambda + \mu\hat{H}_B(\alpha,\theta)$ in $(v_\alpha,v_\theta)$, and $J_F$ the average over $B$ of the convex hull of the first term of $v_\theta$.
Then $J_F$, $J_{\bar{H}} = \Expect_B[J_{\hat{H}_B}]$, and $J_V = \Expect_B[\hat{H}_B J_{\hat{H}_B}]$ are convex-valued conservative Jacobians of $F$, $\bar{H}$, and $V(\alpha,\theta) = (1/2)\Expect_B[\hat{H}_B(\alpha,\theta)^2]$.}

\section{Convergence analysis of the augmented Lagrangian algorithm}
\label{appx:converge_stationary}
We prove Theorem~\ref{thm:convergence}.
\RC{We first relate the inner-loop stationarity assumption of the theorem to Algorithm~\ref{alg:learning_inner}.
Step 1 then shows that the multipliers $\lambda^k$ converge, and Step 2 shows that, almost surely, every accumulation point of $\{(\alpha^k,\theta^k)\}$ is feasible for~\eqref{eq:bilevel_sampled} or a stationary point of $V$, that is, $0 \in J_V$.}

\paragraph{\RC{Inner loop.}} \RC{We fix $\lambda$ and $\mu$ and write $L(\alpha,\theta) = L(\alpha,\theta,\lambda,\mu)$.
Each solution map $\theta \mapsto z(\theta,x^i)$ is definable, since by Theorem~\ref{thm:subdiff} its graph is the projection onto $(\theta,z)$ of the definable set $\{(\theta,z,r) \mid R(z,r;\pi(\theta,x^i)) = 0\}$, so each $\hat{L}_B$ is definable in $(\alpha,\theta)$.
Theorem 3 of~\citet{path-diff} gives the convergence of SGD on a finite average of definable path-differentiable functions with a uniformly drawn index.
In Algorithm~\ref{alg:learning_inner}, $L = \Expect_B[\hat{L}_B]$ is indexed by the batch $B^t$, the step $(v_\alpha^t,v_\theta^t)$ has conditional mean in $J_L$, the iterates are bounded by Assumption~\ref{ass:bounded}, and the step sizes satisfy Assumption~\ref{ass:step_size}.
Therefore, for almost all initializations and scalings of the step sizes, $L(\alpha^t,\theta^t)$ converges almost surely, and every accumulation point of $\{(\alpha^t,\theta^t)\}$ satisfies $0 \in \partial^c L \subseteq J_L$, so it is a stationary point of $L$.}

\paragraph{\RC{Step 1, convergence of $\lambda$.}} \RC{The multiplier changes only in line 5 of Algorithm~\ref{alg:learning}, where $\lambda^k = \Pi_{[\lambda^{\min},\lambda^{\max}]}(\lambda^{k-1} + \mu^{k-1}\bar{H}(\alpha^k,\theta^k))$.
Since $\lambda^0 \in [\lambda^{\min},\lambda^{\max}]$, $\bar{H} \geq 0$, $\mu^{k-1} > 0$, and the projection onto an interval is monotone, $\{\lambda^k\}$ is nondecreasing and bounded above by $\lambda^{\max}$, so it converges.}

\paragraph{\RC{Step 2, convergence of the outer loop.}} \RC{By the assumption of Theorem~\ref{thm:convergence}, the iterates $(\alpha^k,\theta^k)$ lie in the compact set $\mathcal{D}_\alpha \times \mathcal{D}_\theta$, and we fix a subsequence $(\alpha^{k_j},\theta^{k_j})$ converging to an arbitrary accumulation point $(\bar\alpha,\bar\theta)$.

Suppose first that $\{\mu^k\}$ is bounded.
Since line 9 of Algorithm~\ref{alg:learning} multiplies $\mu$ by $\sigma > 1$, it runs only finitely often.
Thus, there is $k_0$ such that, for every $k \geq k_0$, the condition of line 4 holds and line 6 sets $H_{\rm best}$ to $\bar{H}(\alpha^k,\theta^k)$.
Therefore, $\bar{H}(\alpha^{k+1},\theta^{k+1}) \leq \tau \bar{H}(\alpha^k,\theta^k)$ for $k \geq k_0$, and $\bar{H}(\alpha^k,\theta^k) \to 0$ since $\tau < 1$.
The function $\bar{H}$ is continuous, since it is built from the locally Lipschitz maps $\theta \mapsto z(\theta,x^i)$ (Theorem~\ref{thm:subdiff}) and $g$ (Assumption~\ref{ass:lip}) by positive parts, averages, and affine maps in $\alpha$.
Hence $\bar{H}(\bar\alpha,\bar\theta) = 0$, and since $0 \leq H \leq \bar{H}$, the point $(\bar\alpha,\bar\theta)$ is feasible for~\eqref{eq:bilevel_sampled} and $\bar\theta$ satisfies~\eqref{eq:sampled_cvar}.

Suppose now that $\{\mu^k\}$ is unbounded.
Since $\mu^k$ never decreases, $\mu^k \to \infty$.
By the stationarity assumption, $0 \in J_L(\alpha^k,\theta^k)$, with $J_L$ taken at $(\lambda^{k-1},\mu^{k-1})$.
Each element of $J_{\hat{L}_B}$ is a convex combination of selections $(0,a_j) + (\lambda + \mu\hat{H}_B)h_j$, and since $\hat{H}_B$ is a common scalar it equals $(0,\bar{a}_B) + (\lambda + \mu\hat{H}_B)\bar{h}_B$ with $\bar{h}_B \in J_{\hat{H}_B}$.
Averaging over $B$ gives $0 = a^k + \lambda^{k-1} b^k + \mu^{k-1} c^k$, with $a^k = (0,a^k_\theta)$, $a^k_\theta \in J_F(\theta^k)$, $b^k \in J_{\bar{H}}(\alpha^k,\theta^k)$, and $c^k \in J_V(\alpha^k,\theta^k)$.
Conservative Jacobians are locally bounded and the iterates lie in a compact set, so $a^k$ and $b^k$ are bounded, and so is $\lambda^{k-1}$ by Step 1.
Dividing by $\mu^{k_j-1}$ gives $c^{k_j} = -(a^{k_j} + \lambda^{k_j-1} b^{k_j})/\mu^{k_j-1} \to 0$.
Since $J_V$ has a closed graph, $0 \in J_V(\bar\alpha,\bar\theta)$, so $(\bar\alpha,\bar\theta)$ is a stationary point of $V$.}

\section{\RC{Proof of the conditional guarantee}}
\label{appx:cond_proof}
\RC{We prove Theorem~\ref{thm:fin_sample_cond}.}
\begin{proof}
\RC{The proof mirrors that of Theorem~\ref{thm:fin_sample}.
We split the error of the $k$NN estimate into a variance term and a bias term, bound each uniformly over $\mathcal{W}$, and combine the two bounds with a union bound.

Let $k\mathrm{NN}(\bar{x})$ be the set of indices of the $k$ nearest neighbors of $\bar{x}$ among $x^1,\dots,x^N$.
With the $k$NN weights, $\sum_{i=1}^N w(\bar{x},x^i)\psi_w(u^i,x^i) = (1/k)\sum_{i \in k\mathrm{NN}(\bar{x})}\psi_w(u^i,x^i)$, so for every $w \in \mathcal{W}$
\[
\sum_{i=1}^N w(\bar{x}, x^i) \psi_w(u^i, x^i) - \Expect_{(u \mid x = \bar{x})}[\psi_w(u, \bar{x})] = S_{\mathrm{var}} + S_{\mathrm{bias}},
\]
with the variance and bias terms
\begin{align*}
S_{\mathrm{var}} &= \frac{1}{k}\sum_{i \in k\mathrm{NN}(\bar{x})} \left(\psi_w(u^{i}, x^{i}) - \Expect_{(u \mid x = x^{i})}[\psi_w(u, x^{i})]\right),\\
S_{\mathrm{bias}} &= \frac{1}{k}\sum_{i \in k\mathrm{NN}(\bar{x})} \Expect_{(u \mid x = x^{i})}[\psi_w(u, x^{i})] - \Expect_{(u \mid x = \bar{x})}[\psi_w(u, \bar{x})].
\end{align*}}

\paragraph{Variance bound.}
\RC{Conditionally on the contexts $X^N$, the index set $k\mathrm{NN}(\bar{x})$ is fixed, and the uncertain parameters $u^{i}$, $i \in k\mathrm{NN}(\bar{x})$, are independent with $u^i \sim \prob_{(u \mid x = x^i)}$.
The summands of $S_{\mathrm{var}}$ are therefore $k$ independent centered terms, which need not be identically distributed, and each $\psi_w(u^i,x^i)$ lies in $[0,1]$ by Assumption~\ref{assup:bound_range} and the normalization~\eqref{eq:g_setting}.
Since $\psi_w$ is $L_\psi$-Lipschitz in $w$, Lemma~\ref{lem:unif_dev} with $k$ in place of $N$ gives
\[
\prob^{N}\left (\sup_{w \in \mathcal{W}} |S_{\mathrm{var}}| \leq \tau_{\mathrm{var}} \;\middle|\; X^N\right ) \geq 1 - \beta(k,\tau_{\mathrm{var}}) = 1 - \beta_{\mathrm{var}},
\]
and the same bound holds unconditionally by taking the expectation over $X^N$.}

\paragraph{Bias bound.}
\RC{Let $d_k(\bar{x}) = \max_{i \in k\mathrm{NN}(\bar{x})} \|x^{i} - \bar{x}\|_2$ be the $k$NN distance of $\bar{x}$.
Since $\bar{x}$ and, almost surely, every $x^i$ lie in $\mathcal{D}_x$, Assumption~\ref{ass:holder} bounds the difference of the conditional expectations at $x^i$ and $\bar{x}$ by $L\|x^i - \bar{x}\|_2^s \leq L\, d_k(\bar{x})^s$, and averaging over $k\mathrm{NN}(\bar{x})$ gives $\sup_{w \in \mathcal{W}} |S_{\mathrm{bias}}| \leq L\, d_k(\bar{x})^s$.
Under Assumption~\ref{ass:positive}, $\prob_x$ satisfies at $\bar{x}$ the lower half of the full dimension condition of~\citet[Definition 3]{singh2016knn}, with dimension $p$, constant $\mu_{\min} V_p$, and locality radius $R$.
The upper tail bound of~\citet[Lemma 6]{singh2016knn} uses only this lower half, so for any $r$ in the interval of Theorem~\ref{thm:fin_sample_cond} it gives
\[
\prob^N(d_k(\bar{x}) > r) \leq \exp\left(-\mu_{\min}\, V_p\, r^p\, N\right)\left(\frac{e\,\mu_{\min}\, V_p\, r^p\, N}{k}\right)^k = \beta_{\mathrm{bias}}.
\]
On the event $d_k(\bar{x}) \leq r$, which has probability at least $1 - \beta_{\mathrm{bias}}$, we have $\sup_{w \in \mathcal{W}} |S_{\mathrm{bias}}| \leq L r^{s} = \tau_{\mathrm{bias}}$.}

\paragraph{Combining the bounds.}
\RC{Let $\mathcal{E}$ be the intersection of the variance and bias events, which has probability at least $1 - \beta_{\mathrm{var}} - \beta_{\mathrm{bias}}$ by a union bound.
Let $w^\star = (\theta^\star,\alpha^\star) \in \mathcal{W}$ be any point, possibly data-dependent, with $\tilde{H}(\alpha^\star,\theta^\star,\bar{x}) \leq \epsilon$.
By the definition of $\tilde{H}$ and the normalization~\eqref{eq:g_setting},
\[
\sum_{i=1}^N w(\bar{x}, x^i) \psi_{w^\star}(u^i, x^i) \leq \frac{c^{\min} + \kappa + \epsilon}{C},
\]
so on $\mathcal{E}$ the decomposition at $w^\star$ gives
\[
\Expect_{(u \mid x = \bar{x})}[\psi_{w^\star}(u, \bar{x})] \leq \frac{c^{\min} + \kappa + \epsilon}{C} + \tau_{\mathrm{var}} + \tau_{\mathrm{bias}}.
\]
Undoing the normalization as in the proof of Theorem~\ref{thm:fin_sample}, the choice of $\kappa$ gives
\[
\Expect_{(u \mid x = \bar{x})}[r_{w^\star}(u,\bar{x})] \leq \kappa + \epsilon + C(\tau_{\mathrm{var}} + \tau_{\mathrm{bias}}) \leq 0.
\]
Since the $\CVaR$ is the infimum of this expectation over $\alpha$, the $\CVaR$ of $g(z(\theta^\star,\bar{x}),u,\bar{x})$ at level $\eta$ under $\prob_{(u \mid x = \bar{x})}$ is nonpositive on $\mathcal{E}$, so the conditional chance constraint holds.
The event $\mathcal{E}$ does not depend on $w^\star$, so it covers every such point simultaneously, which gives~\eqref{eq:prob_guarantee_cond}.}
\end{proof}

\ifpreprint
\providecommand{\ecref}[1]{\ifpreprint Appendix~\ref{#1}\else \ref{#1}\fi}\section{Direct single-level reformulation}
\label{appx:bilinear}
\RC{We can reformulate the bi-level problem~\eqref{eq:bilevel_sampled} as a single-level problem by replacing each inner problem with its KKT conditions.
We do so for the $1$-norm uncertainty set with a polyhedral support $S = \{u \in \reals^m \mid Cu \leq d\}$,
\begin{equation*}
    \uncset(\theta,x) = \left\{u = A(x)v + b(x) \mid \|v\|_1 \leq \rho,~ u \in S\right\},
\end{equation*}
whose support function is $\sigma_S(\omega) = \min\{d^T\mu \mid C^T\mu = \omega,~\mu \geq 0\}$ by linear programming duality.
Substituting $\omega_l = C^T\mu_l$ with $\mu_l \geq 0$ into the norm-ball reformulation of \ecref{appx:gen_sets} with $\|\cdot\| = \|\cdot\|_1$, whose dual norm is $\|\cdot\|_\infty$, the inner problem~\eqref{eq:robust_prob} at a fixed $\theta$ and context $x^i$ becomes
\begin{equation*}
  \begin{array}[t]{ll}
		\underset{z_i,\,h_i,\,\mu_i,\,\lambda_i}{\mbox{minimize}} & z_{i,n+1} + \varepsilon\|z_{i,1:n}\|^2 \\
		\mbox{subject to} & [-g_l]^\star(z_i,h_{il}) + d^T\mu_{il} - b(x^i)^T(h_{il} + C^T\mu_{il}) + \rho\lambda_{il} \leq 0, \quad l = 1,\dots,L, \\
        & \|A(x^i)^T(h_{il} + C^T\mu_{il})\|_\infty \leq \lambda_{il}, \quad l = 1,\dots,L,\\
        & \mu_{il} \geq 0, \quad l = 1,\dots,L,
	\end{array}
\end{equation*}
with variables $z_i \in \reals^{n+1}$, whose first $n$ entries $z_{i,1:n}$ are the decision variables penalized as in Assumption~\ref{ass:unique}, and, for each $l$, $h_{il} \in \reals^m$, $\lambda_{il} \in \reals$, and $\mu_{il}$ of the dimension of $d$.
Any deterministic constraints on $z$ are included among the $g_l$ as functions that do not depend on $u$.
We write the constraints as $\tilde{g}_k \leq 0$, where for $l = 1,\dots,L$ and $j = 1,\dots,\hat{m}$,
\begin{align*}
    \tilde{g}_{1,il} &= [-g_l]^\star(z_i,h_{il}) + d^T\mu_{il} - b(x^i)^T(h_{il} + C^T\mu_{il}) + \rho\lambda_{il}, \\
    \tilde{g}_{2,ilj} &= \big(A(x^i)^T(h_{il} + C^T\mu_{il})\big)_j - \lambda_{il},\\
    \tilde{g}_{3,ilj} &= -\big(A(x^i)^T(h_{il} + C^T\mu_{il})\big)_j - \lambda_{il},\\
    \tilde{g}_{4,il} &= -\mu_{il}.
\end{align*}
This problem is convex in its variables. The conjugates $[-g_l]^\star$ are extended-valued, and for an affine piece, $[-g_l]^\star(z_i,h_{il})$ is finite only if $h_{il} = -P_l(z_i)$, which we impose as an equality constraint. We suppose that the remaining inequality constraints hold strictly at some point satisfying these equalities for every $\theta \in \mathcal{D}_\theta$ and every sample, so that its KKT conditions are necessary and sufficient for optimality.
With multipliers $\delta_{1,il}, \delta_{2,ilj}, \delta_{3,ilj} \in \reals$ and $\delta_{4,il}$ of the dimension of $\mu_{il}$, these conditions are, in addition to the primal feasibility constraints $\tilde{g}_k \leq 0$,
\begin{equation*}
	\label{eq:KKT}
	\tag{KKT}
\begin{aligned}
& 0 \in \partial_{(z_i,h_i,\mu_i,\lambda_i)} \Big( z_{i,n+1} + \varepsilon\|z_{i,1:n}\|^2 + \sum_{l=1}^{L}\delta_{1,il}\tilde{g}_{1,il}\\
&\qquad\qquad + \sum_{l=1}^{L}\sum_{j=1}^{\hat{m}}\big(\delta_{2,ilj}\tilde{g}_{2,ilj}+ \delta_{3,ilj}\tilde{g}_{3,ilj}\big) + \sum_{l=1}^{L}\delta_{4,il}^T\tilde{g}_{4,il}\Big),\\
&\delta_{1,il} \geq 0, \quad \delta_{2,ilj} \geq 0, \quad \delta_{3,ilj} \geq 0, \quad \delta_{4,il} \geq 0,\\
& \delta_{1,il}\tilde{g}_{1,il} = 0, \quad \delta_{2,ilj}\tilde{g}_{2,ilj} = 0, \quad \delta_{3,ilj}\tilde{g}_{3,ilj} = 0, \quad \delta_{4,il} \circ \tilde{g}_{4,il} = 0,
\end{aligned}
\end{equation*}
for $l = 1,\dots,L$ and $j = 1,\dots,\hat{m}$, where $\partial$ is the convex subdifferential, which replaces the gradient where the conjugates $[-g_l]^\star$ are not differentiable, and $\circ$ is the componentwise product.
Replacing each inner problem by these conditions gives the single-level problem
\begin{equation}
	\label{eq:singlelevel-bilinear}
	\begin{array}{ll}
 \underset{\theta,\,\alpha,\,z,\,h,\,\mu,\,\lambda,\,\delta}{\mbox{minimize}} & \frac{1}{N}\sum_{i=1}^N \big(\gamma f(z_i,u^i,x^i) + z_{i,n+1}\big)\\
		\mbox{subject to}
        &  \frac{1}{N}\sum_{i=1}^N \Big(\frac{(g(z_i,u^i,x^i)- \alpha)_+}{\eta} + \alpha\Big) \leq \kappa, \\
        & \tilde{g}_{1,il} \leq 0,\quad \tilde{g}_{2,ilj} \leq 0,\quad \tilde{g}_{3,ilj} \leq 0,\quad \tilde{g}_{4,il} \leq 0,\\
        & \qquad i = 1,\dots,N,~l = 1,\dots,L,~j = 1,\dots,\hat{m},\\
        & \eqref{eq:KKT}, \quad i = 1,\dots,N,
	\end{array}
\end{equation}
where $z$, $h$, $\mu$, $\lambda$, and $\delta$ collect the inner variables and multipliers of all $N$ samples, and the first constraint is $H(\alpha,\theta) = 0$.
Under Assumption~\ref{ass:unique} the inner solution $z(\theta,x^i)$ is unique, so a point is feasible for~\eqref{eq:singlelevel-bilinear} exactly when $z_i = z(\theta,x^i)$ for every $i$ and $(\alpha,\theta)$ is feasible for~\eqref{eq:bilevel_sampled}, and the two problems have the same optimal set parameters $\theta$.
The set parameter $\theta$ enters through $A(x^i)$ and $b(x^i)$.
For the linear class~\eqref{eq:lin_map}, the terms $b(x^i)^T(h_{il} + C^T\mu_{il})$ and $A(x^i)^T(h_{il} + C^T\mu_{il})$ are therefore bilinear in $\theta$ and the auxiliary variables, and the complementarity conditions multiply them further by the multipliers $\delta$.
For the MLP class, the dependence of these terms on $\theta$ is piecewise polynomial in the network weights, and the reformulation is no easier.
The single-level problem is thus a nonconvex program with complementarity constraints, whose number grows with $N$, $L$, and $\hat{m}$.
Global solvers handle such problems only at small sizes, which is why we search for stationary points of the bi-level formulation instead.}

 \fi

\ifpreprint
\section{General reformulation proofs}
\label{appx:reform}
\providecommand{\ecref}[1]{\ifpreprint Appendix~\ref{#1}\else \ref{#1}\fi}\ifpreprint
\subsection{Conjugate of a general function}
\else 
\section{Conjugate of a general function}
\fi 
We first derive the conjugate of general functions, which can then be used in the upcoming reformulations.
\paragraph{\RC{Affine function.}}
\RC{For an affine piece $g_l(z,u) = P_l(z)^Tu + a_l(z)$, as in the maximum of affine form, we compute the conjugate as}
\begin{equation*}
	[-g\RC{_l}]^\star(z,h) =  \RC{\sup_u\, h^Tu + P_l(z)^Tu + a_l(z)} =
	\begin{dcases}
            \RC{a_l(z)}  & \text{ if}\; h = \RC{-P_l(z)} \\
		\infty & \text{ otherwise}.
	\end{dcases}
\end{equation*}
\RC{Substituting $-P_l(z)$ for the argument of $[-g_l]^\star$ in the reformulations below gives the maximum of affine reformulations of \ecref{appx:gen_sets}.}
\paragraph{Sum of concave functions.}
For a general function $g(z,u) = \sum_{j=1}^J g_j(\RC{z,u})$, where each $g_j$ is concave in $u$, by the theorem of infimal convolutions~\citep[Theorem 11.23(a), p. 493]{rockafellar_wets_1998}, we compute the conjugate as
\begin{equation*}
	[-g]^\star(z,h) = \begin{array}[t]{ll}
		\underset{\RC{h_1},\dots,h_J}{\inf} & \sum_{j=1}^J [-g_j]^\star(z,h_j)\\
		\mbox{subject to} & \sum_{j=1}^J h_j = h.
	\end{array}
\end{equation*}

\paragraph{Maximum of concave functions.}
We also consider the maximum of concave functions, $$g(z,u) = \max_{l=1,\dots,L} g_l(z,u).$$ In practice, the functions $g_l(z,u)$ are often affine. For all the reformulations below, we assume $g(z,u)$ to be of this maximum of concave form.

\ifpreprint
\subsection{Box and ellipsoidal reformulations}
\else
\section{Box and ellipsoidal reformulations}
\fi
\label{apx:boxproof}
Rewriting the robust constraint as an optimization problem, we have
{\allowdisplaybreaks \begin{align*}
& \begin{array}[t]{ll}
	\underset{Av+b\in S}{\sup} & g(z,Av+b) \\
	\mbox{subject to} & \RC{\|v\|} \le 1
\end{array}\\
&= \begin{array}[t]{ll}
	\underset{Av+b\in S}{\sup}  ~\underset{\lambda \geq 0}{\inf}~g(z,Av+b)  - \lambda\RC{\|v\|}  + \lambda\\
	\end{array}\\
	&= \begin{array}[t]{ll}
	\underset{Av+b\in S}{\sup} ~ \underset{\lambda \geq 0}{\inf}  ~\underset{l}{\max}~g_l(z,Av+b)  - \lambda\RC{\|v\|}  + \lambda\\
		\end{array}\\
		&= \begin{array}[t]{ll}
			\underset{l}{\max}~ \underset{Av+b\in S}{\sup} ~ \underset{\lambda \geq 0}{\inf}  ~g_l(z,Av+b)  - \lambda\RC{\|v\|}  + \lambda\\
				\end{array}\\
			&= \begin{array}[t]{ll}
					\underset{l}{\max}~ \underset{\lambda \geq 0}{\inf} ~ \underset{Av+b\in S}{\sup}~ g_l(z,Av+b)  - \lambda\RC{\|v\|}  + \lambda\\
						\end{array}\\
&= \begin{array}[t]{ll}
	\underset{\lambda_l \geq 0, \tau}{\inf} & \tau\\
	\mbox{subject to} & \lambda_l + \underset{Av+b\in S}{\sup } g_l(z,Av+b) - \lambda_l  \RC{\|v\|} \leq \tau, \quad l = 1,\dots,L,
\end{array}
\end{align*}}
\RC{where $\|\cdot\|$ is the norm defining the set, which is the $\infty$-norm for box uncertainty.} We used the Lagrangian in the first equality, and decomposed $g(z,Av+b)$ into its constituents for the second equality. As $l$ and $\lambda$ are separable, in the third step we interchanged the $\inf$ and the $\max$. Next, for each $l$, \RC{we interchanged the $\inf$ and the $\sup$ by Lagrangian duality, which holds, with the infimum attained, under the Slater condition that some $v$ with $\|v\| < 1$ satisfies $Av+b \in {\rm ri}\, S \cap {\rm ri}\, {\rm dom}\, g_l(z,\cdot)$}. We then express the $\max$ in epigraph form by introducing an auxiliary variable $\tau$.

Now, if we define the function $c(v) = \RC{\|v\|}$, then using the definition of conjugate functions, we have, for each $l$, the left-hand side of the constraint equal to
\begin{equation*}
	\lambda_l + [-g_l(z,A\cdot+b) +\chi_S(A\cdot+b) + \lambda_l c]^\star(0).
\end{equation*}
We then borrow results from~\citet[Theorem 11.23(a), p. 493]{rockafellar_wets_1998},~\citet[Lemma C.9]{zhen2021mathematical}, and~\citet[Proposition 13.24]{composites} with regard to the inf-convolution of conjugate functions, the conjugate functions of norms, and the conjugate of linear compositions, to note
\begin{align*}
	[-g_l(z,A\cdot+b)   +\chi_S(A\cdot+b)+ \lambda_l c]^\star (0) =  \underset{\omega_l, h_l, y_l}{\inf} ([-g_l]^\star(z,h_l) &- b^T(h_l + \omega_l) + \sigma_S(\omega_l)\\
  &+ \indicator_{\RC{\|y_l\|_*} \leq \lambda_l} + \indicator_{A^T(h_l + \omega_l)= y_l}),
\end{align*}
where $h_l\in\reals^{m}$, $y_l \in \reals^{\hat{m}}$, $\omega_l \in \reals^{m}$, and \RC{$\|\cdot\|_*$ is the dual norm}. Substituting these back\RC{, and minimizing over $\lambda_l \geq 0$, which gives $\lambda_l = \|y_l\|_*$}, we arrive at
\begin{equation*}
\begin{aligned}
\begin{array}[t]{ll}
\underset{z,\tau, h, \omega, y}{\inf} & \tau\\
\mbox{subject to} & [-g_l]^\star(z,h_l)+ \sigma_S(\omega_l)  -b^T(h_l + \omega_l) + \RC{\|y_l\|_*} \leq \tau, \quad l = 1,\dots,L \\
 & A^T(h_l + \omega_l)= y_l,  \quad l = 1,\dots,L.
\end{array}
\end{aligned}
\end{equation*}
Recall that this optimization problem upper bounds our original robust constraint, which we constrain to be nonpositive. Therefore, we only need $\tau \leq 0$. For the final step, we can now write the original minimization problem~\eqref{eq:robust_prob} in minimization form for both the objective and constraints, to arrive at
\begin{equation*}
	\begin{aligned}
	\begin{array}[t]{ll}
	\underset{z, h, \omega}{\text{minimize}} & f(z)\\
	\mbox{subject to} & [-g_l]^\star(z,h_l)+ \sigma_S(\omega_l)  -b^T(h_l + \omega_l) + \RC{\|A^T(h_l + \omega_l) \|_*} \leq 0, \quad l = 1,\dots,L.	\end{array}
	\end{aligned}
	\end{equation*}

\ifpreprint
\subsection{Budget reformulation}
\else 
\section{Budget reformulation}
\fi
\label{apx:budproof}
Rewriting the robust constraint as an optimization problem, we have
\begin{equation*}
	\begin{aligned}
& \begin{array}[t]{ll}
	\underset{\RC{Av+b}\in S}{\sup} & g(z,Av+b) \\
	\mbox{subject to} & \| v \|_\infty \le \rho_1\\
	 & \| v \|_1 \le \rho_2\\
\end{array}\\
&= \begin{array}[t]{ll}
	\underset{Av+b\in S}{\sup}~ \underset{\RC{\lambda_1,\lambda_2} \geq 0}{\inf} ~\underset{l}{\max} ~ g_l(z,Av+b)  - \lambda_1\| v \|_\infty  + \lambda_1\rho_1 - \lambda_2\| v\|_1  + \lambda_2\rho_2\\
	\end{array}\\
 &=  \begin{array}[t]{ll}
	\underset{l}{\max} ~\underset{Av+b\in S}{\sup} ~ \underset{\RC{\lambda_1,\lambda_2} \geq 0}{\inf} ~g_l(z,Av+b)  - \lambda_1\| v \|_\infty  + \lambda_1\rho_1 - \lambda_2\| v\|_1  + \lambda_2\rho_2\\
	\end{array}\\
	&= \begin{array}[t]{ll}
		\underset{l}{\max}~\underset{\RC{\lambda_1,\lambda_2} \geq 0}{\inf}~\underset{Av+b\in S}{\sup} ~g_l(z,\RC{Av+b})  - \lambda_1\| v\|_\infty  + \lambda_1\rho_1 - \lambda_2\| v \|_1  + \lambda_2\rho_2,\\
		\end{array}\\
\end{aligned}
\end{equation*}
where we again used the Lagrangian, separated $g$ into its constituents, and \RC{interchanged the $\inf$ and the $\sup$ by Lagrangian duality under the Slater condition}. We now express the max in epigraph form,
\begin{equation*}
	\begin{aligned}
	\begin{array}[t]{ll}
		\underset{\lambda_{l} \geq 0, \tau}{\inf} & \tau \\
		\mbox{subject to} & \lambda_{l1}\rho_1 + \lambda_{l2}\rho_2 + \underset{Av+b\in S}{\sup } g_l(z,Av+b) -\lambda_{l1}\| v \|_\infty - \lambda_{l2}\| v \|_1 \leq \tau, \quad l = 1,\dots,L.
	\end{array}
	\end{aligned}
	\end{equation*}

If we define the functions $c_1(v) = \|v\|_\infty$, $c_2(v) = \|v\|_1$, then using the definition of conjugate functions, we have for each $l$ the left-hand side of the constraint equal to
\begin{equation*}
\begin{aligned}
		& \lambda_{l1}\rho_1 +  \lambda_{l2}\rho_2 +  [-g_l(z,A\cdot+b) +\chi_S(A\cdot+b) + \lambda_{l1} c_1 + \lambda_{l2} c_2]^\star(0)\\
		&= 	\underset{h_l,\omega_l, y_l}{\inf} \lambda_{l1}\rho_1 + \lambda_{l2}\rho_2 + [-g_l]^\star(z,h_l) - b^T(h_l + \omega_l) + \sigma_S(\omega_l)  +  [\lambda_{l1} c_1 + \lambda_{l2} c_2 ]^\star(-y_l) \\
  &\hspace{12cm}+ \indicator_{A^T(h_l + \omega_l) = y_l} \\
		&= \underset{h_l,\omega_l, y_l, r_l}{\inf} \lambda_{l1}\rho_1 + \lambda_{l2}\rho_2 + [-g_l]^\star(z,h_l) - b^T(h_l + \omega_l)+ \sigma_S(\omega_l)  + \indicator_{\|r_l-y_l\|_1 \leq \lambda_{l1}} + \indicator_{\|r_l\|_\infty \leq \lambda_{l2}} \\
  & \hspace{12cm} + \indicator_{A^T(h_l + \omega_l) = y_l}.
\end{aligned}
\end{equation*}
\RC{Minimizing over $\lambda_{l1}, \lambda_{l2} \geq 0$ gives $\lambda_{l1} = \|r_l - y_l\|_1$ and $\lambda_{l2} = \|r_l\|_\infty$, so} we can express this as
\begin{equation*}
	\begin{aligned}
	\begin{array}[t]{ll}
		\underset{z, h_l,\omega_l, y_l, r_l, \tau}{\inf} & \tau \\
		\mbox{subject to} & [-g_l]^\star(z,h_l)  + \sigma_S(\omega_l) -  b^T(h_l + \omega_l) + \rho_1\|r_l - y_l\|_1  + \rho_2\|r_l\|_\infty \leq \tau, \\
  &\hfill l =  1, \dots,L \\
	 & A^T(h_l + \omega_l) = y_l,  \quad l =  1, \dots,L.
	\end{array}
	\end{aligned}
	\end{equation*}
 We can then let $\tau = 0$ and substitute the new constraints in for the robust constraint in~\eqref{eq:robust_prob},
\begin{equation*}
	\begin{array}[t]{ll}
		{\text{minimize}} & f(z) \\
		\mbox{subject to} & [-g_l]^\star(z,h_l)  + \sigma_S(\omega_l) -  b^T(h_l + \omega_l) + \rho_1\|r_l - A^T(h_l + \omega_l)\|_1  + \rho_2\|r_l\|_\infty \leq 0,\\
  &\hfill l =  1, \dots,L.	\end{array}
	\end{equation*}

\ifpreprint
\subsection{Polyhedral reformulation}
\else 
\section{Polyhedral reformulation}
\fi 
\label{apx:polyproof}
Rewriting the robust constraint as an optimization problem, we have
{\allowdisplaybreaks \begin{align*}
& \begin{array}[t]{ll}
	\underset{Av+b\in S}{\sup} & g(z,Av+b) \\
	\mbox{subject to} & Dv \le d
\end{array}\\
&= \begin{array}[t]{ll}
	\underset{Av+b\in S}{\sup} ~\underset{h \geq 0}{\inf} ~g(z,Av+b)  - (Dv)^Th   + d^Th\\
	\end{array}\\
	&= \begin{array}[t]{ll}
		\underset{Av+b\in S}{\sup}~\underset{h \geq 0}{\inf}~\underset{l}{\max}~g_l(z,Av+b)  - (Dv)^Th   + d^Th\\
		\end{array}\\
		&= \begin{array}[t]{ll}
			\underset{l}{\max}~\underset{Av+b\in S}{\sup}~\underset{h \geq 0}{\inf}~g_l(z,Av+b)  - (Dv)^Th   + d^Th\\
			\end{array}\\
		&= \begin{array}[t]{ll}
				\underset{l}{\max}~\underset{h \geq 0}{\inf}~\underset{Av+b\in S}{\sup}~g_l(z,Av+b)  - (Dv)^Th   + d^Th\\
				\end{array}\\
&= \begin{array}[t]{ll}
	\underset{h\geq 0,\tau}{\inf} & \tau \\
	\mbox{subject to}& d^Th_l + \underset{Av+b\in S}{\sup } g_l(z,Av+b) - (Dv)^Th_l \leq \tau, \quad l = 1,\dots,L
\end{array}\\
&= \begin{array}[t]{ll}
	\underset{h \geq 0, \omega,\tau, y, \gamma}{\inf} &\tau\\
	\mbox{subject to}& d^Th_l + [-g_l]^\star(z,\omega_l) - b^T(\omega_l + \gamma_l) + \sigma_S(\gamma_l) + \indicator_{D^Th_l = -y_l} \\
 &\hfill  + \indicator_{A^T(\omega_l+ \gamma_l) = y_l} \leq \tau, \quad l = 1,\dots,L
\end{array}\\
&= \begin{array}[t]{ll}
	\underset{h \geq 0,\omega,\tau,\gamma}{\inf} &\tau\\
	\mbox{subject to}& [-g_l]^\star(z,\omega_l) \RC{{}+ d^Th_l} - b^T(\omega_l + \gamma_l) + \sigma_S(\gamma_l) \leq \tau, \quad l = 1,\dots,L\\
	& A^T(\omega_l+ \gamma_l) = -D^Th_l,\quad l = 1,\dots,L,
\end{array}
\end{align*}}
where we again used the Lagrangian, \RC{Lagrangian duality under the Slater condition}, the infimal convolution of conjugates, and the conjugate of affine compositions. The final reformulation of~\eqref{eq:robust_prob} is then
\begin{equation*}
\begin{array}[t]{ll}
	{\text{minimize}} & f(z)\\
	\mbox{subject to} & [-g_l]^\star(z,\omega_l) \RC{{}+ d^Th_l} - b^T(\omega_l + \gamma_l) + \sigma_S(\gamma_l) \leq 0, \quad l = 1,\dots,L\\
	 &A^T(\omega_l+ \gamma_l) = -D^Th_l,\quad l = 1,\dots,L\\
	 & \RC{h_l \geq 0, \quad l = 1,\dots,L}.
\end{array}
\end{equation*}

\ifpreprint
\subsection{Mean robust optimization reformulation}
\else
\section{Mean robust optimization reformulation}
\fi
\label{apx:mroproof}
We adapt the proof from~\citet[Appendix A.3]{mro}.
To simplify notation, we define $c_k(v_{lk}) \define   \|v_{lk} - \bar{d}_k \|^{\RC{r}} - \rho^{\RC{r}}$.
We also note that $u_{lk} = Av_{lk}+b$. 
Then, we have
\begin{equation*}
	\begin{aligned}
& \begin{array}[t]{ll}
	\underset{u_{11}, \dots, u_{LK} \in \supp, \alpha \in \Gamma}{\sup} &\sum_{k=1}^K \sum_{l=1}^L \alpha_{lk} g_l(z,u_{lk}) \\
	\mbox{subject to} & \sum_{k=1} ^K \sum_{l=1}^L \alpha_{lk} c_k(v_{lk}) \le 0
\end{array}\\
&=\begin{array}[t]{ll}
	\underset{u_{11}, \dots, u_{LK} \in \supp, \alpha \in \Gamma}{\sup} & \underset{\lambda \geq 0}{\inf}~\sum_{k=1}^K \sum_{l=1}^L \alpha_{lk}g_l(z,u_{lk}) - \lambda\sum_{k=1} ^K \sum_{l=1}^L \alpha_{lk} c_k(v_{lk})\\
\end{array}\\
&=\begin{array}[t]{ll}
	\underset{\alpha \in \Gamma}{\sup}~\underset{\lambda \geq 0}{\inf}~\underset{u_{11}, \dots, u_{LK} \in \supp}{\sup} &\sum_{k=1}^K \sum_{l=1}^L \alpha_{lk} g_l(z,u_{lk}) - \lambda\sum_{k=1} ^K \sum_{l=1}^L \alpha_{lk} c_k(v_{lk}).\\
\end{array}\\
\end{aligned}
\end{equation*}
We applied the Lagrangian in the first equality. \RC{Then, for each fixed $\alpha$, we interchanged the inf and the sup over $u_{11},\dots,u_{LK}$ by Lagrangian duality under the Slater condition, since each $g_l(z,u_{lk})$ is concave in $u_{lk}$.
Next, we rewrite the formulation using an epigraph trick.
Since each $\alpha_{lk}$ is nonnegative and the objective is separable in the $v_{lk}$, the weights can be pulled outside the inner supremum, which gives
\begin{equation*}
	\begin{aligned}
&= \begin{array}[t]{ll}
	\underset{\alpha \in \Gamma}{\sup}~\underset{\lambda \geq 0}{\inf} &\sum_{k=1}^K s_k  \\
	\mbox{subject to}& \sum_{l=1}^L \alpha_{lk} \underset{Av_{lk}+b \in \supp}{\sup} \left( g_l(z,Av_{lk}+b) - \lambda c_k(v_{lk}) \right) \le s_k, \quad k = 1,\dots,K\\
\end{array}\\
&= \begin{array}[t]{ll}
	\underset{\alpha \in \Gamma}{\sup}~\underset{\lambda \geq 0}{\inf} &\sum_{k=1}^K s_k  \\
	\mbox{subject to}& \sum_{l=1}^L \alpha_{lk}[-g_l(z, A\cdot + b) +\chi_{\supp}(A\cdot + b) + \lambda c_k]^\star(0) \le s_k, \quad k = 1,\dots,K,
\end{array}\\
\end{aligned}
\end{equation*}
where a term with $\alpha_{lk} = 0$ is zero, since its $v_{lk}$ does not enter the objective, and the last equality is the definition of the conjugate function.
By the conjugates of sums and affine compositions,
\begin{equation*}
	\begin{aligned}
		[-g_l(z, A\cdot + b) +\chi_{\supp}(A\cdot + b) + \lambda c_k]^\star (0) &= \underset{\omega_{lk}, h_{lk}, y_{lk}}{\inf} ([-g_l]^\star(z,h_{lk}) - b^T(h_{lk} + \omega_{lk})\\
  &+ \indicator_{A^T(h_{lk} + \omega_{lk})= y_{lk}} + \sigma_S(\omega_{lk})+[\lambda c_k]^\star (-y_{lk})),
	\end{aligned}
\end{equation*}}
and
\begin{equation*}
\begin{aligned}
	[\lambda c_k]^\star (-y_{lk}) =-y_{lk}^T\bar{d}_k + \RC{\phi(r_*)}\lambda {\|}y_{lk}/\lambda{\|}^{\RC{r_*}}_* + \lambda \rho^{\RC{r}}\RC{,}
 \end{aligned}
\end{equation*}
\RC{where $\phi(r_*)$ and the closed perspective $\lambda\|y/\lambda\|_*^{r_*}$ are as in \ecref{appx:gen_sets}.}
Substituting these back, we arrive at
\begin{equation*}
	\begin{aligned}
\begin{array}[t]{ll}
	\underset{\alpha \in \Gamma}{\sup}~\underset{\lambda \geq 0, h_{lk},\omega_{lk}, s_k}{\inf} & \sum_{k=1}^K s_k\\
	\mbox{subject to} & \sum_{l=1}^L \alpha_{lk}([-g_l]^\star(z,h_{lk}) -b^T(h_{lk} + \omega_{lk})\\
	& + \sigma_{\supp}(\omega_{lk}) - (A^T(h_{lk} + \omega_{lk}))^T \bar{d}_k + \RC{\phi(r_*)}\lambda {\|}A^T(h_{lk} + \omega_{lk})/\lambda{\|}^{\RC{r_*}}_*\\
	& \hfill + \lambda \rho^{\RC{r}} ) \le s_k, \\
	&\hspace{2cm} \quad k = 1,\dots,K.
\end{array}
\end{aligned}
\end{equation*}
Taking the supremum over $\alpha$, noting that $\sum_{l=1}^L \alpha_{lk} = w_k$ for all $k$, we arrive at
\begin{equation*}
\begin{aligned}
\begin{array}[t]{ll}
\underset{\lambda \geq 0, h_{lk},\omega_{lk}, s_k}{\inf} & \lambda \rho^{\RC{r}} + \sum_{k=1}^K w_k s_k\\
\mbox{subject to} & [-g_l]^\star(z,h_{lk}) -b^T(h_{lk} + \omega_{lk}) + \sigma_S(\omega_{lk}) - (A^T(h_{lk} + \omega_{lk}))^T \bar{d}_k \\
& \hfill + \RC{\phi(r_*)}\lambda \left\|A^T(h_{lk} + \omega_{lk})/\lambda \right\|^{\RC{r_*}}_* \le s_k,\\
& \hfill l = 1,\dots,L, \quad k = 1,\dots,K,
\end{array}
\end{aligned}
\end{equation*}
for which the final reformulation of~\eqref{eq:robust_prob} is
\begin{equation*}
	\begin{aligned}
	\begin{array}[t]{ll}
	{\text{minimize}} & f(z)\\
	\mbox{subject to} &\lambda \rho^{\RC{r}} + \sum_{k=1}^K w_k s_k \leq 0\\
	 & [-g_l]^\star(z,h_{lk}) -b^T(h_{lk} + \omega_{lk}) + \sigma_S(\omega_{lk}) - (A^T(h_{lk} + \omega_{lk}))^T \bar{d}_k \\
& \hfill + \RC{\phi(r_*)}\lambda \left\|A^T(h_{lk} + \omega_{lk})/\lambda \right\|^{\RC{r_*}}_* \le s_k,\\
& \hfill l = 1,\dots,L, \quad k = 1,\dots,K\\
	 & \RC{\lambda \geq 0.}
	\end{array}
	\end{aligned}
	\end{equation*}
\RC{For $r = 1$, the term $\phi(r_*)\lambda\|A^T(h_{lk} + \omega_{lk})/\lambda\|_*^{r_*}$ is the indicator of $\|A^T(h_{lk} + \omega_{lk})\|_* \le \lambda$ and $\lambda\rho^r = \lambda\rho$, which gives the order-one form stated in \ecref{appx:gen_sets}.} \fi

\ifpreprint
\section{\RC{Experimental details and hyperparameter ablations}}
\label{appx:ablation}

\RC{We give the full data-generating processes summarized in Section~\ref{sec:examples}, as well as hyperparameter ablations for the \gls{LRO} sets.}

\subsection{\RC{Data generation}}
\label{appx:data}
\RC{We begin with the data generation procedures, which create nonlinear contextual mappings for both the mean and the covariance. }

\subsubsection{\RC{Portfolio optimization}}
\label{appx:data_port}

\RC{We consider $n = 10$ stocks in $k = 3$ sectors and $n = 30$ stocks in $k = 5$ sectors, with $m = n$.
Every sample draws its own context $x = (x_{\rm sec}, x_{\rm vol}, x_{\rm mom})$ uniformly on $[-1,1]^{k+2}$, consisting of one state per sector measuring how favorable current conditions are, a market volatility state, and a momentum state.
Conditional on the context, log returns are Gaussian, and the arithmetic return is their exponential minus one, which keeps the gross return positive.
We clip the arithmetic returns to $[-1,3]$, so that the gross return lies in $[0,4]$ and the support is compact. 
The heat of a sector maps its favorability state to $(0,1)$ through the logistic function
\begin{equation}
	\label{eq:port_heat}
	w_j(x) = \varphi\big(8\,x_{{\rm sec},j}
	+ x_{{\rm sec},j}(Ax_{\rm sec})_j\big),\quad j = 1,\dots,k,
	\qquad \varphi(t) = 1/(1+e^{-t}),
\end{equation}
where $A$ is the ring adjacency of the sectors, so favorable conditions spill over to neighboring sectors.
The vector $h(x) = Sw(x) \in \reals^n$ gives each stock its sector's heat, and the conditional log-return moments are}
\RC{\begin{align*}
	m(x) &= \mu^{\rm nom} + \nu\big(h(x) - \tfrac12 \ones\big)
	+ \nu\big(0.25\,x_{\rm mom}\ones + 0.6\,E^T(x_{\rm vol},x_{\rm mom})\big),\\
	\Sigma(x) &= 0.15^2 \Big( 0.25 \diag\big(0.4\,\ones + 1.2\,h(x)\big)
	+ 0.75 \big(1 + 0.8\,x_{\rm vol}\big)\,\ell\ell^T \Big),
\end{align*}
where $\mu^{\rm nom}$ has entries drawn uniformly from $[0,0.01]$, the matrix $E$ holds fixed mean-zero per-stock loadings on the volatility and momentum states, and $\ell$ holds fixed market-factor loadings with unit mean square.
In other words, a stock's expected return rises with its sector's heat and with the momentum, and returns are correlated through a common market shock whose weight grows with the volatility, on top of stock-specific noise that is larger in hot sectors.
The scale $\nu$ ($0.149$ at $n = 10$, $0.138$ at $n = 30$) is calibrated so that the spread of the conditional mean across contexts is twice the return standard deviation, a context signal that is strong but still plausible for a market.
Centering the heat at $1/2$ keeps the context average of each stock's conditional mean log return at its nominal value $\mu^{\rm nom}$.
With this data-generating procedure, the mean and covariance mappings are clearly contextual and mostly affine, with slight nonlinearity due to the logistic mapping~\eqref{eq:port_heat}.}

\RC{The loss budgets are geometric ladders, permuted so that a stock's budget tier is independent of its sector.
The per-stock budgets $\beta$ have geometric mean $\tau_a$ and spread $3$, and the sector budgets $\beta^{\rm sec}$ have geometric mean $0.10$ and spread $2$.
We set $\tau_a = 0.04$ at $n=10$ and $\tau_a = 0.02$ at $n=30$ on a separate calibration pass, so that about $3$ budget rows are active at the \gls{ETO} set's radius at $n=10$ and about $6$ at $n=30$.}

\subsubsection{\RC{Inventory management}}
\label{appx:data_inv}

\RC{We let $n = m = 10$.
The nominal revenues $\bar r$ and capacities $c$ are geometrically spaced over $[2,4]$ and $[50,30]$ respectively, the nominal demand levels $d$ are uniform on $[10,20]$, the transportation and holding costs $t,h$ are uniform on $[0.1,0.3]$, and the shared capacity $C$ is set as 92\% of the expected total demand.}

\RC{For this uncertain demand, each sample $j$ depends on its own context $x^j = (r^j, \xi^j)$, with $\xi^j$ uniform on $[0,1]^{k+3}$, and $r^j = \bar r + \varepsilon^j$ for $\varepsilon^j \sim N(0,(0.05\bar r)^2)$ truncated to $[-0.3\bar r, 0.3\bar r]$, so that $r^j \in [0.7\bar r, 1.3\bar r]$ has compact support.
We write $\xi = (\xi_{\rm reg}, \xi_{\rm pro}, \xi_{\rm sea}, \xi_{\rm com})$ for the regime, promotion, season, and competitor contexts, and $W \in \{0,1\}^{n\times k}$ for the fixed product-to-regime indicator, with $k = 3$.
The conditional mean and covariance of demand are}
\RC{\begin{align*}
	\mu(x) &= d - 0.1\,r + 15\, W \xi_{\rm reg} + \Big[1.3 \big(\xi_{\rm pro} - \tfrac12\big) + 0.15 \big(\xi_{\rm com} - \tfrac12\big) - 1.3\big((\xi_{\rm pro}-\tfrac12)^2 - \tfrac1{12}\big)\Big] d,\\
	\Sigma(x) &= 0.25\,\diag\!\Big((0.3\,\ones + 6\, W \xi_{\rm reg}) \odot \big(\mu(x)/\bar\mu\big)^{2}\Big) + \big(0.495 + 3.96\,\xi_{\rm sea}\big)\,\ell\ell^T,
\end{align*}
where $\bar\mu \approx 22$ is the expected mean demand.}
\RC{Under this regime, both the mean and covariance are mostly linear in the context, with a small polynomial dependence, to mimic a real-world contextual dependence that would be difficult to exactly recover.}
\RC{The vector $\ell \in \reals^n$ denotes fixed common-factor loadings, so the second covariance term is a rank-one seasonal factor that accounts for approximately 75\% of the total covariance.}
\RC{Demand is drawn as $u^j = \mu(x^j) + e^j + f^j \ell$ and clipped to $[0,120]$, where the product-level noise $e^j \sim N(0, D(x^j))$ has the diagonal first term of $\Sigma(x^j)$ as its covariance $D(x^j)$, and the common factor $f^j$ is Student-$t$ with $5$ degrees of freedom, rescaled to variance $0.495 + 3.96\,\xi^j_{\rm sea}$. In this manner, the conditional covariance is $\Sigma(x^j)$ up to the clipping, which in our instances affects $0.03\%$ of the drawn components. An \gls{ETO} method that relies on the first two latent moments cannot distinguish this distribution from a Gaussian with the same latent moments, up to the clipping, and is therefore misspecified. } 

\subsection{\RC{Hyperparameter ablation}}
\label{appx:hyper_ablation}

\RC{Recall that $\gamma$ is the objective multiplier of~\eqref{eq:bilevel_sampled}, weighting the average case against the $\CVaR$ term, and $\eta$ is the risk level of the $\CVaR$
constraint. 
Each of the first six rows of Tables~\ref{table:port_abl_lin}--\ref{table:inv_abl_nl} restricts the \gls{LRO} hyperparameter selection to a
single $(\eta,\gamma)$. 
The three rows marked $\dagger$ are the settings reported in Section~\ref{sec:examples}. 
The rows after the two summary rows keep these three settings and change only the shape learning-rate multiplier $c$, with $c=0$ the center-only variant.
The tune column follows Table~\ref{table1}, where the base number is the average time for one setting, including its $\rho$ search, and the multiplier is the number of settings.}

\RC{Across both problems, no single $(\eta,\gamma)$ is best everywhere, which is why we tune over a small set.
Increasing $\gamma$ can improve the average case (see the portfolio example), and sometimes even without trading off tail performance, though the effect is not guaranteed. In fact, in the inventory example, increasing $\gamma$ led to worse average-case performance, which shows that increasing the penalty does not guarantee a better solution. The dependency on $\eta$ is also not straightforward.
Although $\rho$ is recalibrated to the same risk target afterwards, the realized $\hat{\eta}$ still varies with the training risk level $\eta$, due to the noisy validation-dependent selection.
Larger $\eta$ resulted in an improved average case in the portfolio example, though again the inventory example exhibited the opposite trend.}

\RC{In general, how much the hyperparameters matter is highly problem-dependent and model-dependent.
On the inventory example, for instance, the $\CVaR$ varies almost three times as much across the six settings for the MLP class as for the linear class.
However, every single setting outperforms the MSE baseline in both tail metrics, and in general also in the average case, except for the inventory MLP class (where 2 of the 6 settings improve on it).
Tuning over all six settings does not meaningfully improve on the three we report, so the restricted search buys a fair comparison with the baselines at little or no cost.
How much $c$ matters is also problem-dependent.
On the inventory example with the linear class, the $\CVaR$ gain over MSE shrinks from about $37\%$ at the selected $c=1$ to about $20\%$ at $c=0.1$ and $9\%$ at $c=0$, a far larger effect than that of $(\eta,\gamma)$.
On the portfolio example, $c$ moves the average case less than $(\eta,\gamma)$ does.
Overall, a limited tuning budget is best spent on the model class and $c$ before $(\eta,\gamma)$, since their effects can be far larger, and every configuration in these tables already improves on MSE in both tail metrics.}

\begin{table}[htp]
\RC{\centering
  \fontsize{10}{12}\selectfont\setlength{\tabcolsep}{4pt}
 \adjustbox{max width=\linewidth}{
	\begin{tabular}{llllllll}\toprule
		$\eta$ & $\gamma$ & $\hat{\eta}$ & infeas & $\VaR_{90}$ & $\CVaR_{90}$ & average-case & tune (s)\\
        \midrule
        $0.10$ & $0.2$$^\dagger$ & $0.1190$ & $0.1146$ & $0.0621 \pm 0.0047$ & $0.0996 \pm 0.0065$ & $-0.0781 \pm 0.0037$ & $1\cdot 184$\\
        $0.10$ & $1$$^\dagger$ & $0.1230$ & $0.1176$ & $\bf{0.0571 \pm 0.0073}$ & $0.0962 \pm 0.0096$ & $-0.0884 \pm 0.0065$ & $1\cdot 169$\\
        $0.10$ & $5$$^\dagger$ & $0.1188$ & $0.1136$ & $0.0614 \pm 0.0078$ & $0.1006 \pm 0.0107$ & $\bf{-0.0948 \pm 0.0074}$ & $1\cdot 166$\\
        $0.15$ & $0.2$ & $0.1240$ & $0.1194$ & $0.0572 \pm 0.0094$ & $\bf{0.0955 \pm 0.0116}$ & $-0.0849 \pm 0.0060$ & $1\cdot 150$\\
        $0.15$ & $1$ & $0.1172$ & $0.1128$ & $0.0574 \pm 0.0063$ & $0.0965 \pm 0.0054$ & $-0.0917 \pm 0.0041$ & $1\cdot 149$\\
        $0.15$ & $5$ & $0.1154$ & $0.1110$ & $0.0619 \pm 0.0079$ & $0.1036 \pm 0.0081$ & $\bf{-0.0948 \pm 0.0019}$ & $1\cdot 135$\\
        \midrule
        \multicolumn{2}{l}{all six, $c=1$} & $0.1210$ & $0.1152$ & $0.0576 \pm 0.0081$ & $0.0962 \pm 0.0098$ & $-0.0891 \pm 0.0079$ & $6\cdot 159$\\
        \multicolumn{2}{l}{reported ($\dagger$), $c=1$} & $0.1230$ & $0.1176$ & $0.0569 \pm 0.0091$ & $0.0943 \pm 0.0077$ & $-0.0897 \pm 0.0107$ & $3\cdot 173$\\
        \midrule
        \multicolumn{2}{l}{$c=0$} & $0.1104$ & $0.1042$ & $0.0557 \pm 0.0089$ & $0.0970 \pm 0.0072$ & $-0.0884 \pm 0.0054$ & $3\cdot 153$\\
        \multicolumn{2}{l}{$c=0.03$} & $0.1072$ & $0.1008$ & $0.0574 \pm 0.0061$ & $0.0985 \pm 0.0067$ & $-0.0876 \pm 0.0048$ & $3\cdot 165$\\
        \multicolumn{2}{l}{$c=0.1$} & $0.1104$ & $0.1054$ & $0.0570 \pm 0.0072$ & $0.0980 \pm 0.0083$ & $-0.0896 \pm 0.0056$ & $3\cdot 170$\\
        \midrule
        \multicolumn{2}{l}{MSE} & $0.1138$ & $0.1054$ & $0.0627 \pm 0.0048$ & $0.1043 \pm 0.0044$ & $-0.0732 \pm 0.0071$ & $6\cdot 40$\\
 \bottomrule
		\end{tabular}}
		\caption{Hyperparameter ablation for \gls{LRO} with the linear head ($c=1$) on the portfolio example at $n=10$, $N=300$, over 10 runs. Values are mean $\pm$ half-IQR. Each of the first six rows restricts the selection to that single $(\eta,\gamma)$; the next two rows tune over all six settings and over the three marked $\dagger$, which is what Table~\ref{table:port_n10} reports. The rows after these change only the shape learning-rate multiplier $c$, with $c=0$ the center-only variant, and tune over the same three settings marked $\dagger$. The last row is the MSE set with the same head, as reported in Table~\ref{table:port_n10}. The best value for each metric among the six single settings is in bold.}
		\label{table:port_abl_lin}
}\end{table}

\begin{table}[htp]
\RC{\centering
  \fontsize{10}{12}\selectfont\setlength{\tabcolsep}{4pt}
 \adjustbox{max width=\linewidth}{
	\begin{tabular}{llllllll}\toprule
		$\eta$ & $\gamma$ & $\hat{\eta}$ & infeas & $\VaR_{90}$ & $\CVaR_{90}$ & average-case & tune (s)\\
        \midrule
        $0.10$ & $0.2$$^\dagger$ & $0.1250$ & $0.1208$ & $0.0654 \pm 0.0097$ & $0.1053 \pm 0.0139$ & $-0.0746 \pm 0.0073$ & $1\cdot 355$\\
        $0.10$ & $1$$^\dagger$ & $0.1262$ & $0.1232$ & $0.0660 \pm 0.0073$ & $0.1058 \pm 0.0050$ & $-0.0827 \pm 0.0041$ & $1\cdot 376$\\
        $0.10$ & $5$$^\dagger$ & $0.1384$ & $0.1326$ & $0.0646 \pm 0.0063$ & $0.1049 \pm 0.0083$ & $-0.0912 \pm 0.0040$ & $1\cdot 233$\\
        $0.15$ & $0.2$ & $0.1238$ & $0.1206$ & $0.0635 \pm 0.0064$ & $\bf{0.1019 \pm 0.0074}$ & $-0.0792 \pm 0.0055$ & $1\cdot 472$\\
        $0.15$ & $1$ & $0.1294$ & $0.1252$ & $0.0670 \pm 0.0039$ & $0.1079 \pm 0.0070$ & $-0.0846 \pm 0.0047$ & $1\cdot 316$\\
        $0.15$ & $5$ & $0.1352$ & $0.1286$ & $\bf{0.0606 \pm 0.0070}$ & $0.1033 \pm 0.0074$ & $\bf{-0.0953 \pm 0.0036}$ & $1\cdot 234$\\
        \midrule
        \multicolumn{2}{l}{all six, $c=1$} & $0.1288$ & $0.1232$ & $0.0620 \pm 0.0077$ & $0.1016 \pm 0.0127$ & $-0.0882 \pm 0.0077$ & $6\cdot 331$\\
        \multicolumn{2}{l}{reported ($\dagger$), $c=1$} & $0.1266$ & $0.1222$ & $0.0622 \pm 0.0069$ & $0.1019 \pm 0.0089$ & $-0.0858 \pm 0.0055$ & $3\cdot 321$\\
        \midrule
        \multicolumn{2}{l}{$c=0$} & $0.1216$ & $0.1106$ & $0.0631 \pm 0.0073$ & $0.1043 \pm 0.0065$ & $-0.0818 \pm 0.0108$ & $3\cdot 278$\\
        \multicolumn{2}{l}{$c=0.03$} & $0.1224$ & $0.1156$ & $0.0618 \pm 0.0091$ & $0.1020 \pm 0.0082$ & $-0.0805 \pm 0.0082$ & $3\cdot 255$\\
        \multicolumn{2}{l}{$c=0.1$} & $0.1226$ & $0.1156$ & $0.0629 \pm 0.0063$ & $0.1042 \pm 0.0050$ & $-0.0833 \pm 0.0076$ & $3\cdot 279$\\
        \midrule
        \multicolumn{2}{l}{MSE} & $0.1226$ & $0.1102$ & $0.0688 \pm 0.0094$ & $0.1102 \pm 0.0097$ & $-0.0690 \pm 0.0067$ & $6\cdot 38$\\
 \bottomrule
		\end{tabular}}
		\caption{Hyperparameter ablation for \gls{LRO} with the MLP head ($c=1$) on the portfolio example at $n=10$, $N=300$. Columns and conventions as in Table~\ref{table:port_abl_lin}.}
		\label{table:port_abl_nl}
}\end{table}

\begin{table}[htp]
\RC{\centering
  \fontsize{10}{12}\selectfont\setlength{\tabcolsep}{4pt}
 \adjustbox{max width=\linewidth}{
	\begin{tabular}{llllllll}\toprule
		$\eta$ & $\gamma$ & $\hat{\eta}$ & infeas & $\VaR_{90}$ & $\CVaR_{90}$ & average-case & tune (s)\\
        \midrule
        $0.10$ & $0.2$$^\dagger$ & $0.0970$ & $0.0970$ & $-301.31 \pm 5.52$ & $-251.41 \pm 8.75$ & $\bf{-430.02 \pm 4.82}$ & $1\cdot 968$\\
        $0.10$ & $1$$^\dagger$ & $0.0890$ & $0.0890$ & $-295.17 \pm 7.26$ & $-242.68 \pm 11.07$ & $-427.76 \pm 5.59$ & $1\cdot 938$\\
        $0.10$ & $5$$^\dagger$ & $0.0944$ & $0.0944$ & $-298.60 \pm 17.65$ & $-247.94 \pm 13.72$ & $-425.21 \pm 10.02$ & $1\cdot 901$\\
        $0.15$ & $0.2$ & $0.0998$ & $0.0998$ & $\bf{-305.79 \pm 3.36}$ & $\bf{-254.18 \pm 5.43}$ & $-427.58 \pm 8.28$ & $1\cdot 883$\\
        $0.15$ & $1$ & $0.0956$ & $0.0956$ & $-302.62 \pm 11.41$ & $-250.26 \pm 10.40$ & $-427.10 \pm 9.07$ & $1\cdot 903$\\
        $0.15$ & $5$ & $0.0960$ & $0.0960$ & $-302.16 \pm 10.80$ & $-250.92 \pm 5.83$ & $-419.02 \pm 5.49$ & $1\cdot 886$\\
        \midrule
        \multicolumn{2}{l}{all six, $c=1$} & $0.1038$ & $0.1038$ & $-306.47 \pm 6.71$ & $-255.83 \pm 10.36$ & $-431.52 \pm 10.24$ & $6\cdot 913$\\
        \multicolumn{2}{l}{reported ($\dagger$), $c=1$} & $0.0972$ & $0.0972$ & $-302.10 \pm 6.32$ & $-251.26 \pm 10.34$ & $-430.78 \pm 6.29$ & $3\cdot 936$\\
        \midrule
        \multicolumn{2}{l}{$c=0$} & $0.0906$ & $0.0906$ & $-253.69 \pm 11.61$ & $-199.36 \pm 12.46$ & $-423.54 \pm 4.81$ & $3\cdot 886$\\
        \multicolumn{2}{l}{$c=0.03$} & $0.0870$ & $0.0870$ & $-258.49 \pm 11.20$ & $-203.64 \pm 13.69$ & $-424.30 \pm 2.09$ & $3\cdot 853$\\
        \multicolumn{2}{l}{$c=0.1$} & $0.0906$ & $0.0906$ & $-274.03 \pm 13.55$ & $-219.88 \pm 17.21$ & $-428.28 \pm 2.80$ & $3\cdot 853$\\
        \midrule
        \multicolumn{2}{l}{MSE} & $0.0938$ & $0.0938$ & $-236.55 \pm 6.37$ & $-183.21 \pm 6.29$ & $-417.52 \pm 1.79$ & $6\cdot 85$\\
 \bottomrule
		\end{tabular}}
		\caption{Hyperparameter ablation for \gls{LRO} with the linear head ($c=1$) on the inventory example, over 10 runs. Columns and conventions as in Table~\ref{table:port_abl_lin}.}
		\label{table:inv_abl_lin}
}\end{table}

\begin{table}[htp]
\RC{\centering
  \fontsize{10}{12}\selectfont\setlength{\tabcolsep}{4pt}
 \adjustbox{max width=\linewidth}{
	\begin{tabular}{llllllll}\toprule
		$\eta$ & $\gamma$ & $\hat{\eta}$ & infeas & $\VaR_{90}$ & $\CVaR_{90}$ & average-case & tune (s)\\
        \midrule
        $0.10$ & $0.2$$^\dagger$ & $0.1006$ & $0.1006$ & $\bf{-309.13 \pm 10.49}$ & $\bf{-241.20 \pm 12.25}$ & $\bf{-453.40 \pm 4.58}$ & $1\cdot 895$\\
        $0.10$ & $1$$^\dagger$ & $0.0812$ & $0.0812$ & $-286.24 \pm 18.01$ & $-217.52 \pm 10.64$ & $-447.35 \pm 3.59$ & $1\cdot 906$\\
        $0.10$ & $5$$^\dagger$ & $0.0938$ & $0.0938$ & $-280.01 \pm 18.65$ & $-211.14 \pm 11.30$ & $-449.25 \pm 5.88$ & $1\cdot 917$\\
        $0.15$ & $0.2$ & $0.0948$ & $0.0948$ & $-306.54 \pm 11.13$ & $-236.16 \pm 11.16$ & $-452.10 \pm 5.08$ & $1\cdot 930$\\
        $0.15$ & $1$ & $0.0850$ & $0.0850$ & $-284.96 \pm 17.11$ & $-215.66 \pm 17.72$ & $-445.99 \pm 6.75$ & $1\cdot 924$\\
        $0.15$ & $5$ & $0.0954$ & $0.0954$ & $-281.14 \pm 8.22$ & $-209.61 \pm 7.99$ & $-444.91 \pm 5.61$ & $1\cdot 942$\\
        \midrule
        \multicolumn{2}{l}{all six, $c=1$} & $0.1004$ & $0.1004$ & $-317.51 \pm 10.13$ & $-245.67 \pm 10.47$ & $-452.80 \pm 4.51$ & $6\cdot 919$\\
        \multicolumn{2}{l}{reported ($\dagger$), $c=1$} & $0.0992$ & $0.0992$ & $-308.18 \pm 5.69$ & $-239.79 \pm 11.21$ & $-452.80 \pm 5.65$ & $3\cdot 906$\\
        \midrule
        \multicolumn{2}{l}{$c=0$} & $0.1032$ & $0.1032$ & $-312.76 \pm 6.42$ & $-246.60 \pm 5.98$ & $-448.96 \pm 4.02$ & $3\cdot 924$\\
        \multicolumn{2}{l}{$c=0.03$} & $0.1080$ & $0.1080$ & $-313.77 \pm 7.68$ & $-248.21 \pm 10.36$ & $-449.84 \pm 4.55$ & $3\cdot 951$\\
        \multicolumn{2}{l}{$c=0.1$} & $0.1048$ & $0.1048$ & $-312.77 \pm 9.63$ & $-247.24 \pm 11.70$ & $-450.45 \pm 4.18$ & $3\cdot 893$\\
        \midrule
        \multicolumn{2}{l}{MSE} & $0.0850$ & $0.0850$ & $-264.33 \pm 5.64$ & $-200.55 \pm 4.27$ & $-449.49 \pm 1.41$ & $6\cdot 88$\\
 \bottomrule
		\end{tabular}}
		\caption{Hyperparameter ablation for \gls{LRO} with the MLP head ($c=1$) on the inventory example. Columns and conventions as in Table~\ref{table:port_abl_lin}.}
		\label{table:inv_abl_nl}
}\end{table}

\subsection{\RC{DFL baseline}}
\label{appx:dfl}

\RC{Tables~\ref{table:port_dfl} and~\ref{table:inv_dfl} report the non-robust DFL baseline of Section~\ref{sec:intro_sets} with only the average cost penalized, for the linear and MLP classes.
It uses the same data, training settings, and model classes as \gls{LRO}, with $\rho$ fixed at $10^{-6}$.
Since DFL has no set to size, its violation probability cannot be calibrated to the target.
On the portfolio example, its $\hat{\eta}$ is about $0.5$ and its infeasibility about $0.2$, twice that of the methods in Tables~\ref{table:port_n10} and~\ref{table:port_n30}.
Its objective values are computed over the feasible instances only, so they are not comparable to those of the other methods.
On the inventory example, no test-set instance is feasible.}

\begin{table}[ht]
\RC{\centering
  \fontsize{10}{12}\selectfont\setlength{\tabcolsep}{4pt}
 \adjustbox{max width=\linewidth}{
	\begin{tabular}{llllllll}   \toprule
		Method & $\hat{\eta}$ & infeas & $\VaR_{90}$ & $\CVaR_{90}$ & average-case & tune (s) & inf. (ms)\\
        \midrule
        \multicolumn{8}{l}{$n=10$, $N=300$}\\
        DFL (linear)  & $0.5282$ & $0.2144$ & $0.0041 \pm 0.0042$ & $0.0250 \pm 0.0035$ & $-0.1742 \pm 0.0048$ & $1\cdot 104$ & $16.1$\\
        DFL (MLP)  & $0.5360$ & $0.2124$ & $0.0052 \pm 0.0038$ & $0.0237 \pm 0.0025$ & $-0.1648 \pm 0.0074$ & $1\cdot 103$ & $17.3$\\
        \midrule
        \multicolumn{8}{l}{$n=10$, $N=600$}\\
        DFL (linear)  & $0.5173$ & $0.2157$ & $0.0019 \pm 0.0036$ & $0.0235 \pm 0.0023$ & $-0.1731 \pm 0.0067$ & $1\cdot 81$ & $14.3$\\
        DFL (MLP)  & $0.5223$ & $0.2133$ & $0.0045 \pm 0.0040$ & $0.0250 \pm 0.0025$ & $-0.1615 \pm 0.0055$ & $1\cdot 80$ & $14.5$\\
        \midrule
        \multicolumn{8}{l}{$n=30$, $N=300$}\\
        DFL (linear)  & $0.4992$ & $0.1838$ & $-0.0252 \pm 0.0035$ & $-0.0039 \pm 0.0027$ & $-0.2137 \pm 0.0059$ & $1\cdot 498$ & $78.3$\\
        DFL (MLP)  & $0.5356$ & $0.2094$ & $-0.0197 \pm 0.0026$ & $-0.0007 \pm 0.0022$ & $-0.1969 \pm 0.0045$ & $1\cdot 494$ & $76.6$\\
        \midrule
        \multicolumn{8}{l}{$n=30$, $N=600$}\\
        DFL (linear)  & $0.4949$ & $0.1905$ & $-0.0258 \pm 0.0019$ & $-0.0025 \pm 0.0014$ & $-0.2133 \pm 0.0041$ & $1\cdot 437$ & $70.1$\\
        DFL (MLP)  & $0.5221$ & $0.2121$ & $-0.0181 \pm 0.0011$ & $0.0024 \pm 0.0015$ & $-0.1939 \pm 0.0048$ & $1\cdot 444$ & $73.0$\\
		\bottomrule
		\end{tabular}}
		\caption{\RC{Non-robust DFL on the portfolio example, averaged over 10 runs. Columns and conventions as in Table~\ref{table:port_n10}. $\VaR_{90}$, $\CVaR_{90}$, and average-case are computed over the feasible test-set instances only.}}
		\label{table:port_dfl}
}\end{table}

\begin{table}[ht]
\RC{\centering
  \fontsize{10}{12}\selectfont\setlength{\tabcolsep}{4pt}
 \adjustbox{max width=\linewidth}{
	\begin{tabular}{llllllll}   \toprule
		Method & $\hat{\eta}$ & infeas & $\VaR_{90}$ & $\CVaR_{90}$ & average-case & tune (s) & inf. (ms)\\
        \midrule
        \multicolumn{8}{l}{$n=10$, $N=300$}\\
        DFL (linear)  & $1.0000$ & $1.0000$ & -- & -- & -- & $1\cdot 655$ & $33.2$\\
        DFL (MLP)  & $1.0000$ & $1.0000$ & -- & -- & -- & $1\cdot 540$ & $33.7$\\
		\bottomrule
		\end{tabular}}
		\caption{\RC{Non-robust DFL on the inventory example, averaged over 10 runs. Columns and conventions as in Table~\ref{table1}. No test-set instance is feasible, so the objective columns are empty.}}
		\label{table:inv_dfl}
}\end{table}
 \fi

\FloatBarrier
\end{appendices}


\begin{thebibliography}{}

\bibitem[Adelh{\"u}tte et~al., 2023]{pareto3}
Adelh{\"u}tte, D., Biefel, C., Kuchlbauer, M., and Rolfes, J. (2023).
\newblock Pareto robust optimization on {E}uclidean vector spaces.
\newblock {\em Optimization Letters}, 17(3):771--788.

\bibitem[Agrawal et~al., 2019]{cvxpylayers2019}
Agrawal, A., Amos, B., Barratt, S., Boyd, S., Diamond, S., and Kolter, J.~Z. (2019).
\newblock Differentiable convex optimization layers.
\newblock In {\em Advances in Neural Information Processing Systems}, volume~32.

\bibitem[Amos and Kolter, 2017]{amos2021optnet}
Amos, B. and Kolter, J.~Z. (2017).
\newblock {OptNet}: Differentiable optimization as a layer in neural networks.
\newblock In {\em Proceedings of the 34th International Conference on Machine Learning}, volume~70 of {\em Proceedings of Machine Learning Research}, pages 136--145. PMLR.

\bibitem[Angelopoulos and Bates, 2023]{conformal}
Angelopoulos, A.~N. and Bates, S. (2023).
\newblock Conformal prediction: A gentle introduction.
\newblock {\em Foundations and Trends in Machine Learning}, 16(4):494--591.

\bibitem[Bandi and Bertsimas, 2012]{bandi_tractable_2012}
Bandi, C. and Bertsimas, D. (2012).
\newblock Tractable stochastic analysis in high dimensions via robust optimization.
\newblock {\em Mathematical Programming}, 134(1):23--70.

\bibitem[Bauschke and Combettes, 2017]{composites}
Bauschke, H.~H. and Combettes, P.~L. (2017).
\newblock {\em Convex Analysis and Monotone Operator Theory in Hilbert Spaces}.
\newblock Springer, Cham, second edition.

\bibitem[Ben-Tal et~al., 2009]{ben-tal_robust_2009}
Ben-Tal, A., El~Ghaoui, L., and Nemirovski, A. (2009).
\newblock {\em Robust {Optimization}}.
\newblock Princeton University Press.

\bibitem[Ben-Tal and Nemirovski, 2000]{ben-tal_robust_2000}
Ben-Tal, A. and Nemirovski, A. (2000).
\newblock Robust solutions of {Linear} {Programming} problems contaminated with uncertain data.
\newblock {\em Mathematical Programming}, 88(3):411--424.

\bibitem[Ben-Tal and Nemirovski, 2008]{robustconvexopt}
Ben-Tal, A. and Nemirovski, A. (2008).
\newblock Selected topics in robust convex optimization.
\newblock {\em Math. Program.}, 112:125--158.

\bibitem[Bertsimas et~al., 2011]{bertsimassurvey}
Bertsimas, D., Brown, D.~B., and Caramanis, C. (2011).
\newblock Theory and applications of robust optimization.
\newblock {\em SIAM Review}, 53(3):464--501.

\bibitem[Bertsimas and den Hertog, 2022]{robustadaptopt}
Bertsimas, D. and den Hertog, D. (2022).
\newblock {\em Robust and Adaptive Optimization}.
\newblock Dynamic Ideas.

\bibitem[Bertsimas et~al., 2021]{bertsimas_probabilistic_2019}
Bertsimas, D., den Hertog, D., and Pauphilet, J. (2021).
\newblock Probabilistic {Guarantees} in {Robust} {Optimization}.
\newblock {\em SIAM Journal on Optimization}, 31(4):2893--2920.

\bibitem[Bertsimas et~al., 2018]{bertsimas_data-driven_2018}
Bertsimas, D., Gupta, V., and Kallus, N. (2018).
\newblock Data-driven robust optimization.
\newblock {\em Mathematical Programming}, 167(2):235--292.

\bibitem[Bertsimas and Kallus, 2020]{bertsimas2020predictive}
Bertsimas, D. and Kallus, N. (2020).
\newblock From predictive to prescriptive analytics.
\newblock {\em Management Science}, 66(3):1025--1044.

\bibitem[Bertsimas et~al., 2022]{wassinf}
Bertsimas, D., Shtern, S., and Sturt, B. (2022).
\newblock Technical note---two-stage sample robust optimization.
\newblock {\em Operations Research}, 70(1):624--640.

\bibitem[Bertsimas and Sim, 2004]{bertsimas_price_2004}
Bertsimas, D. and Sim, M. (2004).
\newblock The {Price} of {Robustness}.
\newblock {\em Operations Research}, 52(1):35--53.

\bibitem[Bertsimas et~al., 2024]{pareto1}
Bertsimas, D., ten Eikelder, S. C.~M., den Hertog, D., and Trichakis, N. (2024).
\newblock Pareto adaptive robust optimality via a {F}ourier--{M}otzkin elimination lens.
\newblock {\em Mathematical Programming}, 205(1--2):485--538.

\bibitem[Birgin and Mart{\'\i}nez, 2014]{auglag_safe}
Birgin, E.~G. and Mart{\'\i}nez, J.~M. (2014).
\newblock {\em Practical Augmented Lagrangian Methods for Constrained Optimization}, volume~10 of {\em Fundamentals of Algorithms}.
\newblock SIAM, Philadelphia, PA.

\bibitem[Boero et~al., 2026]{boero2025alcole}
Boero, I., Hounie, I., and Ribeiro, A. (2026).
\newblock {AL-CoLe}: Augmented {L}agrangian for constrained learning.
\newblock In {\em IEEE International Conference on Acoustics, Speech and Signal Processing (ICASSP)}, pages 1526--1530.

\bibitem[Bolte et~al., 2021]{path-diff}
Bolte, J., Le, T., Pauwels, E., and Silveti-Falls, A. (2021).
\newblock Nonsmooth implicit differentiation for machine-learning and optimization.
\newblock In {\em Advances in Neural Information Processing Systems}, volume~34, pages 13537--13549.

\bibitem[Bolte and Pauwels, 2021]{conservative}
Bolte, J. and Pauwels, E. (2021).
\newblock Conservative set valued fields, automatic differentiation, stochastic gradient methods and deep learning.
\newblock {\em Mathematical Programming}, 188(1):19--51.

\bibitem[Bonnans and Ram\'{\i}rez~C., 2005]{bonnans2005}
Bonnans, J.~F. and Ram\'{\i}rez~C., H. (2005).
\newblock Perturbation analysis of second-order cone programming problems.
\newblock {\em Mathematical Programming}, 104(2--3):205--227.

\bibitem[Bonnans and Shapiro, 2000]{bonnans2000}
Bonnans, J.~F. and Shapiro, A. (2000).
\newblock {\em Perturbation Analysis of Optimization Problems}.
\newblock Springer Series in Operations Research. Springer, New York.

\bibitem[Busseti et~al., 2019]{cones}
Busseti, E., Moursi, W., and Boyd, S. (2019).
\newblock Solution refinement at regular points of conic problems.
\newblock {\em Computational Optimization and Applications}, 74(3):627--643.

\bibitem[Calafiore and Campi, 2005]{scenario}
Calafiore, G. and Campi, M.~C. (2005).
\newblock Uncertain convex programs: randomized solutions and confidence levels.
\newblock {\em Mathematical Programming}, 102(1):25--46.

\bibitem[Cameron et~al., 2022]{end-to-end}
Cameron, C., Hartford, J., Lundy, T., and Leyton-Brown, K. (2022).
\newblock The perils of learning before optimizing.
\newblock {\em Proceedings of the AAAI Conference on Artificial Intelligence}, 36:3708--3715.

\bibitem[Chenreddy and Delage, 2024]{cond-cov}
Chenreddy, A. and Delage, E. (2024).
\newblock End-to-end conditional robust optimization.
\newblock In {\em Proceedings of the Fortieth Conference on Uncertainty in Artificial Intelligence}, volume 244 of {\em Proceedings of Machine Learning Research}, pages 736--748. PMLR.

\bibitem[Clarke, 1975]{Clarke}
Clarke, F.~H. (1975).
\newblock Generalized gradients and applications.
\newblock {\em Transactions of the American Mathematical Society}, 205:247--262.

\bibitem[Clarke, 1990]{clarke1990}
Clarke, F.~H. (1990).
\newblock {\em Optimization and Nonsmooth Analysis}.
\newblock Classics in Applied Mathematics. SIAM, Philadelphia.

\bibitem[Coste, 2000]{definable2}
Coste, M. (2000).
\newblock {\em An Introduction to O-minimal Geometry}.
\newblock Dottorato di ricerca in matematica / Universit{\`a} di Pisa, Dipartimento di Matematica. Istituti editoriali e poligrafici internazionali.

\bibitem[Davis et~al., 2020]{Davis}
Davis, D., Drusvyatskiy, D., Kakade, S.~M., and Lee, J.~D. (2020).
\newblock Stochastic subgradient method converges on tame functions.
\newblock {\em Foundations of Computational Mathematics}, 20(1):119--154.

\bibitem[Demirovi{\'c} et~al., 2019]{end1}
Demirovi{\'c}, E., Stuckey, P.~J., Bailey, J., Chan, J., Leckie, C., Ramamohanarao, K., and Guns, T. (2019).
\newblock An investigation into prediction + optimisation for the knapsack problem.
\newblock In {\em Integration of Constraint Programming, Artificial Intelligence, and Operations Research}, Lecture Notes in Computer Science, pages 241--257. Springer.

\bibitem[Donti et~al., 2017]{amos1}
Donti, P.~L., Amos, B., and Kolter, J.~Z. (2017).
\newblock Task-based end-to-end model learning in stochastic optimization.
\newblock In {\em Advances in Neural Information Processing Systems}, volume~30.

\bibitem[Elmachtoub and Grigas, 2022]{predict_opt}
Elmachtoub, A.~N. and Grigas, P. (2022).
\newblock Smart ``{P}redict, {T}hen {O}ptimize''.
\newblock {\em Manage. Sci.}, 68(1):9--26.

\bibitem[Elmachtoub et~al., 2025]{elmachtoub2025eto}
Elmachtoub, A.~N., Lam, H., Zhang, H., and Zhao, Y. (2025).
\newblock Estimate-then-optimize versus integrated-estimation-optimization versus sample average approximation: A stochastic dominance perspective.
\newblock {\em arXiv preprint arXiv:2304.06833v4}.

\bibitem[Facchinei and Pang, 2003]{facchinei2003}
Facchinei, F. and Pang, J.-S. (2003).
\newblock {\em Finite-Dimensional Variational Inequalities and Complementarity Problems}.
\newblock Springer, New York.

\bibitem[Fioretto et~al., 2021]{fioretto2020lagrangian}
Fioretto, F., Van~Hentenryck, P., Mak, T. W.~K., Tran, C., Baldo, F., and Lombardi, M. (2021).
\newblock Lagrangian duality for constrained deep learning.
\newblock In {\em Machine Learning and Knowledge Discovery in Databases: Applied Data Science and Demo Track (ECML PKDD 2020)}, volume 12461 of {\em Lecture Notes in Computer Science}, pages 118--135, Cham. Springer.

\bibitem[Gao, 2023]{DBLP:gao2020}
Gao, R. (2023).
\newblock Finite-sample guarantees for {W}asserstein distributionally robust optimization: Breaking the curse of dimensionality.
\newblock {\em Operations Research}, 71(6):2291--2306.

\bibitem[Gao and Kleywegt, 2023]{gao2016distributionally}
Gao, R. and Kleywegt, A. (2023).
\newblock Distributionally robust stochastic optimization with {W}asserstein distance.
\newblock {\em Mathematics of Operations Research}, 48(2):603--655.

\bibitem[Goerigk and Kurtz, 2023]{GOERIGK2023106087}
Goerigk, M. and Kurtz, J. (2023).
\newblock Data-driven robust optimization using deep neural networks.
\newblock {\em Computers \& Operations Research}, 151:106087.

\bibitem[Gupta, 2019]{gupta}
Gupta, V. (2019).
\newblock Near-optimal {B}ayesian ambiguity sets for distributionally robust optimization.
\newblock {\em Management Science}, 65(9):4242--4260.

\bibitem[Guzman et~al., 2017]{guarantees1}
Guzman, Y.~A., Matthews, L.~R., and Floudas, C.~A. (2017).
\newblock New a priori and a posteriori probabilistic bounds for robust counterpart optimization: {III}. {E}xact and near-exact a posteriori expressions for known probability distributions.
\newblock {\em Computers \& Chemical Engineering}, 103:116--143.

\bibitem[{Homem-de-Mello} et~al., 2026]{homemdemello2026}
{Homem-de-Mello}, T., Valencia, J., Lagos, F., and Lagos, G. (2026).
\newblock Forecasting outside the box: application-driven optimal pointwise forecasts for a class of two-stage stochastic programs.
\newblock {\em Mathematical Programming}.

\bibitem[Hong et~al., 2021]{hong2022ms}
Hong, L.~J., Huang, Z., and Lam, H. (2021).
\newblock Learning-based robust optimization: Procedures and statistical guarantees.
\newblock {\em Management Science}, 67(6):3447--3467.

\bibitem[Hu et~al., 2022]{hu2022fast}
Hu, Y., Kallus, N., and Mao, X. (2022).
\newblock Fast rates for contextual linear optimization.
\newblock {\em Management Science}, 68(6):4236--4245.

\bibitem[Iancu and Trichakis, 2014]{pareto}
Iancu, D.~A. and Trichakis, N. (2014).
\newblock Pareto efficiency in robust optimization.
\newblock {\em Management Science}, 60(1):130--147.

\bibitem[Kanzow et~al., 2009]{cones1}
Kanzow, C., Ferenczi, I., and Fukushima, M. (2009).
\newblock On the local convergence of semismooth {N}ewton methods for linear and nonlinear second-order cone programs without strict complementarity.
\newblock {\em SIAM Journal on Optimization}, 20(1):297--320.

\bibitem[Kingma and Ba, 2015]{kingma2015adam}
Kingma, D.~P. and Ba, J. (2015).
\newblock {Adam}: A method for stochastic optimization.
\newblock In {\em Proceedings of the 3rd International Conference on Learning Representations}.

\bibitem[Kong et~al., 2011]{kong2011jacobian}
Kong, L., Tun{\c{c}}el, L., and Xiu, N. (2011).
\newblock Equivalent conditions for {J}acobian nonsingularity in linear symmetric cone programming.
\newblock {\em Journal of Optimization Theory and Applications}, 148(2):364--389.

\bibitem[Kotary and Fioretto, 2024]{kotary2024deepal}
Kotary, J. and Fioretto, F. (2024).
\newblock Learning constrained optimization with deep augmented {L}agrangian methods.
\newblock {\em arXiv preprint arXiv:2403.03454}.

\bibitem[Kotary et~al., 2022]{kotary2022fast}
Kotary, J., Fioretto, F., and Van~Hentenryck, P. (2022).
\newblock Fast approximations for job shop scheduling: A {L}agrangian dual deep learning method.
\newblock {\em Proceedings of the AAAI Conference on Artificial Intelligence}, 36(7):7239--7246.

\bibitem[Kotary et~al., 2021]{end-survey}
Kotary, J., Fioretto, F., Van~Hentenryck, P., and Wilder, B. (2021).
\newblock End-to-end constrained optimization learning: A survey.
\newblock In {\em Proceedings of the Thirtieth International Joint Conference on Artificial Intelligence (IJCAI-21)}, pages 4475--4482.

\bibitem[Li et~al., 2012]{guarantees_2012}
Li, Z., Tang, Q., and Floudas, C.~A. (2012).
\newblock A comparative theoretical and computational study on robust counterpart optimization: {II}. {P}robabilistic guarantees on constraint satisfaction.
\newblock {\em Industrial \& Engineering Chemistry Research}, 51(19):6769--6788.

\bibitem[Long et~al., 2023]{long2023robust}
Long, D., Sim, M., and Zhou, M. (2023).
\newblock Robust satisficing.
\newblock {\em Operations Research}, 71(1):61--82.

\bibitem[Mandi et~al., 2024]{dfl}
Mandi, J., Kotary, J., Berden, S., Mulamba, M., Bucarey, V., Guns, T., and Fioretto, F. (2024).
\newblock Decision-focused learning: Foundations, state of the art, benchmark and future opportunities.
\newblock {\em Journal of Artificial Intelligence Research}, 80:1623--1701.

\bibitem[Messoudi et~al., 2022]{knn}
Messoudi, S., Destercke, S., and Rousseau, S. (2022).
\newblock Ellipsoidal conformal inference for multi-target regression.
\newblock In {\em Proceedings of the Eleventh Symposium on Conformal and Probabilistic Prediction with Applications}, volume 179 of {\em Proceedings of Machine Learning Research}, pages 294--306. PMLR.

\bibitem[{Mohajerin Esfahani} and Kuhn, 2018]{mohajerin_esfahani_data-driven_2018}
{Mohajerin Esfahani}, P. and Kuhn, D. (2018).
\newblock Data-driven distributionally robust optimization using the {W}asserstein metric: performance guarantees and tractable reformulations.
\newblock {\em Mathematical Programming}, 171:115--166.

\bibitem[Pataki and Tun{\c{c}}el, 2001]{pataki2001generic}
Pataki, G. and Tun{\c{c}}el, L. (2001).
\newblock On the generic properties of convex optimization problems in conic form.
\newblock {\em Mathematical Programming}, 89(3):449--457.

\bibitem[Per{\v{s}}ak and Anjos, 2023]{persak2023contextual}
Per{\v{s}}ak, E. and Anjos, M.~F. (2023).
\newblock Contextual robust optimisation with uncertainty quantification.
\newblock In {\em Integration of Constraint Programming, Artificial Intelligence, and Operations Research (CPAIOR 2023)}, volume 13884 of {\em Lecture Notes in Computer Science}, pages 124--132, Cham. Springer.

\bibitem[Rockafellar and Uryasev, 2002]{cvar}
Rockafellar, R.~T. and Uryasev, S. (2002).
\newblock Conditional value-at-risk for general loss distributions.
\newblock {\em Journal of Banking \& Finance}, 26(7):1443--1471.

\bibitem[Rockafellar and Wets, 1998]{rockafellar_wets_1998}
Rockafellar, R.~T. and Wets, R. J.~B. (1998).
\newblock {\em Variational Analysis}, volume 317 of {\em Grundlehren der mathematischen Wissenschaften}.
\newblock Springer, Berlin.

\bibitem[Roos and den Hertog, 2020]{reduce_conserve}
Roos, E. and den Hertog, D. (2020).
\newblock Reducing conservatism in robust optimization.
\newblock {\em INFORMS Journal on Computing}, 32(4):1109--1127.

\bibitem[Sadana et~al., 2025]{contex}
Sadana, U., Chenreddy, A., Delage, E., Forel, A., Frejinger, E., and Vidal, T. (2025).
\newblock A survey of contextual optimization methods for decision-making under uncertainty.
\newblock {\em European Journal of Operational Research}, 320(2):271--289.

\bibitem[Sarykalin et~al., 2008]{cvar_var}
Sarykalin, S., Serraino, G., and Uryasev, S. (2008).
\newblock Value-at-risk vs. conditional value-at-risk in risk management and optimization.
\newblock In {\em State-of-the-Art Decision-Making Tools in the Information-Intensive Age}, Tutorials in Operations Research, pages 270--294. INFORMS.

\bibitem[Schutte et~al., 2024]{schutte2024robust}
Schutte, N., Postek, K., and Yorke-Smith, N. (2024).
\newblock Robust losses for decision-focused learning.
\newblock In {\em Proceedings of the Thirty-Third International Joint Conference on Artificial Intelligence (IJCAI)}, pages 4868--4875.

\bibitem[Schutte et~al., 2026]{schutte2025proxies}
Schutte, N., Veviurko, G., Postek, K., and Yorke-Smith, N. (2026).
\newblock Sufficient decision proxies for decision-focused learning.
\newblock In {\em Proceedings of the Thirty-Fifth International Joint Conference on Artificial Intelligence (IJCAI)}, pages 2337--2345.

\bibitem[Schuurmans and Patrinos, 2023]{schuurmans2023distributionally}
Schuurmans, M. and Patrinos, P. (2023).
\newblock Distributionally robust optimization using cost-aware ambiguity sets.
\newblock {\em IEEE Control Systems Letters}, 7:1855--1860.

\bibitem[Singh and P{\'o}czos, 2016]{singh2016knn}
Singh, S. and P{\'o}czos, B. (2016).
\newblock Analysis of $k$-nearest neighbor distances with application to entropy estimation.
\newblock {\em arXiv preprint arXiv:1603.08578}.

\bibitem[Stooke et~al., 2020]{stooke2020responsive}
Stooke, A., Achiam, J., and Abbeel, P. (2020).
\newblock Responsive safety in reinforcement learning by {PID} {L}agrangian methods.
\newblock In {\em Proceedings of the 37th International Conference on Machine Learning}, volume 119 of {\em Proceedings of Machine Learning Research}, pages 9133--9143. PMLR.

\bibitem[Tanoumand et~al., 2023]{tanoumanddata}
Tanoumand, N., Bodur, M., and Naoum-Sawaya, J. (2023).
\newblock Data-driven distributionally robust optimization: Intersecting ambiguity sets, performance analysis and tractability.
\newblock Optimization Online preprint.

\bibitem[Uryasev and Rockafellar, 2001]{cvaropt}
Uryasev, S. and Rockafellar, R.~T. (2001).
\newblock Conditional value-at-risk: optimization approach.
\newblock In {\em Stochastic optimization: algorithms and applications}, pages 411--435. Springer.

\bibitem[van~den Dries and Miller, 1996]{definable}
van~den Dries, L. and Miller, C. (1996).
\newblock Geometric categories and o-minimal structures.
\newblock {\em Duke Mathematical Journal}, 84(2):497--540.

\bibitem[{Van Parys} et~al., 2021]{vanparys2021data}
{Van Parys}, B. P.~G., {Mohajerin Esfahani}, P., and Kuhn, D. (2021).
\newblock From data to decisions: Distributionally robust optimization is optimal.
\newblock {\em Management Science}, 67(6):3387--3402.

\bibitem[Vershynin, 2018]{vershynin2018}
Vershynin, R. (2018).
\newblock {\em High-Dimensional Probability: An Introduction with Applications in Data Science}.
\newblock Cambridge University Press.

\bibitem[Wang et~al., 2025]{mro}
Wang, I., Becker, C., Van~Parys, B., and Stellato, B. (2025).
\newblock Mean robust optimization.
\newblock {\em Mathematical Programming}, 213(1--2):1235--1277.

\bibitem[Weed and Bach, 2019]{wass_rate}
Weed, J. and Bach, F. (2019).
\newblock {Sharp asymptotic and finite-sample rates of convergence of empirical measures in {W}asserstein distance}.
\newblock {\em Bernoulli}, 25(4A):2620--2648.

\bibitem[Wilder et~al., 2019a]{end2}
Wilder, B., Dilkina, B., and Tambe, M. (2019a).
\newblock Melding the data-decisions pipeline: Decision-focused learning for combinatorial optimization.
\newblock {\em Proceedings of the AAAI Conference on Artificial Intelligence}, 33(01):1658--1665.

\bibitem[Wilder et~al., 2019b]{end3}
Wilder, B., Ewing, E., Dilkina, B., and Tambe, M. (2019b).
\newblock End to end learning and optimization on graphs.
\newblock In {\em Advances in Neural Information Processing Systems}, volume~32.

\bibitem[Yeh et~al., 2025]{yeh}
Yeh, C., Christianson, N., Wu, A., Wierman, A., and Yue, Y. (2025).
\newblock End-to-end conformal calibration for optimization under uncertainty.
\newblock {\em Transactions on Machine Learning Research}.

\bibitem[Zhen et~al., 2025]{zhen2021mathematical}
Zhen, J., Kuhn, D., and Wiesemann, W. (2025).
\newblock A unified theory of robust and distributionally robust optimization via the primal-worst-equals-dual-best principle.
\newblock {\em Operations Research}, 73(2):862--878.

\end{thebibliography}
\end{document}